\documentclass[reqno,10pt]{amsart}
\usepackage{amssymb}
\usepackage{amsmath}
\usepackage{amsthm}
\usepackage{pgf}
\usepackage{color}
\usepackage{graphicx}   
\usepackage{float}

\usepackage{amsmath}
\usepackage{amssymb}
\usepackage{amsthm}
\usepackage{amsfonts}
\usepackage{graphicx}
\usepackage{ esint }
\usepackage[abs]{overpic}
\usepackage{caption}
\usepackage{wrapfig}
\usepackage{float}
\usepackage{graphicx}

\newtheorem{theorem}{Theorem}[section]
\newtheorem{lemma}[theorem]{Lemma}

\newtheorem{corollary}[theorem]{Corollary}
\newtheorem{definition}[theorem]{Definition}
\newtheorem{example}[theorem]{Example}

\newcommand{\Proof}{\par\noindent{\em Proof. }}
\newcommand{\eop}{\nopagebreak\hspace*{\fill}$\Box$\smallskip}

\numberwithin{equation}{section}

\newtheorem{rem}[theorem]{Remark}

\newcommand{\N}{\Bbb N}
\newcommand{\Z}{\Bbb Z}
\newcommand{\R}{\Bbb R}

\def\sat{{\rm sat}}

\def\eps{\varepsilon}
\def\diam{{\rm diam}}

\def\e{\mathbf{e}}

\def\dist{\operatorname{dist}}

\def\XXint#1#2#3{{\setbox0=\hbox{$#1{#2#3}{\int}$ }
\vcenter{\hbox{$#2#3$ }}\kern-.6\wd0}}

\usepackage{color}

\begin{document}

\title[Decomposition of regular domains into John domains]{On a decomposition of regular domains into John domains with uniform constants}

\author{Manuel Friedrich}
\address[Manuel Friedrich]{Faculty of Mathematics, University of Vienna, Oskar-Morgenstern-Platz 1, A-1090 Vienna, Austria}


\begin{abstract} 
We derive a decomposition result for regular, two-dimensional domains  into John domains with uniform constants. We prove that for every simply connected domain $\Omega \subset \R^2$ with $C^1$-boundary there is a corresponding partition $\Omega = \Omega_1 \cup \ldots \cup \Omega_N$   with $\sum_{j=1}^N \mathcal{H}^1(\partial \Omega_j \setminus \partial \Omega) \le \theta$ such that each component is a John domain with a John constant only depending on $\theta$. The result implies that many inequalities in Sobolev spaces such as Poincar\'e's or Korn's inequality hold on the partition of $\Omega$ for uniform constants, which are independent of $\Omega$.
\end{abstract}


\maketitle
              
                              %
      
      
\pagestyle{myheadings}

\section{Introduction}\label{sec: intro}

It is a fundamental question to identify classes of domains for which the existence of solutions for partial differential equations or the validity of inequalities in Sobolev spaces can be guaranteed. The last decades have witnessed a tremendous process in establishing results for different assumptions on the domains. 

For instance, one of the first proofs of Korn's inequality, being a widely studied inequality due to its importance in the analysis of elasticity equations, was given by Friedrichs \cite{Friedrichs1} for domains allowing for a finite number of corners or edges on the boundary. Subsequently, generalizations  appeared including versions for star-shaped sets \cite{Kondratiev}, general Lipschitz domains \cite{Nit}, and more recently results \cite{Duran1}  were obtained for the broader class of \emph{uniform domains}  using a modification of the extension operator by Jones \cite{Jones}.

On the other hand, it has been known for a long time that many inequalities are false on domains with external cusps.  Several arguments have been provided for this fact (see \cite{Geymonat-Gilardi, Weck}), but the oldest is due to Friedrichs  \cite{Friedrichs2},  who studied an inequality for analytic complex functions (cf. also \cite{Acosta}). 

Recently Acosta, Dur\'an, and Muschietti \cite{Acosta} investigated the existence of solutions of the divergence operator on \emph{John domains} (see \cite{John:1961, Martio-Sarvas:1978, Nakki-Vaisala:1991}). Apart from its application to the study of the Stokes equation the result is of interest due to its connection to Poincar\'e's and Korn's inequality, which may be deduced herefrom. Roughly speaking, a  domain is a John domain if it satisfies a \emph{twisted cone condition} such that each two points can be connected by a curve not  getting too close to the boundary of the domain in terms of a corresponding \emph{John constant} (we refer to Section \ref{sec john-pre}  below for an exact definition).

 John domains represent a very general class allowing for sets with fractal boundary (e.g. Koch's \emph{snowflake}), but at the same time excluding the formation of external cusps. They may be regarded as a very natural and in some sense most general  notion of sets for the investigation of problems alluded to above since  in \cite{Acosta, Buckley}  it has been shown that for domains satisfying the separation property (e.g. for simply connected planar domains) the validity of Poincar\'e's or Korn's inequality implies that the set is a John domain. Moreover, as already observed by Bojarski \cite{Bojarski}, the constant involved in the estimates essentially only depends on the John constant. 

Difficulties concerning the properties and regularity  of domains become even more challenging in models dealing with varying domains, e.g. \emph{free boundary} or  \emph{shape optimization problems}, where the best shape of a set in dependence of a cost functional is identified as the solution of a variational problem (we refer to \cite{Bucur} for an introduction). Another important class is given by  \emph{free discontinuity problems} in the language  of Ambrosio and De Giorgi \cite{DeGiorgi-Ambrosio:1988} with various applications in fields of fracture mechanics or digital image segmentation, where the set of discontinuities of the function of interest is not preassigned, but determined from an energy minimization principle (cf. \cite{Ambrosio-Fusco-Pallara:2000}). 

Obviously without additional conditions there is no hope to derive uniform estimates being independent of the set shape as can be seen, e.g., by considering a sequence of smooth sets converging to a domain with external cusp. Moreover, one may think of Neumann sieve type phenomena (see \cite{Murat}) where the set is only connected by a small periodically distributed contact zone. 

Therefore, many works appeared analyzing the behavior of constants in terms of the domain (cf. \cite{Horgan}  and the references therein) or investigating  special structures as convex, star-shaped   or thin domains (see e.g. \cite{Duran2, Harutyunyan, Lewicka}).  Another approach particularly used in the study of free discontinuity problems is based on the idea to establish results for a certain class of admissible (discontinuity) sets  for which uniform estimates can be shown (we refer e.g. to \cite{Lazzaroni, NegriToader:2013, Rondi}).

Also the present article is devoted to the derivation of uniform estimates being independent of the particular set shape. However, we will not restrict ourselves to a specific class of sets with certain properties, but rather show that for a generic domain one may construct a partition of the set such that the shape of each component can be controlled. The main result of this contribution is the following.

\begin{theorem}\label{th: main result smooth}
Let $\theta>0$.  Then there is  $\varrho=\varrho(\theta)>0$  such that the following holds: For all open, bounded and simply connected sets $\Omega \subset \R^2$ with $C^1$-boundary there is a  partition  $\Omega = \Omega_1 \cup \ldots \cup \Omega_N$ (up to a set of negligible measure) such that the sets  $\Omega_1,\ldots,\Omega_{N}$ are $\varrho$-John domains with Lipschitz boundary and 
\begin{align}\label{eq: addition-mainXXXXX}
\sum^N_{j=1}\mathcal{H}^1(\partial \Omega_j) \le (1+ \theta)\mathcal{H}^1(\partial \Omega).
\end{align}
\end{theorem}

Loosely speaking, the result states that in spite of the fact that there is no uniform control of the John constant for generic domains, it is at least possible to establish uniform estimates  locally in certain regions of the set. Here it is essential that the fineness of the partition can be bounded in terms of the length of the boundary of  $\Omega$. The original motivation for the derivation of Theorem \ref{th: main result smooth} is a piecewise Korn inequality \cite{Friedrich:15-4}  for special functions of bounded deformation (see \cite{Ambrosio-Coscia-Dal Maso:1997, Bellettini-Coscia-DalMaso:98}). We hope, however, that the result may be also applied in various other situations  due to the fact that John domains are a very general class and indeed many estimates only depend on the John constant (cf. \cite{Diening}).

It is a natural question if it is  possible to derive a partition of the form  \eqref{eq: addition-mainXXXXX} into sets satisfying more specific properties, e.g. convexity. By constructing an example related to Koch's snowflake we see, however, that in general this is not the case and similarly as in the results for the validity of Poincar\'e's and Korn's inequality (again see \cite{Acosta, Buckley}) also in the present context John domains appear to be an appropriate notion. 

Let us remark the the regularity assumption in Theorem \ref{th: main result smooth} is no real restriction as in many applications domains can be approximated by smooth sets (see \cite[Theorem 3.42]{Ambrosio-Fusco-Pallara:2000}) or discontinuities can be regularized by density arguments (see \cite{Chambolle:2004, Cortesani-Toader:1999}). Moreover, the result may be generalized to sets with Lipschitz boundary whose complements have a uniformly bounded number of connected components (see Theorem \ref{th: main part2}), which is a frequently used condition  for various models in fracture mechanics or shape optimization (cf.  \cite{Bucur2, Chambolle:2003, DM-Toa, Sverak}). However, the limitation to sets with a specific topology is crucial as without a requirement of this type the problem is essentially, again up to a density argument, equivalent to the derivation of a version of Theorem \ref{th: main result smooth} in the space of functions of bounded variation. This  is an even more challenging issue and we refer to \cite{Friedrich:15-4}  for a deeper analysis.

The essential step in the proof of Theorem \ref{th: main result smooth} is the derivation of a version  for polygons and the general case then  follows by approximation of regular sets. Although the methods we apply are rather elementary, the proof is comparably long and technical. Therefore, we restrict our decomposition scheme and analysis to a planar setting as in higher dimensions an analogous treatment of polyhedra leads to further technical difficulties. Let us remark, however, that based on Theorem \ref{eq: addition-mainXXXXX} in \cite{Friedrich:15-4}  various estimates of Korn and Korn-Poincar\'e type are derived, which hold in arbitrary space dimension.

Our strategy is twofold. We introduce two special subclasses of polygons, which we call \emph{semiconvex polygons} and \emph{rotund polygons}. We then show that (1) each polygon can be partitioned into semiconvex and rotund polygons  and (2) the specific characteristics of these subclasses of polygons are essentially equivalent to the property of John domains. 

Loosely speaking, in semiconvex polygons concave vertices are not `too close to opposite segments of the boundary' (see Definition \ref{def: semi}) and rotund polygons contain a ball whose diameter is comparable to the diameter of the polygon (see Definition \ref{def: qussc}). The decomposition scheme presented below is based on the idea to separate the domain by segments and in this context it is crucial that (1) by an iterative partition we do not violate properties which have already been established in a previous step and (2) the overall length of added segments is controllable in terms of $\mathcal{H}^1(\partial \Omega)$.

The proof that semiconvex,  rotund polygons are John domains for a John constant only depending on $\theta$ is constructive by defining appropriate piecewise affine curves between  generic points of the domain. Hereby we crucially exploit the fact that concave vertices are not `too close to opposite parts of the boundary' and that polygons are not `too thin'. Despite the specific properties of the subclasses of polygons we still have to face additional difficulties concerning the geometry of the curves, which may, e.g., partially have the form of a  helix.

The paper is organized as follows. In Section \ref{sec john-pre} we first recall the definition of John domains and state fundamental properties. In Section \ref{sec main-pre} we present a version of Theorem \ref{th: main result smooth} for polygons and give a more thorough overview of the proof. Here we also discuss an example  giving some intuition why John domains appear to be the appropriate notion for the formulation of the problem. In Section \ref{sec: prep} we introduce basic notation.

The subsequent sections are then devoted to the derivation of the result for polygons. In Section \ref{sec: semi} we introduce the notion of semiconvex polygons, prove basic properties and present a decomposition scheme. Afterwards, in Section \ref{sec: circ} 
we provide a fine analysis on the position of concave vertices and see that semiconvex polygons essentially coincide with convex polygons up to at most two small regions. In spite of their special structure, convex polygons are not necessarily rotund and we therefore discuss a further method to partition convex polygons. Finally, in Section \ref{sec: equi}  we prove   that semiconvex and rotund polygons are John domains with controllable John constant.   

In Section \ref{sec: proofi} we extend our findings to sets with $C^1$-boundary and in Section \ref{sec: gen} we discuss a variant of Theorem \ref{th: main result smooth} for sets with Lipschitz boundary allowing for a bounded number of components of the complement. Here we also present a piecewise Korn inequality as an application of our main result.

\section{Preliminaries}\label{sec: preli}

\subsection{John domains}\label{sec john-pre}

We first introduce the notion of \emph{John domains} and state some basic properties. Consider rectifiable curves $\gamma:[0,l(\gamma)]\to \R^d$ with length $l(\gamma)$ and assume
that they are parameterized by arc length. 
For $0 < \eta < 1$ we define the $\eta$-\emph{cigar} by
\begin{align}\label{eq: cigar}
{\rm cig}(\gamma, \eta) := \bigcup_{t \in [0,l(\gamma)]} B(\gamma(t),\eta\min\lbrace t, l(\gamma)-t \rbrace),
\end{align}
where $B(x,r) \subset \R^d$ denotes the open  ball with radius $r \ge 0$ and midpoint $x \in \R^d$. Likewise, we define 
the $\eta$-\emph{carrot} by
\begin{align}\label{eq: car}
{\rm car}(\gamma, \eta) := \bigcup_{t \in [0,l(\gamma)]} B(\gamma(t),\eta t ).
\end{align}

\begin{definition}\label{def: John}
Let $\varrho >0$. We say a bounded domain $\Omega \subset \R^d$ is a $\varrho$-John domain if there is a point $p \in \Omega$ such that for all $x \in \Omega \setminus \lbrace p \rbrace$ there is a rectifiable curve $\gamma: [0, l(\gamma)] \to \Omega$ with $\gamma(0) = x$ and $\gamma(l(\gamma)) = p$ such that ${\rm car}(\gamma, \varrho)  \subset \Omega$. 
\end{definition}

The point $p$ will be called the  \emph{John center} and $\varrho$ is the \emph{John constant}. Domains of this form were introduced by John  \cite{John:1961} to study problems in elasticity theory. The term was first used by Martio and Sarvas   \cite{Martio-Sarvas:1978}. Roughly speaking, a domain is a John domain if it is possible to connect two arbitrary points without getting too close to the boundary of the domain.

\begin{rem}\label{rem: john}
{\normalfont
A lot of different equivalent definitions can be found in \cite{Nakki-Vaisala:1991}. We will  also use the following characterization: a bounded domain $\Omega$ is a $\varrho$-John domain if for each pair of distinct points $x_1,x_2 \in \Omega$ there   is a curve $\gamma: [0, l(\gamma)] \to \Omega$ with $\gamma(0) = x_1$ and $\gamma(l(\gamma)) = x_2$ such that ${\rm cig}(\gamma, \varrho)  \subset \Omega$. Such a curve will be called \emph{John curve between $x_1$ and $x_2$.} 
}
\end{rem}

The class of John domains is much larger than Lipschitz domains and contains sets with fractal boundaries or internal cusps,  while the formation  of external cusps is excluded.  For instance the interior of Koch's \emph{snowflake} is a John domain. 
We state a simple property (see e.g. \cite{Vaisala:2000}).

\begin{lemma}\label{lemma: plump}
Let $\Omega$ be a $\varrho$-John domain. Then for each $x \in \Omega$ and $r >0$ with $\Omega \setminus B(x,r) \neq \emptyset$, there is $z \in \overline{B(x,r)}$ with $B(z,\frac{1}{2}\varrho r) \subset \Omega$. 
\end{lemma}

Our main result will be first established for polygons. To prove Theorem \ref{th: main result smooth}, we then need to combine different John domains so that the unions are still John domains. In  \cite{Vaisala:2000} we find the following lemma.

\begin{lemma}\label{lemma: Johna}
Let $\varrho,c_0>0$. There is $\varrho' = \varrho'(\varrho,c_0)$ such that the following holds:

(i) If $D_1,D_2 \subset \R^d$ are $\varrho$-John domains with $\min\lbrace|D_1|, |D_2|\rbrace   \le c_0 |D_1 \cap D_2|$, then $D_1 \cup D_2$ is a $\varrho'$-John domain.

(ii) If $D_0, D_1, \ldots,$ is a sequence of $\varrho$-John domains in $\R^d$ with $|D_j| \le c_0|D_0 \cap D_j|$ for all $j \ge 1$, then  $\bigcup_{j \ge 0} D_j$ is a $\varrho'$-John domain.

\end{lemma}

\subsection{Formulation of the main result for polygons}\label{sec main-pre}

The general strategy in this article is to derive the partition result first for polygons, which is easier due to the specific geometry of the boundary. In this section we present the main result for polygons and give an overview of the proof. 

Our partition technique for polygons will differ from widely used algorithms as triangulation, trapezoidalization or the Hertel and Mehlhorn Algorithm (see \cite{{ORourke:94}})  in the sense that we do not provide an optimal partition (concerning number of pieces or runtime), but one where the length of the boundary of all polygons is comparable to the length of the boundary of the original polygon. 

We consider sets $P \subset \R^2$ being the region enclosed by a simple polygon. For convenience sets of this form will be called \emph{polygons} in the following although the notion typically refers only to the boundary of such sets. We always assume that polygons are closed.  We notice that, according to our
definition, every polygon $P$ is simply connected and coincides with the closure of its
interior, which is nonempty. In particular the Lebesgue measure $|P|$ of $P$ is strictly positive.

 We intent to prove the following theorem.

\begin{theorem}\label{th: main part}
Let $\eps, \theta>0$.  Then  it exists  $\varrho=\varrho(\theta)>0$  such that  for all  polygons $P$ there is a  partition  $P = P_0 \cup \ldots \cup P_N$ with $\mathcal{H}^1(\partial P_0) \le \eps$ and the polygons $P_1,\ldots,P_{N}$ are $\varrho$-John domains satisfying 
\begin{align}\label{eq: addition-main}
\sum^N_{j=1}\mathcal{H}^1(\partial P_j) \le  (1+ \theta) \mathcal{H}^1(\partial P).
\end{align}
\end{theorem}

We start with a short outline of the proof. In particular, we indicate how an arbitrary polygon may be partitioned to satisfy the condition in Definition \ref{def: John} for a  John constant $\varrho$. 

First of all, the   property of $\varrho$-John domains  may be violated if the polygon has a `star shape', i.e. there are concave vertices for which the distance to other concave vertices or opposite segments of the boundary is small. We see that if this distance is too small, we can partition the polygon by introducing a short segment between a concave vertex and another point of the boundary. By this procedure we construct what we call \emph{semiconvex polygons} (see Section \ref{sec: semi}, in particular Definition \ref{def: semi}). Intuitively, such sets have the property that, separating the set by a short segment between a concave vertex and another point of the boundary, the `bulk part' of the polygon lies on one side.

Clearly, for convex sets it is much easier to satisfy  the condition in Definition \ref{def: John}. It turns out, however, that even  a convex polygon is possibly not a $\varrho$-John domains if the set is long and thin or has small interior angles. The presence of the latter phenomenon cannot be avoided and therefore the introduction of the set $P_0$ in Theorem \ref{th: main part} is possibly necessary. To tackle the first problem, we introduce so called \emph{rotund polygons} (see Section \ref{sec: circ}) which are sets containing a ball whose size is comparable to the diameter of the set. We then show that convex polygons can be partitioned into rotund polygons up to a small exceptional set (see Lemma \ref{lemma: convpart}). Finally, this kind of partition can also be performed for semiconvex polygons, which is related to the fact that a semiconvex polygon, which is not already rotund, coincides with a convex polygon up to at most two small regions (see Theorem \ref{th: new}).

After combining the above described partitions we show in Section \ref{sec: equi} that semiconvex and rotund polygons are indeed $\varrho$-John domains for a constant $\varrho=\varrho(\theta)$, which essentially only depends on the length of the additional  boundary induced by the partition (cf. \eqref{eq: addition-main}). The basic idea is to take a shortest path between two points (which will `touch' the boundary of the polygon in concave vertices) and to modify this path in such a way that the condition in Definition \ref{def: John} is satisfied. To do this, it is essential that (1) the polygons contain a ball whose size is comparable to the diameter of the set and (2) concave vertices are `not too close to opposite parts of the boundary'.

We remark that the definitions and terms of the  subclasses of polygons introduced in the following sections (see Section \ref{sec: semi}, Section \ref{sec: circ}) are not taken from the literature but tailored for the present exposition in order to avoid the ongoing repetition of technical assumptions. Let us also remark that, once the basic definition of semiconvex and rotund polygons have been internalized, Section \ref{sec: semi}--Section  \ref{sec: equi} can be read rather independently from each other.

\smallskip

Before we start to prove Theorem \ref{th: main part}, let us note that it does not appear to be possible to provide a partition for which the sets satisfy a stronger property than the one given in Definition \ref{def: John}. To give some intuition, we consider the following example being a modification of \emph{Koch's snowflake}.

\begin{example}\label{eq:ex}
{\normalfont

 Let $0<\eta<1$. Let $S_0$ be an equilateral triangle. As in the construction of Koch's snowflake we replace the middle third of each segment by two segments of equal length which enclose an angle $\frac{\pi}{3}$ with the original segment. Hereby, we obtain $S_1$. Then $S_2$ is obtained by replacing the middle third of each segment of $S_1$ by two segments which enclose an angle $\frac{\pi}{3}\eta$ with the original one. We continue with this construction where in the definition of $S_i$ the new segments enclose an angle  $\frac{\pi}{3}\eta^{i-1}$ with the original ones. 
 }
 \end{example}

Although the construction is very similar to the one of  Koch's snowflake, we find $\mathcal{H}^1(\partial S_i) \le C$ for all $i \in \N$ for some $C=C(\eta)$. Moreover, one can show that all $S_i$ are $\varrho$-John domains for some $\varrho>0$. Let us assume that the polygon $S_i$ for $i$ large could be partitioned into sets with `better properties' (e.g. convexity). Due to the geometry of $S_i$ we note that after separating $S_i$ into two sets by a segment there is one set which essentially has the same shape as $S_i$. Consequently, to derive a partition into sets with more specific properties, it appears to be necessary to introduce all boundaries $\bigcup_{j \le i-1} \partial S_j$. This, however, violates \eqref{eq: addition-main}.

\subsection{Notation}\label{sec: prep}

Let us fix the main notations for polygons which will be used in the following proof of Theorem \ref{th: main part}. Recall that polygons $P$ are always assumed to be closed subsets of $\R^2$. We denote the \emph{vertices} of $P$ by $\mathcal{V}_ P$ and for $v \in \mathcal{V}_P$ we let $\sphericalangle(v,P)$ be the corresponding  interior angle.     A vertex $v\in \mathcal{V}_ P$ with  $\sphericalangle(v,P)>\pi$ is called \emph{concave}, otherwise \emph{convex}. Denote the subset of concave vertices by $\mathcal{V}_P'$.    

 Sometimes we will understand vertices $v$ as complex numbers and let ${\rm arg}(v) \in [0,2\pi)$ be the phase of the complex number so that $v = |v| e^{i\,{\rm arg}(v)}$.  For $\varphi \in [0,2\pi)$ we denote by $v + \R_+ e^{i\varphi}$ open half lines with initial point $v$, where $\R_+ = (0,\infty)$.   The line segment between two given points $p_1,p_2 \in \R^2$ is denoted by $[p_1;p_2]$ and $|[p_1;p_2]|$ is its length.  For a segment $|[p_1;p_2]|$ we also introduce the notation (recall \eqref{eq: cigar})
$$
{\rm cig}([p_1;p_2],\eta) := {\rm cig}(\gamma,\eta),
$$
where $\gamma:[0,l(\gamma)] \to [p_1;p_2]$ is the (affine) curve, parametrized by arc length, with $\gamma(0) = p_1$, $\gamma(l(\gamma)) = p_2$ and length $l(\gamma) = |[p_1;p_2]|$.  Moreover, we define the \emph{visible region} of $[v;w]$ by  $${\rm cig}_P([v;w],\eta)=  \big\{ x \in \overline{{\rm cig}([v;w],\eta)}:  \ \exists \, p \in [v;w] \text{ s.t. } [p;x] \subset P\big\}$$
(see Figure \ref{semi1} below). We define an  \emph{intrinsic metric}  on $P$ by
$$d_P(p,p') = \min\lbrace l(\gamma): \gamma:[0,l(\gamma)]\to P \text{ Lipschitz curve with } \gamma(0)=p, \  \gamma(l(\gamma))=p'\rbrace $$
for $p,p' \in P$, where the curves are always assumed to be parameterized by arc length.  We notice that the minimum exists as $P$ is   closed and that it is attained by a piecewise affine curve, where the endpoints of each segment lie in $\mathcal{V}'_P \cup \lbrace p,p' \rbrace$. Likewise, for $p \in P$ and $S \subset P$ we let $\dist_P(p,S) = \inf_{p' \in S} d_P(p,p')$. Let the  \emph{intrinsic diameter}  of a polygon be given by
$$
d(P) = \max_{p,p' \in P} d_P(p,p').
$$
We find $d(P) \le \frac{1}{2}\mathcal{H}^1(\partial P)$ by considering a pair $p,p'$ maximizing $d_P(p,p')$ and the corresponding piecewise affine curve. The following definition will be used frequently.

\begin{definition}\label{def: induce}
Let $P$ be a polygon. We say a segment $[p;q] \subset P$ with $p,q \in \partial P$ induces a partition of $P$ if there are two polygons $Q_1,Q_2$ with $P = Q_1 \cup Q_2$  and   $[p;q] = Q_1 \cap Q_2$. 
\end{definition}

Note that, according to our definition of polygon, we have $|Q_1|, |Q_2|>0$. Moreover, $[p; q] = \partial Q_1 \cap \partial Q_2$ and every continuous path connecting a point of $Q_1$ with
a point of $Q_2$ must meet the segment $[p; q]$.

\section{Semiconvex polygons}\label{sec: semi}

We first refine Definition \ref{def: induce}.

\begin{definition}\label{def: induce2}
Let $\eta>0$ and $P$ be a polygon. We say a segment $[v;w]$  between a concave vertex $v \in \mathcal{V}_P'$ and some $w \in \partial P$ which induces a partition of $P = Q_1 \cup Q_2$ according to Definition \ref{def: induce} satisfies the \textit{segmentation property} {\rm (SP)} if 
\begin{align}\label{eq: cap}
\mathcal{V}_P' \cap {\rm cig}_P([v;w],\eta) \subset \lbrace v,w \rbrace 
\end{align}
 and  for $i=1,2$ 
\begin{align}\label{eq: capXX}
Q_i \text{ is a triangle } \ \ \ \Rightarrow  \ \ \    \sphericalangle(v,Q_i)> \frac{1}{2}\arcsin\eta.
\end{align}
\end{definition}

These technical conditions  are necessary to avoid the formation of geometrical artefacts in the partition process  in Section \ref{sec: semi2}  such as degenerated triangles and polygons where a concave vertex is very close to an opposite side.   We now introduce the notion of \emph{semiconvexity}.

\vspace{-0.20cm}

\begin{figure}[H]
\centering
\begin{overpic}[width=0.8\linewidth,clip]{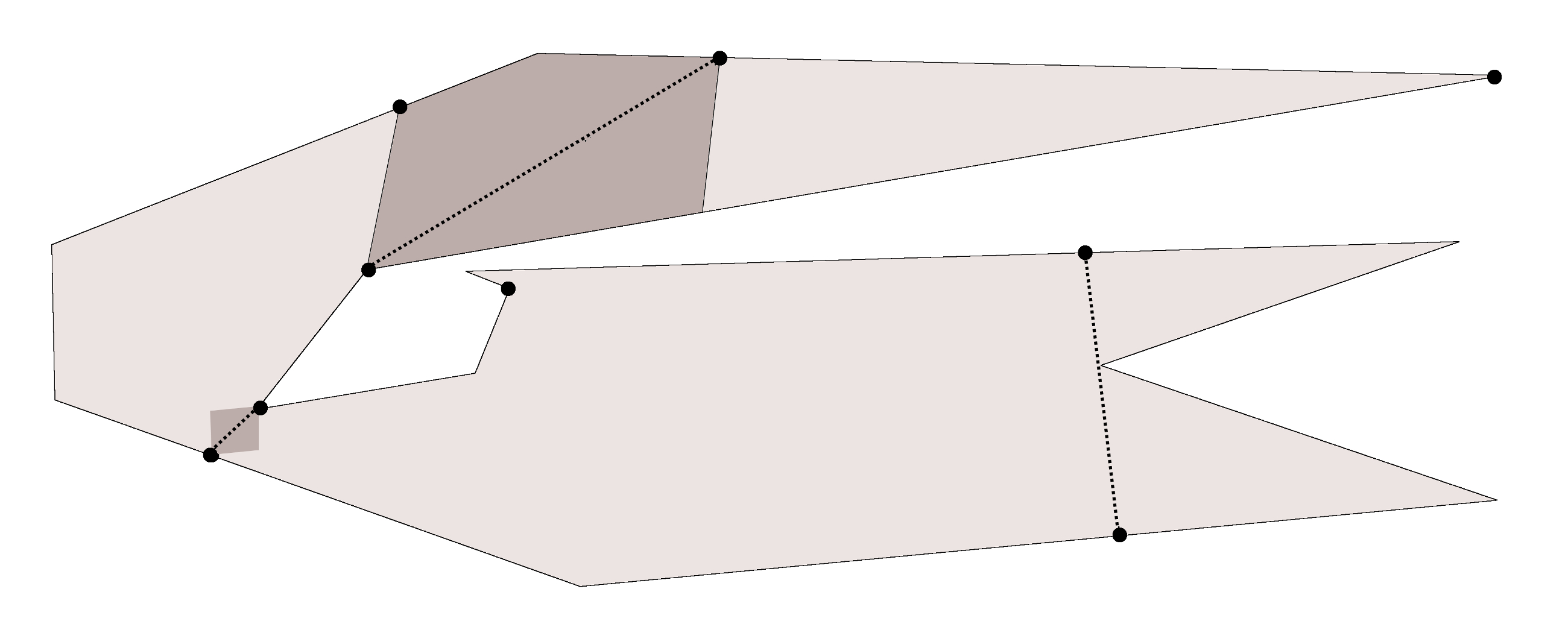}
\put(66,72){{$v$}}
\put(72,108){{$w'$}}
\put(150,119){{$w$}}
\put(109,66){{$u$}}
\put(315,112){{$p$}}

\put(47,47){{$\bar{v}$}}
\put(39,23){{$\bar{w}$}}
\put(223,82){{$\bar{p}$}}
\put(222,21){{$\bar{q}$}}

\end{overpic}
\caption{In dark gray we depicted ${\rm cig}_P([v;w],\eta)$ and  ${\rm cig}_P([\bar{v};\bar{w}],\eta)$. Observe that $u \notin {\rm cig}_P([v;w],\eta)$ although $u \in {\rm cig}([v;w],\eta)$. Consequently, $[v;w]$ satisfies \eqref{eq: cap} but not (SP) due to the small interior angle at $v$ in the triangle formed by $v$,$w$,$p$. The segment $[v;w']$ satisfies (SP). The segment $[\bar{v},\bar{w}]$ is a typical example of a segment which satisfies (SP) such that condition \eqref{eq: semi} is violated. Note that $[\bar{p};\bar{q}]$ does not induce a partition into two simple polygons.} \label{semi1}
\end{figure}
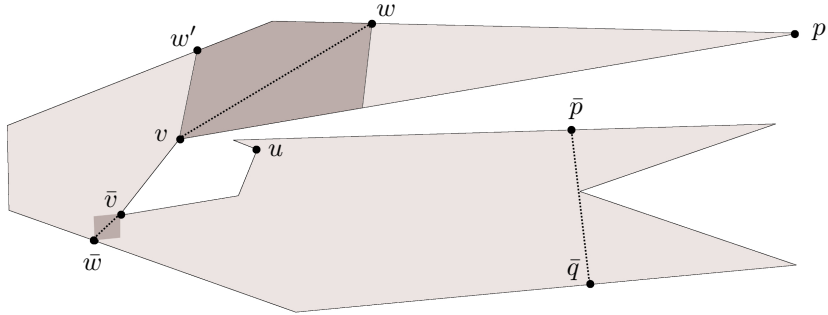

\begin{definition}\label{def: semi}
Let $0< \vartheta, \eta <1$. 

(i) We say a polygon $P$ is $\vartheta$-semiconvex if for each segment $[v;w]$ between a concave vertex $v \in \mathcal{V}'_P$ and some $w \in \partial P$ which induces a partition $P = Q_1 \cup Q_2$  one has
\begin{align}\label{eq: semi}
|[v;w]| \ge \vartheta \min_{k=1,2} d(Q_k).
\end{align}  

(ii) We say a polygon $P$ is {\rm (SP)}-$\vartheta$-semiconvex if for each segment $[v;w]$ between a concave vertex $v \in \mathcal{V}'_P$ and some $w \in \partial P$ which induces a partition $P = Q_1 \cup Q_2$ and satisfies  {\rm (SP)}  one has \eqref{eq: semi}.
\end{definition}

For simplicity we will often drop the parameters and will call a polygon semiconvex and (SP)-semiconvex if no confusion arises.

\begin{rem}
{\normalfont
Intuitively, the definition states that, separating the set by a short segment between a concave vertex and another point of the boundary, the `bulk part' of the polygon lies on one side.  The semiconvexity of a polygon together with rotundness considered in Section \ref{sec: circ} is the essential property to control the John constant of polygons. We note that in \eqref{eq: semi} the  intrinsic diameter is the suitable notion and cannot be replaced by the length of the boundary although it seems to be another natural choice. To see this,  consider Koch's \emph{snowflake} which is a John domain with finite { intrinsic diameter} but whose boundary is of infinite $\mathcal{H}^1$-measure.
}
\end{rem}

In Section \ref{sec: semi1} we study the relation between semiconvex and (SP)-semiconvex polygons deriving that the notions are very similar. In Section \ref{sec: circ}--Section \ref{se: mainproof} we will only need the concept of  semiconvex polygons. However, for the partition of polygons  into semiconvex polygons performed in Section \ref{sec: semi2} it is  convenient to consider also the more technical notion in Definition \ref{def: semi}(ii).

\subsection{Properties of semiconvex polygons}\label{sec: semi1}

By definition we clearly have that each semiconvex polygon is also (SP)-semiconvex. We now investigate the reverse direction.

\begin{theorem}\label{lemma: partpol*}
Let $0< \vartheta, \eta <1$ with $\vartheta \le \frac{1}{2}\eta$. Then for $\eta>0$ small enough there is some $\bar{\vartheta} = \bar{\vartheta}(\vartheta) \le \vartheta$ such that each {\rm (SP)-}$\vartheta$-semiconvex polygon is $\bar{\vartheta}$-semiconvex.
\end{theorem}

\Proof Let $P$ be a {\rm (SP)-}$\vartheta$-semiconvex polygon.  Let   $v \in \mathcal{V}'_P$ and some $w \in \partial P$ be given inducing a partition of $P$. The  goal is to confirm \eqref{eq: semi} for $[v;w]$. To this end, we will construct a chain of segments consisting of concave  vertices and combining  $v$ with $w$  such that each segment satisfies (SP) and therefore \eqref{eq: semi} is applicable by assumption.\\

\noindent  \smallskip
\textit{Step 1: Cigar condition}

\noindent We first assume that $[v;w]$ induces a partition $P = Q_1 \cup Q_2$ and that 
\begin{align}\label{eq: {eq: cap2}}
\mathcal{V}_P'   \cap   {\rm cig}_P([v;w], 2\eta)\subset \lbrace v,w\rbrace.
\end{align}
(Compare with \eqref{eq: cap} and note that in contrast to Definition \ref{def: induce2} we do not require \eqref{eq: capXX}.) We show that
\begin{align}\label{eq: direct connec}
|[v;w]|  \ge  \frac{\vartheta}{2} \min_{k=1,2} d(Q_k).
\end{align}
We distinguish the following cases: \smallskip \\
(a) If each $Q_k$  is either not a triangle or a triangle where the interior angle at $v$ exceeds $\alpha_\eta := \frac{1}{2}\arcsin\eta$, we find that \eqref{eq: cap}-\eqref{eq: capXX} hold and thus $|[v;w]|  \ge \vartheta \min_{k=1,2} d(Q_k)$ by Definition \ref{def: semi}(ii). \smallskip \\ 
(b) Otherwise, we can suppose without restriction that $Q_1$ is a triangle consisting of the vertices $v,w,p$ with $\sphericalangle(v,Q_1)\le \alpha_\eta$ (cf. Figure \ref{semi1}). Thus, if $\eta$ small enough, we get $d(Q_1)=\max\lbrace|[v;p]|,|[v;w]|\rbrace$ and may assume $|[v;p]| \ge |[v;w]|$ as otherwise \eqref{eq: direct connec} follows directly. If $\sphericalangle(w,Q_1) \le \pi-2\alpha_\eta$, we apply the sine rule $\frac{|[v;p]|}{\sin\sphericalangle(w,Q_1)} = \frac{|[v;w]|}{\sin\sphericalangle(p,Q_1)}$ and the fact that $\sphericalangle(p,Q_1) \ge \alpha_\eta$ to see 
\begin{align*}
d(Q_1) = |[v;p]|  =  \frac{\sin\sphericalangle(w,Q_1)}{\sin\sphericalangle(p,Q_1)}|[v;w]| \le \frac{|[v,w]|}{\sin\alpha_\eta} \le \frac{4}{\eta}|[v,w]| \le \frac{2}{\vartheta}|[v,w]|
\end{align*}
for $\eta$ small, where we used $\sin\alpha_\eta \ge \frac{1}{4}\eta$ by a Taylor expansion and $\vartheta \le \frac{1}{2}\eta$.
\smallskip  \\
(c)  Otherwise, we have $\sphericalangle(w,Q_1) > \pi-2\alpha_\eta$. First suppose $w \in \mathcal{V}_P'$, which means that  we can change the roles of $v$ and $w$. We see that \eqref{eq: cap} holds by assumption. Moreover, we have  $\sphericalangle(w,Q_1) > \pi - 2\alpha_\eta > \alpha_\eta$   for $\eta$ small  and that $Q_2$ is not a triangle since $P$ has at least five vertices due to $\lbrace v,w\rbrace \subset \mathcal{V}_P'$. Consequently, also \eqref{eq: capXX} holds and we can proceed as in (a) to find $|[v;w]|  \ge \vartheta \min_{k=1,2} d(Q_k)$.\smallskip  \\
Observe that  in (b) we used a purely geometrical argument and in (a),(c) we only showed that \eqref{eq: capXX} holds, whereby Definition \ref{def: semi}(ii) was applicable. In the following last case, however, we will explicitly  use \eqref{eq: {eq: cap2}}.  \smallskip  \\
(d)  Finally, we suppose that $\sphericalangle(w,Q_1) > \pi-2\alpha_\eta$ and that $w$ is not a concave vertex. Understanding the vertices as complex numbers we define the phase $\varphi_0 = {\rm arg}(w-v)$. Let $f: { D} \to \R^2$ so that $f(\varphi)$ denotes the closest point to $v$ on $(v + \R_+ e^{i(\varphi_0+\varphi)}) \cap \partial P$, {  where $D \subset [-\pi,\pi)$ contains a neighborhood of $0$ and satisfies $|D| = \sphericalangle(v,P)$. (Recall $\R_+ = (0,\infty)$.)} For $\varphi>0$ small let $\triangle_\varphi$ the triangle formed by $v$, $p$, $f(\varphi)$ and up to changing the sign of $\varphi$ we may assume that $\sphericalangle(v,\triangle_\varphi) > \sphericalangle(v,Q_1)$ for $\varphi>0$ small. Observe that due to the fact that $w$ is not a concave vertex and $\sphericalangle(w,Q_1) > \pi-2\alpha_\eta = \pi -  \arcsin\eta$ we have 
$$f(\varphi) \in {\rm cig}_P([v;w],2\eta)$$ 
for $\varphi$ small. This then implies $f(\varphi) \notin \mathcal{V}_P'$ for $\varphi \in [0,2\alpha_\eta]$ since otherwise \eqref{eq: {eq: cap2}} would be violated. Consequently, letting $w' = f(2\alpha_\eta)$ we find that $[v;w']$ induces a partition $P = Q_1' \cup Q_2'$, where the sets are labeled such that $p \in Q_1'$. Moreover,  $[v;w']$ satisfies (SP). In fact, the angle condition \eqref{eq: capXX} follows directly by construction. Moreover, we get ${\rm cig}_P([v;w'],\eta) \subset {\rm cig}_P([v;w],2\eta)$ and thus \eqref{eq: cap} follows from \eqref{eq: {eq: cap2}} and the fact that $w \notin \mathcal{V}_P'$. Consequently, as $P$ is {\rm (SP)-}$\vartheta$-semiconvex,   we obtain by \eqref{eq: semi}
$$|[v;w']|  \ge \vartheta \min_{k=1,2} d(Q_k').$$
Consider the convex polygon $\hat{P} := Q'_1 \cap Q_2$ and note that for $\eta$ small $d(\hat{P}) = |[v;w]|$ (see Figure \ref{semi1}) as well as $|[v;w']| \le |[v;w]|$. Moreover,   $\min_{k=1,2} d(Q_k)   \le \min_{k=1,2} d(Q_k') + d(\hat{P})$
 and therefore  we obtain since $\vartheta < 1$ 
\begin{align*}
\begin{split}
|[v;w]| &\ge \frac{1}{2}|[v;w']| + \frac{1}{2}|[v;w]|   \ge \frac{\vartheta}{2} \min_{k=1,2} d(Q_k) -  \frac{\vartheta}{2} |[v;w]| + \frac{1}{2}|[v;w]|\ge \frac{\vartheta}{2} \min_{k=1,2} d(Q_k).
\end{split}
\end{align*}

\noindent  \smallskip
\textit{Step 2: Chains of vertices}

\begin{figure}[H]
\centering

\begin{overpic}[width=0.75\linewidth,clip]{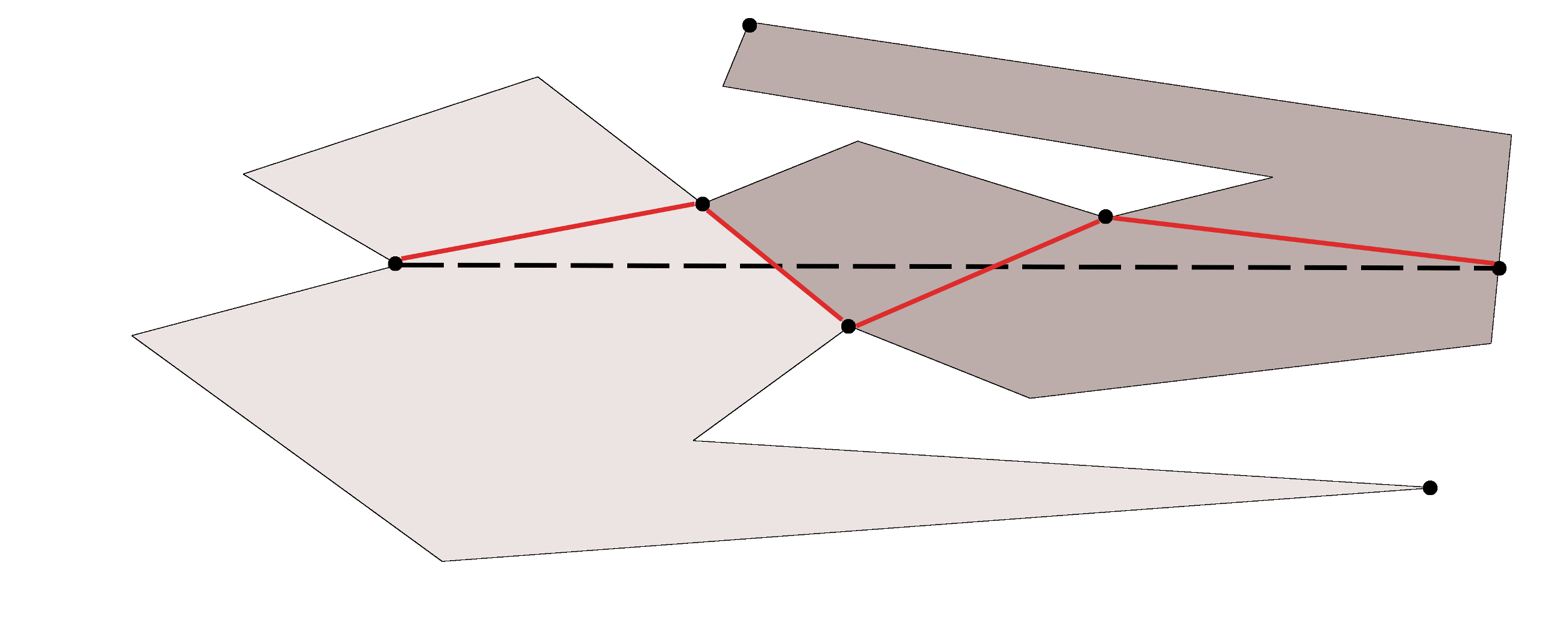}
\put(38,68){{$y_1=v$}}
\put(295,68){{$y_5=w$}}
\put(133,87){{$y_2$}}
\put(80,36){{$Q^{(2)}_1$}}
\put(260,83){{$Q^{(2)}_2$}}
\put(214,83){{$y_4$}}
\put(162,49){{$y_3$}}
\put(283,23){{$p_1^1$}}
\put(133,115){{$p_2^1$}}

\end{overpic}
\caption{The segments $[y_i;y_{i+1}]$ inducing partitions of $P$ are depicted in red. In light and dark gray the partition $Q_1^{(2)} \cup Q_2^{(2)}$ is sketched, where $p^1_j \in Q^{(2)}_j$ for $j=1,2$.} \label{semi2}
\end{figure}
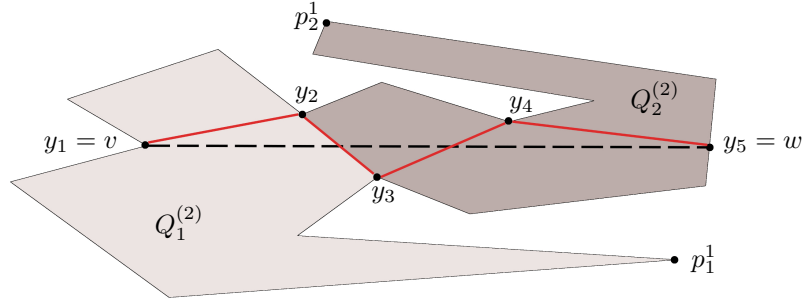

\noindent Now we only assume that $[v;w]$, $v \in \mathcal{V}_P'$, $w \in \partial P$, induces a partition of $P$. We construct a chain  $(y_1,\ldots,y_{n})$ between $v$ and $w$ with $y_1 = v$, $y_{n} = w$ and  $y_i \in \mathcal{V}_P'$ for $i=2,\ldots,n-1$ such that 
\begin{align}\label{eq: partpol223}
|[y_i;y_{i+1}]|\le 3|[v;w]| , \ \ \ \ \  d_P(v,y_{i})\le \frac{3}{2}|[v;w]| , \ \ \ i=1,\ldots,n-1
\end{align}
and the segments $[y_i;y_{i+1}] \subset P$ induce a partition satisfying \eqref{eq: {eq: cap2}} with $y_i$, $y_{i+1}$ in place of $v,w$ (cf. Figure \ref{semi2}). (See Section \ref{sec: prep}  for the definition of $d_P(v,y_{i})$.)

The strategy is to define the chain between $v$ and $w$ inductively. Let $\mathcal{C}_0 = (y^0_1,y^0_{2}) = (v,w)$ and assume $\mathcal{C}_k = (y^k_1,\ldots,y^k_{2+k})$ with $y^k_1 = v$, $y^k_{2+k} = w$  and $[y^k_j;y^k_{j+1}] \subset P$ for $j=1,\ldots,k+1$  has been constructed. If 
\begin{align}\label{eq: all j}
\mathcal{V}_P' \cap   {\rm cig}_P([y^k_j;y_{j+1}^k],2\eta)  \subset \lbrace y^k_j,y_{j+1}^k \rbrace \ \ \ \text{ for all }  j=1,\ldots,k+1,
\end{align}
we stop. Otherwise, we find some $J \in \lbrace 1, \ldots, k+1\rbrace$ and $\hat{v}_k \in \mathcal{V}_P' \setminus \lbrace y^k_{J}, y^k_{J+1} \rbrace$ such that $\hat{v}_k \in {\rm cig}_P([y^k_{J};y_{J+1}^k],2\eta)$ and $[y^k_{J};\hat{v}_k] \cup [\hat{v}_k;y^k_{J+1}] \subset P$. (Choose $\hat{v}_k$ as the concave vertex  in ${\rm cig}_P([y^k_{J};y_{J+1}^k],2\eta)$  with minimal distance to $[y^k_{J};y^k_{J+1}]$.) We define
$$\mathcal{C}_{k+1} = (y_1^k, \ldots, y_{J}^k, \hat{v}_k,y_{J + 1}^k, \ldots, y^k_{k+2}).$$
 Note that the triangle formed by $[y^k_{J};y^k_{J+1}]$ and $\hat{v}_k$ is contained in $P$ since $P$ is simply connected.   As in each step we choose a different $\hat{v}_k$ and $\# \mathcal{V}_P' < \infty$, after a finite number of steps we find a chain $(y_1,\ldots,y_n)$ such that \eqref{eq: all j} is satisfied. 

We now show that \eqref{eq: partpol223} holds. To this end, we fix $y_i$, $i=2,\ldots,n-1$, and identify the iteration steps that `led to the definition of $y_i$'. Let $k_0$ be the   index such that $\hat{v}_{k_0} := y_i \in \mathcal{C}_{k_0+1}$. Choose $J_0$ such that $y_i \in {\rm cig}_P(S_0,2\eta)$ with $S_0 = [y^{k_0}_{J_0};y_{J_0+1}^{k_0}]$. 

Assume steps $k_0>k_1 > \ldots >k_n$ and $(J_i)_{i=0}^n$ have been found with corresponding $\hat{v}_{k_i}$ such that $\hat{v}_{k_i} \in {\rm cig}_P(S_i,2\eta)$ with $S_i := [y^{k_i}_{J_i};y_{J_{i}+1}^{k_i}]$. 

We then choose the largest value $k_{n+1} < k_n$ such that one of the points $y^{k_n}_{J_n},y_{J_n+1}^{k_n}$ is not contained in $\mathcal{C}_{k_{n+1}}$, e.g. $y^{k_n}_{J_n} =: \hat{v}_{k_{n+1}}$. We then find $J_{n+1}$ such that $\hat{v}_{k_{n+1}} \in {\rm cig}_P(S_{n+1},2\eta)$ with $S_{n+1} = [y^{k_{n+1}}_{J_{n+1}};y_{J_{n+1}+1}^{k_{n+1}}]$, where one of the endpoints of $S_{n+1}$ coincides with $y^{k_n}_{J_n+1}$. For later purpose we note that $S_n$, $S_{n+1}$ have a common endpoint and $\hat{v}_{k_{n+1}}$ is an endpoint of $S_{n}$. Finally, after a finite number of steps, denoted by $N$, we arrive at $S_N = [v;w]$.  

\vspace{-0.2cm}
\begin{figure}[H]
\centering

\begin{overpic}[width=0.58\linewidth,clip]{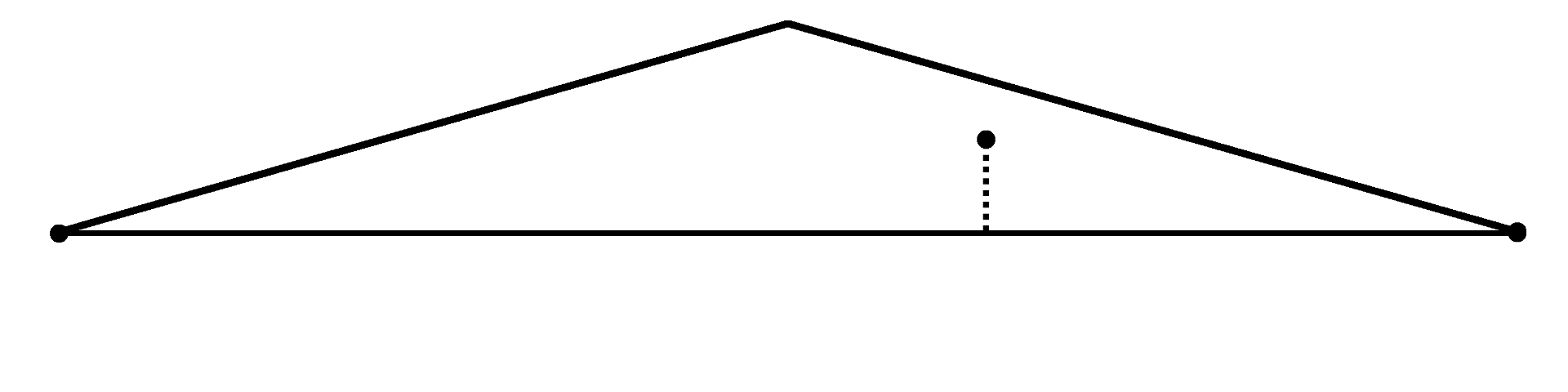}
\put(130,39){{$\hat{v}_{k_n}$}}
\put(-9,17){{$y^{k_n}_{J_n}$}}
\put(234,17){{$y^{k_n}_{J_n+1}$}}
\put(150,29){{$\big\}$}}
\put(149,22){\small{$\underbrace{ \ \  \ \ \ \ \ \ \ \ \ \ \ \ \ \ \ \ \ \ \ \ \ \ \   }$}}
\put(12,10){\small{$\underbrace{ \ \  \ \ \ \ \ \ \ \ \ \ \ \ \ \ \ \ \ \ \ \ \ \  \ \ \ \ \ \ \ \ \  \ \ \ \ \ \ \ \ \ \ \ \  \ \ \ \  \ \ \ \ \ \ \ \ \ \ \ \  \ \ \ \ \ \ \ \ \   }$}}
\put(156,29){{$a$}}
\put(186,10){{$c$}}
\put(119,-5){{$b$}}
\put(234,34){{$4\alpha_\eta$}}
\put(224,24) {\line(1,1){10}}

\end{overpic}
\caption{We set $a=  \dist(S_{n},\hat{v}_{k_n})$ and $b = \mathcal{H}^1(S_n)$. Elementary trigonometry yields $c \ge (\tan(4\alpha_\eta))^{-1} a$. Note that $S_{n-1}$ is the segment between $\hat{v}_{k_n}$ and the left or right endpoint of $S_n$.} \label{2}
\end{figure}
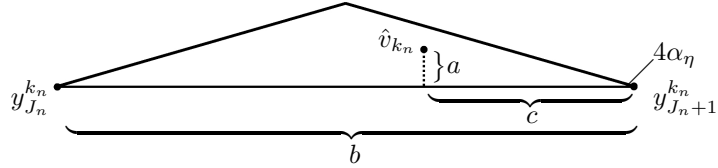

Recalling the geometry of ${\rm cig}(S_n,2\eta)$ an elementary computation yields that the angles at the endpoints of $S_n$ in the triangle formed by $S_n$ and $\hat{v}_{k_n}$ are larger than (cf. Figure \ref{2})
$$\varphi_n := \arctan(g_n) \ \ \ \ \text{with} \ \ \ \ g_n := \frac{ \dist(S_{n},\hat{v}_{k_n})}{\mathcal{H}^1(S_n) -  (\tan(4\alpha_\eta))^{-1}\dist(S_{n},\hat{v}_{k_n})}.$$
 We note $\dist(S_{n},\hat{v}_{k_n}) \le \frac{1}{2}\tan(4\alpha_\eta)\mathcal{H}^1(S_{n})$ and thus $g_n \le \tan(4\alpha_\eta)$. Recalling that the segments $S_{n-1}$, $S_n$ have one common endpoint (either $y^{k_n}_{J_n}$ or $ y^{k_n}_{J_n+1}$) and $\hat{v}_{k_n}$ is an endpoint of $S_{n-1}$, we find $\mathcal{H}^1(S_{n-1}) \le (\sin\varphi_n)^{-1} \dist(S_{n},\hat{v}_{k_n})$ for all $n=1,\ldots, N$. Then we obtain  by a Taylor expansion for $\eta$ small and some large $C>0$ independent of $\eta$ (observe that for small $x$ one has $\arctan(x), \sin(x) \approx x$)
\begin{align*}
\mathcal{H}^1(S_{n-1}) &\le \frac{1}{g_n - Cg_n^2}\dist(S_{n},\hat{v}_{k_n}) \le (1+Cg_n)\Big(\mathcal{H}^1(S_n) -  \frac{\dist(S_{n},\hat{v}_{k_n})}{\tan(4\alpha_\eta)}\Big)\\
& = \mathcal{H}^1(S_n) -   (\tan(4\alpha_\eta))^{-1} \dist(S_{n},\hat{v}_{k_n}) + C\dist(S_{n},\hat{v}_{k_n}).
\end{align*}
Note that $\dist(x,S_{n} ) \le \dist(\hat{v}_{k_n},S_n)$ for all $x \in S_{n-1}$. Then using the previous estimate and summing over all $n$ we find for $\eta$ small (such that $0<\tan(4\alpha_\eta)(1-C\tan(4\alpha_\eta))^{-1}\le \frac{1}{2}$)
\begin{align*}
d_P(y_i,[v;w]) &= d_P(y_i,S_N)\le \dist(\hat{v}_{k_0}, S_{0}) + \sum_{n=1}^N \max_{x \in S_{n-1}} \dist(x,S_n) \le \sum_{n=0}^N \dist(\hat{v}_{k_n}, S_{n}) \\& \le \frac{\tan(4\alpha_\eta)}{1-C\tan(4\alpha_\eta)}\sum_{n=1}^N(\mathcal{H}^1(S_{n}) - \mathcal{H}^1(S_{n-1})) + \dist(\hat{v}_{k_0}, S_{0})\\&\le \frac{\tan(4\alpha_\eta)}{1-C\tan(4\alpha_\eta)}(|[v;w]| - \mathcal{H}^1(S_{0})) + \frac{1}{2}\tan(4\alpha_\eta)\mathcal{H}^1(S_{0}) \le \frac{1}{2}|[v;w]|,
\end{align*}
where we used  $\dist(S_0, \hat{v}_{k_0}) \le \frac{1}{2} \tan(4\alpha_\eta)\mathcal{H}^1(S_0)$.  This together with the triangle inequality yields \eqref{eq: partpol223}.\\

\vspace{-0.3cm}

\noindent  \smallskip
\textit{Step 3: Semiconvexity}

\noindent We now show that $P$ is semiconvex by confirming \eqref{eq: semi} for the segment $[v;w]$ with $\bar{\vartheta} = (3+ 12\vartheta^{-1})^{-1}$. As each of the segments $[y_i;y_{i+1}]$ satisfies  \eqref{eq: {eq: cap2}} (with $y_i,y_{i+1}$ in place of $v,w$), we obtain by \eqref{eq: direct connec}
\begin{align}\label{eq: theta con}
|[y_i;y_{i+1}]|  \ge \frac{\vartheta}{2} \min_{k=1,2} d(Q^{(i)}_k)
\end{align}
for $i=1,\ldots,n-1$, where $P = Q^{(i)}_1 \cup Q^{(i)}_2$ is the corresponding partition. Let $P = Q_1 \cup Q_2$ be the partition induced by $[v;w]$. It suffices to consider the case $|[v;w]| \le \frac{1}{8} \min_{j=1,2} d(Q_j)$ as otherwise the assertion is clear provided we choose $\bar{\vartheta} \le \frac{1}{8}$. We choose $p^1_j, p^2_j \in Q_j$ with $\dist_{P}(p^1_j,p_j^2) =d(Q_j)$ and as $\dist_{P}(p^1_j,p^2_j) \le \dist_{P}(p^1_j,v) + \dist_{P}(p^2_j,v)$, we obtain possibly after relabeling $\dist_{P}(p^1_j,v) \ge \frac{1}{2}d(Q_j)$ for $j=1,2$. 

 Let $B = \lbrace x \in P: \dist_P(x,v) < 4 |[v;w]| \rbrace$.   We now show that two arbitrary points $q_1 \in Q_1 \setminus B$, $q_2 \in Q_2 \setminus B$ do not lie in the same connected component   of $P \setminus \bigcup^{n-1}_{i=1} [y_i;y_{i+1}]$.
 
 Indeed, let $T$ be a connected component of $P \setminus \bigcup^{n-1}_{i=1} [y_i;y_{i+1}]$. It suffices to show that $(T \setminus B) \cap Q_j = \emptyset$ for some $j=1,2$. If $T \subset Q_j$ for some $j=1,2$, this is clear. Otherwise, we find some $j=1,2$ such that  $T' := T \cap Q_j$ satisfies $\partial T' \subset [v;w] \cup \bigcup^{n-1}_{i=1} [y_i;y_{i+1}]$ (see also Figure \ref{semi2}). Now combining the two inequalities in \eqref{eq: partpol223}, we get $d_P(v,x) \le 3|[v;w]|$ for all $x \in \partial T'$. Then also   $d_P(v,x) \le 3|[v;w]|$ for all $x \in T'$ and this shows $(T \setminus B) \cap Q_j = T' \setminus B =\emptyset$.

Therefore, recalling $|[v;w]| \le \frac{1}{8} \min_{j} d(Q_j) \le \frac{1}{4} \min_{j}\dist_{P}(p^1_j,v)$, we find that $p_1^1 \in Q_1 \setminus B$ and $p_2^1 \in Q_2 \setminus B$ lie in different connected components of $P \setminus \bigcup^{n-1}_{i=1} [y_i;y_{i+1}]$. Thus, there is at least one $i=1,\ldots,i-1$ such that possibly after relabeling we have $p_1^1 \in Q^{(i)}_1$ and  $p_2^1 \in Q^{(i)}_2$. Using  \eqref{eq: partpol223}  we find
\begin{align*}
\min\nolimits_j d(Q_j^{(i)}) &\ge   \min\nolimits_j\dist_{P}(y_i,p^1_j) \\&\ge\min\nolimits_j \big( \dist_{P}(v,p^1_j) - \dist_{P}(y_i,v) \big)\ge \frac{1}{2}\min\nolimits_j d(Q_j) - \frac{3}{2}|[v;w]|.
\end{align*} 
By  \eqref{eq: partpol223} and  \eqref{eq: theta con} we conclude with $\bar{\vartheta} = (3+ 12\vartheta^{-1})^{-1}$ 
$$\min\nolimits_j d(Q_j) \le  3|[v;w]| + 4\vartheta^{-1} |[y_i;y_{i+1}]|  \le (3+ 12 \vartheta^{-1})|[v;w]| = \bar{\vartheta}^{-1}|[v;w]|.$$   
This shows \eqref{eq: semi} and concludes the proof. \eop

We now show that a similar property may derived if the condition in Definition \ref{def: semi}(ii) only holds on a part of $\partial P$. To this end, we need to introduce a further notion. Suppose $[v;w]$ induces a partition of $P= Q_1 \cup Q_2$ according to Definition \ref{def: induce}. We   define $N'(Q_j) = \# \lbrace u \in \mathcal{V}_P' \setminus \lbrace v,w\rbrace: u \in \partial Q_j \rbrace$ for $j=1,2$ and the auxiliary set
\begin{align}\label{eq: Qvw}
Q_{v,w} = \begin{cases} Q_1   & \text{if }  N'(Q_1) < N'(Q_2) \text{ or } N'(Q_1) = N'(Q_2), \ |Q_1|< |Q_2|, \\ 
Q_1 \cup Q_2 & \text{if }    N'(Q_1) = N'(Q_2), |Q_1|= |Q_2|, \\ 
Q_2   & \text{else}. 
\end{cases}
\end{align}

\begin{definition}\label{def: induce3}
 We say a segment $[v;w]$ satisfies the \textit{weak segmentation property} {\rm (WSP)} if in Definition \ref{def: induce2} condition  \eqref{eq: cap} is replaced by 
\begin{align}\label{eq: capW}
\mathcal{V}_P' \cap {\rm cig}_P([v;w],\eta) \cap Q_{v,w} \subset \lbrace v,w \rbrace. 
\end{align}
\end{definition}

 We note that for (WSP) we still require \eqref{eq: capXX}. Loosely speaking, condition \eqref{eq: capW} only concerns the part of the polygon containing less concave vertices and is thus in general weaker than \eqref{eq: cap}.

\begin{corollary}\label{cor: partpol*}
Let $0< \vartheta, \eta <1$ with $\vartheta \le \frac{1}{2}\eta$. Consider a polygon $P$ and suppose $[v;w]$  induces  a partition $P=Q_1 \cup Q_2$ satisfying {\rm (WSP)} and  
\begin{align}\label{eq: handing}
\text{ either } \  \ \ N'(Q_1) < N'(Q_2) \ \ \ \text{ or } \ \ \ N'(Q_1) = N'(Q_2), \ |Q_1|<|Q_2|.
\end{align}  
Assume that for each pair $v',w' \in  \mathcal{V}_{Q_1}'\cup \lbrace v,w \rbrace$  such that $[v';w']\neq[v;w]$ and $[v';w']$ induces a partition $P=Q_1' \cup Q_2'$ satisfying {\rm (WSP)} one has
\begin{align}\label{eq: semiXXX}
|[v';w']| \ge \vartheta \min_{k=1,2} d(Q_k').
\end{align} 
Then for $\eta>0$ small enough there is $\bar{\vartheta} = \bar{\vartheta}(\vartheta) \le \vartheta$ independent of $P$ such that each pair  $\bar{v},\bar{w}  \in  \mathcal{V}_{Q_1}' \cup \lbrace v,w \rbrace$ inducing a partition of $P = R_1 \cup R_2$ with 
\begin{align}\label{eq: handing2}
\begin{split}
(i)& \ \ [\bar{v};\bar{w}]\neq[v;w],\\ (ii) & \ \ w \in \lbrace \bar{v},\bar{w} \rbrace \ \ \Rightarrow \ \ [\bar{v};\bar{w}] \cap  {\rm cig}([v;w],\eta) = \emptyset
\end{split}
\end{align} 
  fulfills
\begin{align}\label{eq: special two}
|[\bar{v};\bar{w}]|  \ge \bar{\vartheta} \min_{k=1,2} d(R_k).
\end{align}
\end{corollary}

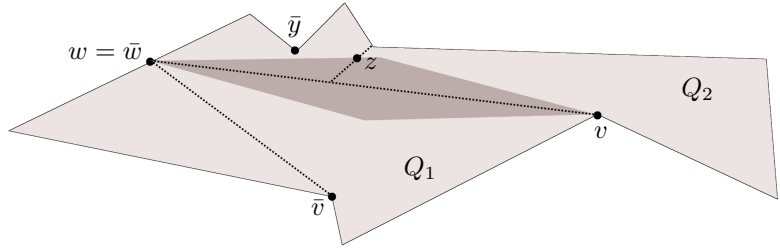
\begin{figure}[H]
\centering

\begin{overpic}[width=0.76\linewidth,clip]{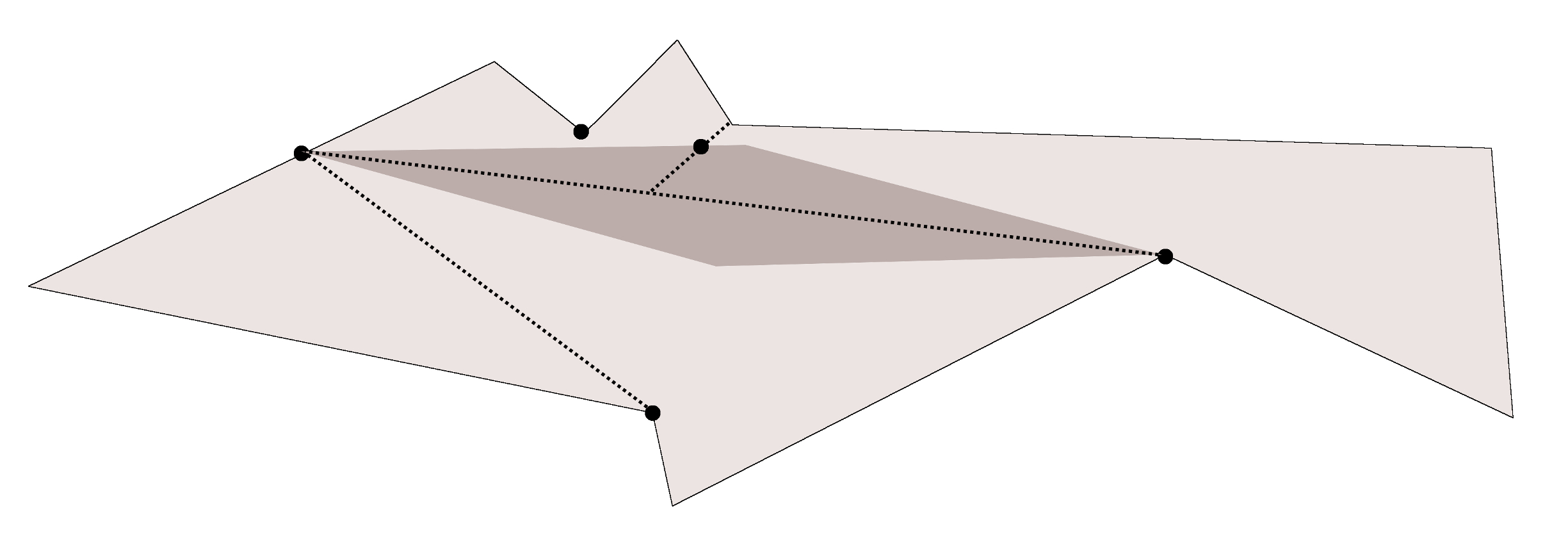}
\put(120,21){{$\bar{v}$}}
\put(227,50){{$v$}}
\put(155,36){{$Q_1$}}
\put(260,66){{$Q_2$}}
\put(28,80){{$w = \bar{w}$}}
\put(111,90){{$\bar{y}$}}
\put(140,76){{$z$}}

\end{overpic}
\caption{A situation with $N'(Q_1)< N'(Q_2)$ is depicted, where $[\bar{v};\bar{w}]$ does not intersect the (open) set ${\rm  cig}([v;w],\eta)$. For the proof of Lemma \ref{lemma: partpol} below we note that $|[w;z]|<|[v;z]|$. Therefore, $\bar{y}$ is `nearer to $w$ than to $v$' and thus $\bar{y} \notin {\rm cig}([v;w],\eta)$ implies $[w;\bar{y}] \cap {\rm cig}([v;w],\eta) = \emptyset$.} \label{semi3}
\end{figure}

For partitions of polygons into semiconvex polygons described in Section \ref{sec: semi2} below we will use this corollary to show that $Q_1$ is   semiconvex. The essential point is that for a segment $[\bar{v};\bar{w}]$ as in \eqref{eq: handing2}  we do not assume the validity of  (SP) and that \eqref{eq: semiXXX} is only required for the vertices contained in $Q_1$. For an illustration of \eqref{eq: handing2}(ii) we refer to Figure \ref{semi3}. 

\Proof We follow the proof of Theorem \ref{lemma: partpol*} and only indicate the necessary changes. Fix $\bar{v},\bar{w}  \in  \mathcal{V}_{Q_1}'\cup \lbrace v,w \rbrace$ such that $[\bar{v};\bar{w}]$ induces a partition $P = R_1 \cup R_2$ fulfilling \eqref{eq: handing2}. Note that one of the sets, say $R_1$, satisfied $R_1 \subset Q_1$ and thus $N'(R_1) \le N'(Q_1)$, $|R_1| \le |Q_1|$. This yields $Q_{\bar{v},\bar{w}} = R_1 \subset Q_1$ (see \eqref{eq: Qvw} and \eqref{eq: handing}). We first suppose 
\begin{align}\label{eq: cap4}
\mathcal{V}_P'\cap {\rm cig}_P([\bar{v};\bar{w}],2\eta) \cap R_1  \subset \lbrace \bar{v},\bar{w} \rbrace
\end{align} 
(compare to \eqref{eq: {eq: cap2}})   and show that under this assumption we have  
 \begin{align}\label{eq: new}
 |[\bar{v};\bar{w}]|  \ge  \frac{\vartheta}{2} \min_{k=1,2} d(R_k).
 \end{align}
 The idea is to proceed as in Step 1 of the previous proof using   \eqref{eq: semiXXX} in place of \eqref{eq: semi}. To this end, we notice that conditions \eqref{eq: semiXXX} and \eqref{eq: cap4} are sufficient to treat the cases (a)-(c). Indeed, as remarked below case (c), case (b) was a purely geometrical argument and in (a),(c) we have only shown \eqref{eq: capXX}. As by \eqref{eq: new} and $Q_{\bar{v},\bar{w}} = R_1$ also condition \eqref{eq: capW} holds (with $\bar{v},\bar{w}$ in place of ${v},{w}$), we derive that in case (a),(c) $[\bar{v};\bar{w}]$ satisfies (WSP). This then  implies \eqref{eq: new} by \eqref{eq: semiXXX}. In cases (a)-(c) we therefore obtain \eqref{eq: new}.  We now show that case (d) never occurs, which concludes the proof of \eqref{eq: new}. 

Suppose case (d) occurs. Then we have that, e.g., $R_1$ is a triangle with vertices $\bar{v}, \bar{w}, p$ such that $\sphericalangle(\bar{v},R_1) \le \alpha_\eta= \frac{1}{2}\arcsin\eta$, $\sphericalangle(\bar{w},R_1) > \pi - 2\alpha_\eta$ and $\bar{w} \notin \mathcal{V}_P'$. The latter immediately implies $\bar{w} = w$ since $\bar{w} \in \mathcal{V}_{Q_1}' \cup \lbrace v,w\rbrace \subset \mathcal{V}_P' \cup \lbrace w\rbrace$. Then  $\bar{v} \neq w$ and $\bar{v} \neq v$ by \eqref{eq: handing2}(i) and \eqref{eq: handing2}(ii) yields that the angle enclosed by the segments $[v;w]$ and $[w;\bar{v}]$ is at least $2\alpha_\eta$ (cf. Figure \ref{semi3}). This, however, contradicts the assumptions $\sphericalangle(\bar{w},R_1) > \pi - 2\alpha_\eta$  and $\bar{w} \notin \mathcal{V}_P'$. Consequently, case (d) never occurs.

\vspace{-0.4cm}

Now we consider an arbitrary segment $[\bar{v};\bar{w}]$ with $\bar{v},\bar{w}  \in  \mathcal{V}_{Q_1}'\cup \lbrace  v, w \rbrace$  which satisfies \eqref{eq: handing2} and induces a partition $P = R_1 \cup R_2$ with $R_1 \subset Q_1$.  Note that  \eqref{eq: handing2} implies
\begin{align}\label{eq: hand3}
w \in \partial R_1 \ \  \ \Rightarrow \ \ \  v \notin \partial R_1 \ \text{ and } \ R_1 \cap {\rm cig}([v;w],\eta) = \emptyset.
\end{align}  
As in Step 2 of the proof of Theorem \ref{lemma: partpol*} we find a chain $(y_1,\ldots,y_{n})$ between $\bar{v}$ and $\bar{w}$ with $y_1 = \bar{v}$, $y_{n} = \bar{w}$ and  $y_i \in \mathcal{V}_{R_1}' \setminus \lbrace \bar{v}, \bar{w} \rbrace$ for $i=2,\ldots,n-1$ such that 
\begin{align}\label{exact}
|[y_i;y_{i+1}]|\le 3|[\bar{v};\bar{w}]|  , \ \ \ \ \  d_P(\bar{v},y_{i})\le \frac{3}{2}|[\bar{v};\bar{w}]| , \ \ \ i=1,\ldots,n-1
\end{align}
and the segments $[y_i;y_{i+1}] \subset  P$ induce a partition satisfying \eqref{eq: cap4} (with $y_i$, $y_{i+1}$ in place of $\bar{v},\bar{w}$). Note that in repeating the argument in \eqref{eq: all j} we only select concave vertices contained in $R_1$, i.e. the essential difference to the previous proof is given by the fact that due to the replacement of \eqref{eq: {eq: cap2}} by \eqref{eq: cap4} we can ensure that each $y_i$, $i=2,\ldots,n-1$, is contained in $\partial R_1 \cap \partial P$, more precisely in $\mathcal{V}'_{R_1}\setminus \lbrace \bar{v}, \bar{w} \rbrace$.  

Note that $y_i \in \mathcal{V}'_{Q_1} \cup \lbrace v, w \rbrace$ for all $i=1,\ldots,n$. Each segment $[y_i;y_{i+1}]$, $i=1,\ldots,n-1$,  induces a  partition $P = R_1^{(i)} \cup R_2^{(i)}$ such that after relabeling $R_1^{(i)} \subset R_1$. By \eqref{eq: hand3} we have that each $[y_i;y_{i+1}]$ satisfies \eqref{eq: handing2} (with $y_i$, $y_{i+1}$ in place of $\bar{v},\bar{w}$). Consequently, as  also 
$$\mathcal{V}_P'\cap {\rm cig}_P([y_i;y_{i+1}],2\eta) \cap R_1^{(i)} \subset \lbrace y_i, y_{i+1} \rbrace
$$
holds by \eqref{eq: cap4} and $R_1^{(i)} \subset R_1$, for each $[y_i;y_{i+1}]$ we may proceed as above and obtain $ |[y_i;y_{i+1}]|  \ge  \frac{\vartheta}{2} \min_{k=1,2} d(R^{(i)}_k)$ by \eqref{eq: new}. This together with \eqref{exact} allows us to proceed exactly as in Step 3 in the proof of Theorem \ref{lemma: partpol*} and we get \eqref{eq: special two} for $\bar{\vartheta}=(3+ 12\vartheta^{-1})^{-1}$. \eop

\subsection{Partition of semiconvex polygons}\label{sec: semi2}

We now show that each polygon can be partitioned into semiconvex polygons.

\begin{theorem}\label{theorem: partpol}
Let $0 < \theta<1$. Then for $\eta >0$ small there exists $\vartheta = \vartheta(\theta,\eta)$ such that for every polygon $P$ there is a partition $P = P_1 \cup \ldots \cup P_N$ into {\rm (SP)}-$\vartheta$-semiconvex polygons $P_1,\ldots, P_N$ such that
\begin{align}\label{eq: partpol0}
\sum^N_{j=1}\mathcal{H}^1(\partial P_j) \le \Big(1 + \frac{2\theta}{1 - \theta}\Big)\mathcal{H}^1(\partial P).
\end{align}
\end{theorem}

Clearly, by Theorem \ref{lemma: partpol*}  the sets $P_1, \ldots, P_N$ are then also $\bar{\vartheta}$-semiconvex for some $\bar{\vartheta}$ small enough. As a preparation we derive a partition $P = Q_1 \cup Q_2$ into two polygons such that $Q_1$ is (SP)-semiconvex. Then Theorem \ref{theorem: partpol}  follows by iterative application. For the proof of Theorem \ref{theorem: partpol} it is essential that (1) the added boundary is small compared to $\mathcal{H}^1(\partial Q_1)$ (see \eqref{eq: partpol21}) and (2) $Q_1$ does not need to be further modified in subsequent iteration steps since hereby the overall added boundary can be controlled (see \eqref{eq: control} below).

\begin{lemma}\label{lemma: partpol}
Let $0 < \theta <1$. Then for $\eta >0$ sufficiently small there is $\tilde{\vartheta} =  \tilde{\vartheta}(\theta,\eta)$ such that for every polygon $P$, which is not an {\rm (SP)}-$\tilde{\vartheta}$-semiconvex polygon, the following holds: We find a segment $[v;w]$ between a concave vertex $v \in \mathcal{V}'_P$ and some $w \in \partial P$ which satisfies {\rm (WSP)} and induces a partition $P = Q_1 \cup Q_2$  such that $Q_1$ is {\rm (SP)}-$\tilde{\vartheta}$-semiconvex and 
\begin{align}\label{eq: partpol21}
|[v;w]| \le \theta\mathcal{H}^1(\partial Q_1 \setminus [v;w]).
\end{align}
Moreover, if $Q_1$ is a triangle we have $\sphericalangle(v,Q_2) <\sphericalangle(v,P) - \frac{1}{2}\arcsin\eta$.
\end{lemma}

\Proof  Let $0 < \theta < 1$ be given and define $\vartheta = \frac{\theta}{2}$. Let $\bar{\vartheta} \le \vartheta$ and $\eta>0$ small as in Corollary \ref{cor: partpol*}.  Define $\tilde{\vartheta} = \bar{\vartheta} \eta (4\eta + 2)^{-1}$. Let $P$ be a non (SP)-$\tilde{\vartheta}$-semiconvex polygon.  \\

\noindent  \smallskip
\textit{Step 1: Choice of $[v;w]$}

\noindent As $P$ is not (SP)-$\tilde{\vartheta}$-semiconvex, there is at least one segment $[v; w]$, between a concave vertex $v  \in \mathcal{V}'_P$ and some $w \in \partial P$  which satisfies (SP) (and thus also (WSP))   and induces a partition $P = Q_1 \cup Q_2$ with $|[v;w]| < \tilde{\vartheta} d(Q_k) \le \vartheta d(Q_k)$ for $k=1,2$. In the following we label the sets such that we always have $N'(Q_1) \le N'(Q_2)$   (recall \eqref{eq: Qvw}).  Choose (possibly not uniquely) a pair $v,w$ satisfying (WSP) and
\begin{align}\label{eq: partpol222}
|[v;w]| < \vartheta \min_{k=1,2} d(Q_k) = \frac{\theta}{2} \min_{k=1,2} d(Q_k)
\end{align}
in such a way that  $N'(Q_1)$ is minimized among all pairs satisfying (WSP)  and \eqref{eq: partpol222}.  If $N'(Q_1) = N'(Q_2)$, we may suppose $|Q_1| \le |Q_2|$ after possible relabeling. After a small perturbation of the point $w$ we may assume that $|Q_1|<|Q_2|$ and  (WSP), \eqref{eq: partpol222}  are still satisfied. (Recall here that ${\rm cig}_P([v;w],\eta)$ is closed.)  Moreover, we note that  $v,w$ can be selected such that  
\begin{align}\label{eq: perturb}
&[w_*;w] \subset \partial P   \ \text{ and } \ |[w_*;w]| \le \frac{1}{2}|[v;w]|
\end{align}
for all $w_* \in \partial Q_1$ with the property that $[v;w_*]$ induces $P = Q_1^* \cup Q_2^*$  satisfying
\begin{align}\label{eq: perturbXXX}
 \  \text{ (WSP), } \  N'(Q_1^*) = N'(Q_1), \   \text{ and } \  |[v;w_*]| < \vartheta \min_{k=1,2} d(Q^*_k). 
\end{align}
In fact,  if \eqref{eq: perturb} is violated for some $w_*$ which satisfies \eqref{eq: perturbXXX},  we can replace the pair $v,w$ by the pair $v,w_*$ in the above choice (accordingly, we replace $Q_1$ by the smaller set $Q_1^*$). Possibly  repeating this procedure at most $\mathcal{V}_P + \lfloor \frac{2\mathcal{H}^1(\partial P)}{d_v}\rfloor $ times, where $d_v :=\inf\lbrace |[v;w']|: [v;w'] \text{ induces a partition of $P$}\rbrace>0 $,  we obtain a (not relabeled) pair $v,w$ such that \eqref{eq: perturb} holds for all $w_*$ satisfying \eqref{eq: perturbXXX}. 

Choose $p,p' \in \partial Q_1$ with  $d(Q_1) = \dist_{Q_1}(p,p')$. Since  $d(Q_1)\le \dist_{Q_1}(p,v) + \dist_{Q_1}(v,p')$, we can without restriction assume that $\dist_{Q_1}(p,v) \ge \frac{1}{2}d(Q_1)$ and thus by \eqref{eq: partpol222} we get 
$$\dist_{Q_1}(p,v) \ge \theta^{-1}|[v;w]| \ \ \  \text{      and      } \ \ \  \dist_{Q_1}(p,w) \ge \big( \theta^{-1}-1\big)|[v;w]|.$$
Consequently, $\mathcal{H}^1(\partial Q_1 \setminus [v;w]) \ge \dist_{Q_1}(p,v) + \dist_{Q_1}(p,w)$ and in view of $\theta<1$ a short calculation yields  \eqref{eq: partpol21}.  The additional assertion after \eqref{eq: partpol21} follows directly from the fact that $[v;w]$ satisfies (WSP),  particularly \eqref{eq: capXX}, where we use $\sphericalangle(v,Q_2) = \sphericalangle(v,P) - \sphericalangle(v,Q_1)$.  It remains to show that $Q_1$ is (SP)-$\tilde{\vartheta}$-semiconvex.

\vspace{-0.3cm}
\begin{figure}[H]
\centering

\begin{overpic}[width=0.76\linewidth,clip]{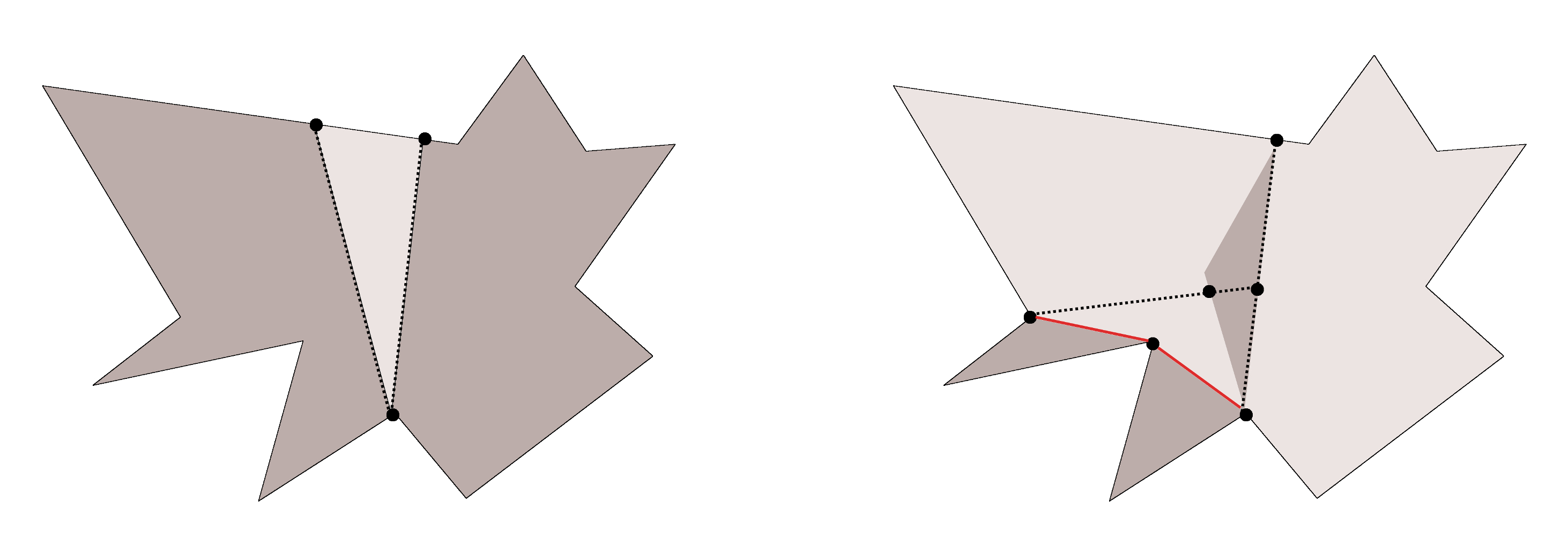}
\put(248.5,25.5){{$v=y_1$}}
\put(247,84){{$w$}}
\put(251,48){{$\bar{w}$}}
\put(230,52){{$z$}}
\put(170,44){{$\bar{v}=y_3$}}
\put(213,32){{$y_2$}}
\put(190,20){{$Q^{(2)}_1$}}
\put(228,5){{$Q^{(1)}_1$}}
\put(188,75){{$R_2$}}
\put(194,30) {\line(0,1){7}}
\put(234,15) {\line(0,1){7}}

\put(81,24){{$v = \bar{v}$}}
\put(84,83){{$w$}}
\put(60,87){{$\bar{w}$}}
\put(67,68){{$R_2$}}
\put(25,69){{$R_1 = T_1$}}
\put(105,68){{$Q_2$}}

\end{overpic}
\caption{On the left case (a) is depicted, where $R_2 = Q_1 \cap T_2$ and $Q_1 = T_1 \cup R_2$. On the right case (b) with $|[v;z]| < |[w;z]|$ is illustrated.} \label{semi4}
\end{figure}
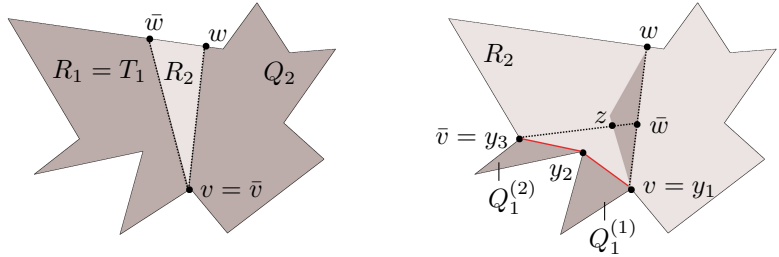
\noindent \smallskip
\textit{Step 2: Semiconvexity of $Q_1$}

\noindent   As a preparation we show that the assumptions of Corollary \ref{cor: partpol*} are satisfied. Consider a pair $v',w' \in \mathcal{V}'_{Q_1} \cup \lbrace v, w \rbrace$ such that the segment $[v'; w'] \neq [v;w]$    induces a partition $P = Q_1' \cup Q_2'$ satisfying (WSP) with $Q_1' \subset Q_1$.  As either (i) $v'\neq v$ and  $w' \neq v$ or (ii) up to  relabeling $v'= v$, $w' \neq w$ with $w' \in \mathcal{V}'_P$, we get $\# \lbrace u \in \mathcal{V}_P' \setminus \lbrace v',w'\rbrace: u \in \partial Q'_1 \rbrace < \# \lbrace u \in \mathcal{V}_P' \setminus \lbrace v,w\rbrace: u \in \partial Q_1 \rbrace$. Thus,  $N'(Q_1') < N'(Q_1)$ (see before \eqref{eq: Qvw}).   As (again up to relabeling of the points) $v' \in \mathcal{V}_P'$ and $w' \in \partial P$, we observe that $N'(Q'_1) < N'(Q_1)$ together with  the choice of $N'(Q_1)$   and \eqref{eq: partpol222} implies
\begin{align}\label{eq: Sstrich}
|[v';w']| \ge   \vartheta\min_{k=1,2} d(Q'_k)
\end{align} 
and thus \eqref{eq: semiXXX} holds. Moreover, we recall that  \eqref{eq: handing} is satisfied by the choice of $Q_1$ (see before \eqref{eq: perturb}).

We now show that $Q_1$ is (SP)-$\tilde{\vartheta}$-semiconvex. To this end, consider a pair $\bar{v} \in \mathcal{V}'_{Q_1}$ and $\bar{w} \in \partial Q_1$ such that $[\bar{v}; \bar{w}]$ induces a partition $Q_1 = R_1 \cup R_2$ satisfying (SP). In particular, $\mathcal{V}_{Q_1}' \cap {\rm cig}_{Q_1}([\bar{v};\bar{w}],\eta) \subset \lbrace \bar{v}, \bar{w} \rbrace$ by \eqref{eq: cap}. We distinguish the cases (a) $\bar{w} \notin [v;w] \setminus \lbrace v,w\rbrace$ and (b) $\bar{w} \in [v;w]\setminus \lbrace v,w\rbrace$.   \smallskip \\
(a) Assume $\bar{w} \notin [v;w] \setminus \lbrace v,w\rbrace$. Clearly, $[\bar{v}; \bar{w}]$ induces also a partition $P=T_1 \cup T_2$, where we label the sets such that $R_1 =  T_1 \subset Q_1$. Since \eqref{eq: handing} holds, we have $Q_{\bar{v},\bar{w}} = T_1$ and therefore
$$\mathcal{V}_{P}' \cap {\rm cig}_{P}([\bar{v};\bar{w}],\eta) \cap Q_{\bar{v},\bar{w}} \subset \mathcal{V}_{Q_1}' \cap {\rm cig}_{Q_1}([\bar{v};\bar{w}],\eta) \subset \lbrace \bar{v}, \bar{w} \rbrace.$$
Consequently,  $[\bar{v}; \bar{w}]$ satisfies (WSP) with respect to the partition $P=T_1 \cup T_2$.  \smallskip \\
(a1) Assume $\bar{w} \notin [v;w] \setminus \lbrace v,w\rbrace$ and $N'(T_1) < N'(Q_1)$. Since $\bar{w} \in \partial P$ and $\bar{v} \in \mathcal{V}_P'$, we may proceed as in \eqref{eq: Sstrich}, particularly using \eqref{eq: partpol222}, to find 
$$|[\bar{v};\bar{w}]| \ge   \vartheta\min_{k=1,2} d(T_k)\ge  \vartheta\min_{k=1,2} d(R_k).$$
(a2) Now suppose $\bar{w} \notin [v;w]\setminus \lbrace v,w\rbrace$ and $N'(T_1) = N'(Q_1)$.  Since $T_1 \subset Q_1$, this is only possible if $\bar{v} = v$.  If $|[\bar{v};\bar{w}]| \ge   \vartheta\min_{k=1,2} d(T_k)$, we proceed as in (a). Otherwise,  by   \eqref{eq: perturb} we obtain $[w;\bar{w}] \subset \partial P$ and $|[w;\bar{w}]| \le \frac{1}{2}|[v;w]|$. Consequently, we get  $|[\bar{v};\bar{w}]|\ge \frac{1}{2}|[v;w]|$. As $R_2 = T_2 \cap Q_1$ is a triangle with vertices $v,w,\bar{w}$ (cf. \eqref{eq: perturb} and Figure \ref{semi4}), we deduce
\begin{align*}
\min_{k=1,2}d(R_k) \le d(R_2) = d(T_2 \cap Q_1)\le |[v;w]| + |[w;\bar{w}]| \le \frac{3}{2}|[v;w]| \le 3|[\bar{v};\bar{w}]|.
\end{align*}
Thus, in both cases (a1), (a2) condition \eqref{eq: semi} holds since $\tilde{\vartheta} \le \min\lbrace \vartheta, \frac{1}{3}\rbrace$  for $\eta$ small. \smallskip \\
(b) It now remains to treat the case $\bar{w} \in [v;w]\setminus \lbrace v,w\rbrace$.   We label the sets such that $v \in R_1$ and $w \in R_2$. As $[v;w]$ satisfies (WSP) and $Q_{v,w} = Q_1$ by \eqref{eq: handing},  $\bar{v} \notin {\rm cig}_P([v;w],\eta)$ holds. Thus,  $\bar{v} \notin {\rm cig}([v;w],\eta)$ since $[\bar{v};\bar{w}] \subset P$. Let $z$ be the intersection point of  $\partial {\rm cig}([v;w],\eta)$ with $[\bar{v};\bar{w}]$ (see Figure \ref{semi4}).  We first treat the case $|[v;z]| \le |[w;z]|$ and present the necessary adaptions for the other case at the end of the proof. Recall that the goal is to show $\tilde{\vartheta} \min_{k=1,2} d(R_k) \le |[\bar{v};\bar{w}]|$. Now $\bar{v} \notin {\rm cig}([v;w],\eta)$ and $|[v;z]| \le |[w;z]|$ imply (cf. Figure \ref{semi4})
\begin{align}\label{eq: shortier}
|[\bar{v};\bar{w}]| \ge |[\bar{w};z]| \ge \eta|[v;z]|.
\end{align}
Choose the (unique) chain $(y_1 = v, y_2, \ldots, y_n = \bar{v})$ with $y_i \in \mathcal{V}_P' \cap \mathcal{V}_{R_1}$ for $i=2,\ldots,n-1$ such that 
$d_P(v,\bar{v}) = \sum_{i=1}^{n-1} |[y_i;y_{i+1}]|$ and $[y_i;y_{i+1}]$ induces a partition $P = Q^{(i)}_1 \cup Q^{(i)}_2$, where the sets are labeled such that $w \in Q^{(i)}_2$. Observe that $Q^{(i)}_1 \subset R_1$ as $[\bar{v};\bar{w}]$ induces a partition of $Q_1 = R_1 \cup R_2$. Then by \eqref{eq: shortier} and  $[v;\bar{w}] \cup [\bar{v};\bar{w}] \subset P$  
\begin{align}\label{eq: lengthi}
d_P(v,\bar{v}) \le |[v;\bar{w}]| + |[\bar{v};\bar{w}]| \le |[v;z]| + |[z;\bar{w}]| + |[\bar{v};\bar{w}]|  \le (2+ \eta^{-1})|[\bar{v};\bar{w}]|. 
\end{align}
Note that each $[y_i;y_{i+1}]$ satisfies \eqref{eq: handing2} as $w \notin \partial R_1$. Then using Corollary \ref{cor: partpol*}, in particular \eqref{eq: special two}, we find for all $i=1,\ldots,n-1$
\begin{align}\label{eq: lengthi2} 
|[y_i;y_{i+1}]|  \ge \bar{\vartheta} \min_{k=1,2} d(Q^{(i)}_k).
\end{align}
First assume there was some $i$ such that $d(Q^{(i)}_2) \le d(Q^{(i)}_1)$. Then we calculate using $Q^{(i)}_2 \supset R_2$ and \eqref{eq: lengthi} 
$$d(R_2) \le  d(Q_2^{(i)}) \le \bar{\vartheta}^{-1}|[y_i;y_{i+1}]|  \le \bar{\vartheta}^{-1}d_P(v,\bar{v}) \le \bar{\vartheta}^{-1}(2+ \eta^{-1})|[\bar{v};\bar{w}]|.$$
Otherwise, we find by \eqref{eq: lengthi}, \eqref{eq: lengthi2}   and $d_P(v,\bar{v}) = \sum_{i=1}^{n-1} |[y_i;y_{i+1}]|$  (cf. Figure \ref{semi4})
\begin{align*}
d(R_1) &\le 2\max_{p \in R_1} d_P(v,p) \le 2\max\lbrace |[v;\bar{w}]| + |[\bar{w};\bar{v}]|, \max_{i=1,\ldots,{n-1}} (d_P(v,y_i) +  d(Q^{(i)}_1)) \rbrace \\
& \le  2\max\lbrace|[v;\bar{w}]| + |[\bar{w};\bar{v}]|,\bar{\vartheta}^{-1}d_P(v,\bar{v})\rbrace \le 2\bar{\vartheta}^{-1}(2+ \eta^{-1})|[\bar{v};\bar{w}]|.
\end{align*}
Collecting the last two estimates and recalling  $\tilde{\vartheta} = \bar{\vartheta} \eta (4\eta + 2)^{-1}$ we  get $\tilde{\vartheta} \min_{k=1,2} d(R_k) \le |[\bar{v};\bar{w}]|$, as desired. 

It remains to treat the case $|[v;z]| > |[w;z]|$. We may proceed as before with $w$ in place of $v$ with the only difference that, due to the fact that the chain $(y_1=w,\ldots,y_n=\bar{v}) \subset \mathcal{V}_{R_2}$ contains $w$,  for the application of \eqref{eq: special two}  we have to check that \eqref{eq: handing2}(ii) holds for $[y_1;y_2]$. Indeed, since $[v;w]$ satisfies (WSP), the fact that $Q_1 = Q_{v,w}$  implies $y_2 \notin {\rm cig}_P([v;w],\eta)$ and then $y_2 \notin {\rm cig}([v;w],\eta)$ since $[y_1;y_2] \subset P$. Then $|[v;z]| > |[w;z]|$ together with $y_2 \in R_2$ yield $[y_1; y_2] \cap {\rm cig}([v;w],\eta) = \emptyset$ (cf. illustration of $\bar{y}$ in Figure \ref{semi3}).  \eop

Now we can give the proof of Theorem \ref{theorem: partpol}.

\smallskip
\noindent\emph{Proof of Theorem \ref{theorem: partpol}.}   We construct the partition inductively. Assume $P_1,\ldots,P_{n}$ have been constructed and set $R_n = \overline{P \setminus \bigcup^{n}_{j=1} P_j}$. (For $n=0$ we set $R_0 = P$.) Moreover, suppose that
\begin{align}\label{eq: partpol1}
\mathcal{H}^1\Big(\partial P_j \setminus \bigcup\nolimits^{j-1}_{i=0} \partial P_i\Big) \le \theta \mathcal{H}^1 \Big(\partial P_j  \cap \bigcup\nolimits^{j-1}_{i=0} \partial P_i \Big)
\end{align}
for $j=1,\ldots,n$, where $P_0 := P$. If $R_n$ is (SP)-semiconvex, we set $P_{n+1} = R_n$ and stop. Otherwise, by Lemma \ref{lemma: partpol} we find a partition $R_n = P_{n+1} \cup R_{n+1}$ such that $P_{n+1}$ is (SP)-semiconvex and $R_{n+1} = \overline{R_n \setminus P_{n+1}} =  \overline{P \setminus \bigcup^{n+1}_{j=1} P_j}$. Furthermore, we obtain by \eqref{eq: partpol21}
\begin{align*}
\mathcal{H}^1\Big(\partial P_{n+1} \setminus \bigcup\nolimits^{n}_{i=0} \partial P_i\Big) &=  \mathcal{H}^1(P_{n+1} \cap R_{n+1})  \le \theta\mathcal{H}^1(\partial  P_{n+1} \setminus (P_{n+1} \cap R_{n+1})) \\ 
&= \theta \mathcal{H}^1\Big(\partial P_{n+1}  \cap \bigcup\nolimits^{n}_{i=0} \partial P_i\Big),
\end{align*}
which gives \eqref{eq: partpol1} for $j=n+1$.

Recall that in each step   the number of vertices of the remaining polygon decreases (namely if $P_{n+1}$ is not a triangle) or the angle of a concave vertex in the remaining polygon decreases by at least $\frac{1}{2}\arcsin\eta$ (if $P_{n+1}$ is a triangle). Thus, there is some $N \in \N$ such that the polygon $P_N := R_{N-1}$ is (SP)-semiconvex since for large $n \in \N$ the polygon $R_{n-1}$  is eventually convex and thus also (SP)-semiconvex.  It remains to show \eqref{eq: partpol0}. First, we note
 \begin{align}\label{eq:nn2}
\mathcal{H}^1(\partial P_{n}  \cap \partial P_i) = \mathcal{H}^1(\partial R_{n-1}  \cap \partial P_{i}) - \mathcal{H}^1(\partial R_{n}  \cap \partial P_i)
 \end{align}
 for $1 \le n \le N$ and $0 \le i \le n-1$, where we set $R_N := \emptyset$ and $R_0 = P$. Moreover, by \eqref{eq: partpol1} we get  for $2 \le n \le N$
 \begin{align}\label{eq:nn1}
 \mathcal{H}^1(\partial R_{n-1} \cap \partial P_{n-1}) = \mathcal{H}^1\Big(\partial P_{n-1} \setminus \bigcup\nolimits_{i=0}^{n-2} \partial P_i\Big) \le  \theta\mathcal{H}^1\Big(\partial P_{n-1}  \cap \bigcup\nolimits^{n-2}_{i=0} \partial P_i\Big).
 \end{align}
Then by \eqref{eq:nn2}-\eqref{eq:nn1} we obtain $\mathcal{H}^1(\partial P_{1}  \cap \partial P) = \mathcal{H}^1(\partial P) - \mathcal{H}^1(\partial R_1 \cap \partial P)$ and for $2 \le n \le N$
\begin{align*}
\mathcal{H}^1\Big(\partial P_{n}  \cap \bigcup\nolimits^{n-1}_{i=0} \partial P_i\Big) &= \sum^{n-1}_{i=0}\mathcal{H}^1\Big(\partial P_{n}  \cap  \partial P_i\Big)  \le \theta\mathcal{H}^1\Big(\partial P_{n-1}  \cap \bigcup\nolimits^{n-2}_{i=0} \partial P_i\Big) \\ & \ \ \ + \sum^{n-2}_{i=0} \mathcal{H}^1(\partial R_{n-1}  \cap \partial P_{i}) - \sum^{n-1}_{i=0} \mathcal{H}^1(\partial R_{n}  \cap \partial P_i).
\end{align*}
By summation and an index shift we  derive
\begin{align*}
\sum^{N}_{n=1}  \mathcal{H}^1 \Big(\partial P_n  \cap \bigcup\nolimits^{n-1}_{i=0} \partial P_i \Big)    & \le  \mathcal{H}^1(\partial P_{1}  \cap \partial P)  + \theta\sum^{N-1}_{n=1}  \mathcal{H}^1 \Big(\partial P_n  \cap \bigcup\nolimits^{n-1}_{i=0} \partial P_i \Big)\\
& \ \ \  +  \sum_{n=1}^{N-1} \sum_{i=0}^{n-1} \mathcal{H}^1(  \partial R_n  \cap \partial P_i) - \sum_{n=2}^{N} \sum_{i=0}^{n-1} \mathcal{H}^1(  \partial R_n  \cap \partial P_i) \\
& \le \mathcal{H}^1(\partial P)  + \theta\sum^{N}_{n=1}  \mathcal{H}^1 \Big(\partial P_n  \cap \bigcup\nolimits^{n-1}_{i=0} \partial P_i \Big),
\end{align*}
where we used that $\partial R_N = \emptyset$ and $\mathcal{H}^1(\partial R_1 \cap \partial P)  + \mathcal{H}^1(\partial P_{1}  \cap \partial P)  = \mathcal{H}^1(\partial P) $. 
This yields $\sum^{N}_{n=1}  \mathcal{H}^1 \Big(\partial P_n  \cap \bigcup\nolimits^{n-1}_{i=0} \partial P_i \Big)  \le ({1- \theta})^{-1}\mathcal{H}^1( \partial P)$. Together with \eqref{eq: partpol1} and the fact that every $x \in \bigcup^N_{n=1}\partial P_n \setminus \partial P$ is contained in the boundary of exactly two sets,  we conclude
\begin{align}\label{eq: control}
\begin{split}
\sum^N_{n=1} \mathcal{H}^1(\partial P_n) & = \mathcal{H}^1(\partial P) +  2\mathcal{H}^1\Big(\bigcup\nolimits^N_{n=1} \partial P_n \setminus \partial P\Big) \\ & = \mathcal{H}^1(\partial P) + 2\sum^N_{n=1} \,   
\mathcal{H}^1\Big(\partial P_n \setminus \bigcup\nolimits^{n-1}_{i=0} \partial P_i\Big)  \\
& \le \mathcal{H}^1(\partial P) +  2\theta\sum^N_{n=1}  \mathcal{H}^1 \Big(\partial P_n  \cap \bigcup\nolimits^{n-1}_{i=0} \partial P_i \Big) \le \Big(1 + \frac{2\theta}{1 - \theta}\Big) \mathcal{H}^1(\partial P).
\end{split}
\end{align}
 \eop

Later in Section \ref{sec: proofi} for the proof of Theorem \ref{th: main result smooth} we will   need the following observations. 
\begin{rem}\label{rem: vertices1}

{\normalfont

(i) Recall that by construction the partition  in Theorem \ref{theorem: partpol} arises from $P$ by introducing a finite number of segments. As by this procedure no additional concave vertices are introduced, we find $v \in \partial P$ for all $v \in \bigcup_{j=1}^N \mathcal{V}'_{P_j}$.   

(ii) By a slight modification of the segments $[v;w]$ introduced in Lemma \ref{lemma: partpol}  (cf. Remark below \eqref{eq: partpol222}) we can always ensure that the  segments $[v_i;w_i] = P_i \cap P_{i+1}$ have the property that the points $w_i$ are not  vertices of $P$ and are pairwise distinct.

(iii) The partition can be chosen with the following additional property: if two convex polygons $P^1,P^2 \subset (P_j)_{j=1}^N$ share some $v \in \mathcal{V}_P'$ with $\sphericalangle(v,P^i)\le \frac{\pi}{4}$ for $i=1,2$, then   $\mathcal{H}^1(\partial P^1 \cap \partial P^2) = 0$. Indeed, otherwise we find $w \in \mathcal{V}_{P^1} \cap \mathcal{V}_{P^2}$ such that $\partial P^1 \cap \partial P^2 =[v;w]$. Then with $P_* = P^1 \cup P^2$ we have $\sphericalangle(v,P_*) \le \frac{\pi}{2}$ and $\sphericalangle(w,P_*) = {\pi}$ by (ii). Consequently, $P_*$ is a convex polygon and we can replace in the partition $P^1$, $P^2$ by $P_*$.
 }
\end{rem} 

We close this section with a further criterion for the partition of a semiconvex polygon.

\begin{lemma}\label{lemma: part-seg}
Let $0 <\alpha, \vartheta  <1$. Then there is $\bar{\vartheta}=\bar{\vartheta}(\alpha,\vartheta)>0$ such that for all $\vartheta $-semiconvex polygons $P$ the following holds: If there is a segment $[u_1; u_2]$  inducing a partition $P = P_1 \cup P_2$ such that for each concave vertex $v \in \mathcal{V}'_{P_1}$ one has that
$$\max_{k=1,2} \sphericalangle(u_k,\triangle_v) \ge \alpha,$$
where $\triangle_v$ is the triangle with vertices $v,u_1,u_2$, then $P_1$ is $\bar{\vartheta}$-semiconvex.
\end{lemma}

\Proof Let $P$ and the partition $P = P_1 \cup P_2$ with the above properties be given. To see that $P_1$ is $\bar{\vartheta}$-semiconvex for some $\bar{\vartheta} \le \vartheta$ to be specified below, it suffices to show that for each segment $[v;w]$ between a concave vertex $v \in \mathcal{V}'_{P_1}$ and some $w \in [u_1; u_2]$, which  induces a partition $P_1 = Q_1 \cup Q_2$, one has
\begin{align}\label{eq: semi2}
|[v;w]| \ge \bar{\vartheta} \min_{k=1,2} d(Q_k).
\end{align} 
Indeed, for $w \in \partial P_1 \setminus [u_1; u_2]$ the property follows directly from the fact that $P$ is $\vartheta$-semiconvex. Without restriction we assume $\sphericalangle(u_1,\triangle_v) \ge \sphericalangle(u_2,\triangle_v)$ and label the sets such that $u_1 \in Q_1$.  Similarly as in the proof of  Lemma \ref{lemma: partpol} we choose the unique chain $(y_1 = v, y_2, \ldots, y_n = u_1)$ with $y_i \in \mathcal{V}_P'$ for $i=2,\ldots,n-1$ such that 
$d_P(v,u_1) = \sum_{i=1}^{n-1} |[y_i;y_{i+1}]|$ and $[y_i;y_{i+1}]$ induces a partition $P = Q^{(i)}_1 \cup Q^{(i)}_2$,  where the sets are labeled such that $u_{2} \in Q^{(i)}_2$.   Since $P$ is $\vartheta$-semiconvex, we get by \eqref{eq: semi}
\begin{align}\label{eq: semi5}
|[y_i;y_{i+1}]|  \ge  \vartheta  \min_{k=1,2} d(Q^{(i)}_k)
\end{align}
for $i=1,\ldots,n-1$. Observe that $Q^{(i)}_1 \subset Q_1$ as $[v;w]$ induces a partition of $P_1$.  Using $\sphericalangle(u_1,\triangle_v) \ge \alpha$ and  the cosine formula we find by an elementary computation 
$$
|[v;w]| \ge \sqrt{|[v;u_1]|^2 + |[u_1;w]|^2 - 2|[v;u_1]||[u_1,w]|\cos\alpha} \ge C_\alpha(|[v;u_1]| + |[u_1;w]|)
$$
for  $C_\alpha>0$ small depending only on $\alpha$. Using that $[v;w] \cup [w;u_1] \subset P$ we then derive
\begin{align}\label{eq: lengthiXX}
d_P(v,u_1) \le |[v;w]| +|[w;u_1]| \le |[v;w]| + C_\alpha^{-1}|[v;w]|  = (1 +C_\alpha^{-1})|[v;w]|. 
\end{align}
We now proceed as in the proof of Lemma \ref{lemma: partpol}. First assume there is some $i$ such that $d(Q^{(i)}_2) \le d(Q^{(i)}_1)$. Then we calculate using $Q^{(i)}_2 \supset Q_2$, \eqref{eq: semi5} and \eqref{eq: lengthiXX}  
\begin{align*}
d(Q_2) \le d(Q^{(i)}_2) \le \vartheta^{-1} |[y_i;y_{i+1}]| \le \vartheta^{-1}  d_P(v,u_1) \le (1 +C_\alpha^{-1})\vartheta^{-1} |[v;w]|. 
\end{align*}
 Otherwise, we find again by \eqref{eq: semi5} and \eqref{eq: lengthiXX} 
\begin{align*}
d(Q_1) &\le 2\max_{p \in Q_1} d_P(v,p) \le 2\max\lbrace |[v;w]| + |[w;u_1]|, \max_{i=1,\ldots,{n-1}} (d_P(v,y_{i}) +  d(Q^{(i)}_1)) \rbrace \notag \\
& \le  2\max\lbrace|[v;w]| + |[w;u_1]|,\vartheta^{-1} d_P(v,u_1)\rbrace \le 2\vartheta^{-1}(1+C_\alpha^{-1})|[v;w]|.
\end{align*}
Consequently, \eqref{eq: semi2} holds for $\bar{\vartheta} = \vartheta C_\alpha (2+2C_\alpha)^{-1}$ and thus  $P_1$ is $\bar{\vartheta}$-semiconvex. \eop

\section{Semiconvex and rotund polygons}\label{sec: circ}

In the section we  introduce a further subclass of polygons. 

\begin{definition}\label{def: qussc}
Let $\omega>0$. We say a polygon $P$ is $\omega$-rotund if there is a ball $B(x,r) \subset P$ with $x \in P$ and $r \ge \omega d(P)$. 
\end{definition} 

Similarly as before, we drop the parameter $\omega$ if no confusion arises. This property together with the semiconvexity will be the main ingredient to show that polygons may be partitioned into John domains with controllable John constant. In Section \ref{sec: quasi} we study the relation between the notions introduced in Definition \ref{def: semi} and Definition \ref{def: qussc}. In Section \ref{sec: part-semrot} we then show that semiconvex polygons can be partitioned into semiconvex and rotund polygons.

\subsection{Properties of semiconvex and rotund polygons}\label{sec: quasi}

To avoid confusion with further subscripts we will from now on denote by $x \e_j$ the $j$-th component of points $x \in \R^2$. For sets $A \subset \R^2$ and $R \in SO(2)$ we let $|A|_{\Pi,R} = \sup_{x,y \in A} |(x-y)R\e_1|$. We will also use the notation  $|A|_{\Pi,j} = \sup_{x,y \in A} |(x-y)\e_j|$  for $j=1,2$. By ${\rm int}(A)$ we denote the interior of a set. Recall also the notions introduced in Section \ref{sec: prep}.   We begin with a simple property of convex polygons.

 \begin{lemma}\label{lemma: convpart2}
Every convex polygon $P$ contains a ball with radius $$\frac{1}{4} \min_{R \in SO(2)} |P|_{\Pi,R}.$$ 
\end{lemma}

\Proof By \cite{Lassak:1993} we find that for each convex polygon $P$ there is a rectangle $S$ and a homothetic copy $S'$ of $S$ such that $S \subset P \subset S'$ and the positive homothety ratio is at most $2$. As $P\subset S'$, both rectangle sides of $S'$ are larger than $\min_{R \in SO(2)} |P|_{\Pi,R}$ and thus each rectangle side of $S$ is larger than $\frac{1}{2}\min_{R \in SO(2)} |P|_{\Pi,R}$.  \eop

We now show that the  intrinsic diameter of semiconvex polygons $P$ can be controlled in terms of  $|P|_{\Pi,R}$.

\begin{lemma}\label{lemma: diam}
Let $0<\vartheta <1$ and let $P$ be a $\vartheta$-semiconvex polygon. Then 
$$\vartheta d(P) \le 2\max_{R \in SO(2)}|P|_{\Pi,R}.$$
\end{lemma}

\Proof If $P$ is convex, the assertion is clear. Otherwise, choose $p_1,p_2 \in P$ with $d_P(p_1,p_2) = d(P)$ and let $\gamma:[0,l(\gamma)] \to P$ be a piecewise affine curve between $p_1,p_2$, parametrized by arc length, with $l(\gamma) = d(P)$. 

Since $P$ is not convex,  there is some $v \in \mathcal{V}_P'$ such that $[v;\gamma(\frac{l(\gamma)}{2})] \subset P$. (Possibly we have to take $v= \gamma(\frac{l(\gamma)}{2})$.) Then we can choose $w \in \partial P$ such that $\gamma(\frac{l(\gamma)}{2}) \in [v;w]$ and $[v;w]$ induces a partition $P=Q_1 \cup Q_2$ according to Definition \ref{def: induce}. The choice of $\gamma$ implies $\min_{k=1,2}d(Q_k) \ge \frac{1}{2} d(P)$ and thus we conclude, using that $P$ is $\vartheta$-semiconvex  
$$\frac{\vartheta}{2} d(P)\le \vartheta\min_{k=1,2} d(Q_k) \le |[v;w]| \le \max_{R \in SO(2)}|P|_{\Pi,R}.$$\eop

We now formulate the first main result of this section stating that semiconvex polygons are rotund if the lengths of shortest and longest extend are comparable. 

\begin{theorem}\label{th: only two}
Let $0 < \vartheta, \lambda <1$. Then there is an $\omega=\omega(\vartheta,\lambda)>0$ such that all $\vartheta$-semiconvex polygons $P$ with 
$$\min_{R \in SO(2)} |P|_{\Pi,R} \ge \lambda \max_{R \in SO(2)} |P|_{\Pi,R} $$
are $\omega$-rotund.
\end{theorem}

Whereas the statement is straightforward for convex polygons by  Lemma \ref{lemma: convpart2}, the argument for nonconvex polygons relies on the observation that concave vertices are `not too close to opposite parts of the boundary'  due to condition \eqref{eq: semi}.

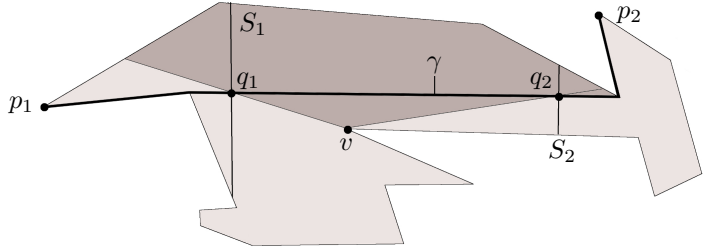
\begin{figure}[H]
\centering

\begin{overpic}[width=0.76\linewidth,clip]{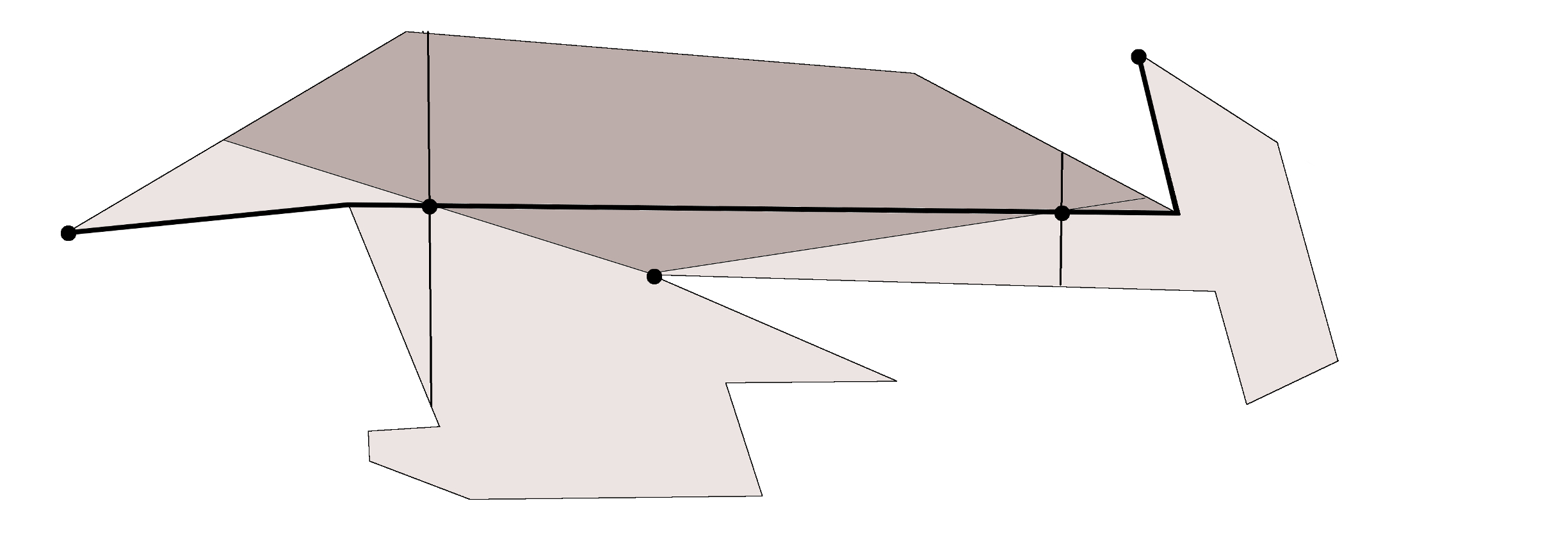}
\put(230,97){{$p_2$}}
\put(204,44){{$S_2$}}
\put(197,72){{$q_2$}}
\put(161,68) {\line(0,1){7}}
\put(158,78) {$\gamma$}
\put(125,47) {$v$}
\put(0,63){{$p_1$}}
\put(87,92){{$S_1$}}
\put(86,72){{$q_1$}}
\end{overpic}
\caption{We sketched case (a2), where $T$ is contained in the dark gray set.} \label{rot1}
\end{figure}

\Proof  Choose $p_1,p_2 \in P$ with $d_P(p_1,p_2) = d(P)$ and let $\gamma:[0,l(\gamma)] \to P$ be a piecewise affine curve between $p_1,p_2$, parametrized by arc length, with $l(\gamma) = d(P)$. As noticed in Section \ref{sec: prep}, recall that the endpoints of each segment of $\gamma$ are contained in $\mathcal{V}_P' \cup  \lbrace p_1,p_2 \rbrace$.    Define $\delta = \frac{1}{14}\vartheta \lambda$ and set for shorthand $q_1 = \gamma(\delta)$ and $q_2 = \gamma(1-\delta)$.  We distinguish two cases: \smallskip \\
(a) First assume $ \gamma([\delta,1-\delta]) =  [q_1;q_2]$ is a segment with $q_1, q_2 \notin \partial P$ and suppose that after translation and rotation we have $q_1 =  (t_1,0)$,  $q_2 = (t_2,0)$ with $t_1 < t_2$.   For $k=1,2$ denote by $S_k$ the  connected component of $(\lbrace t_k \rbrace \times \R) \cap {\rm int}(P)$ containing $q_k$. The segments $\overline{S_1},\overline{S_2}$ induce a partition $P = P_1 \cup P' \cup P_2$ of $P$ with $P_k \cap P' = \overline{S_k}$ for $k=1,2$ (cf. Figure \ref{rot1}). First, by the  fact that $l(\gamma) = d(P)$ and $2 \delta \le \frac{1}{2}$ we get  
\begin{align}\label{eq:n1}
 d(P') \ge |P'|_{\Pi,1} = |[q_1;q_2]|  \ge (1-2\delta)d(P)\ge \frac{1}{2}d(P).
 \end{align}
Moreover, we obtain for $k=1,2$  by the choice of $\gamma$ and $q_1,q_2$ 
\begin{align}\label{eq:n1*}
\dist_P(x,S_k) \le 3\delta d(P) \ \text{for all} \ x \in P_k.
 \end{align}
Indeed, e.g. for $k=1$, we observe $\dist_P(y,q_2) \ge (1-2\delta)d(P)$ for all $y \in S_1$ and $\dist_P(q_2,p_2) = \delta d(P)$. This implies $\dist_P(y,p_2) \ge (1-3\delta)d(P)$ for all $y \in S_1$, from which \eqref{eq:n1*} follows. Consequently, by \eqref{eq:n1*}, Lemma \ref{lemma: diam} and $\delta = \frac{1}{14}\vartheta \lambda$ we obtain
\begin{align}\label{eq:n2}
\min_{R \in SO(2)}|P'|_{\Pi,R} \ge \min_{R \in SO(2)} |P|_{\Pi,R} -6\delta d(P) \ge \lambda \max_{R \in SO(2)} |P|_{\Pi,R} -6\delta d(P) \ge \delta d(P).
\end{align}
(a1) If $P'$ is a convex polygon, we find by \eqref{eq:n2} and Lemma \ref{lemma: convpart2} that $P'$ contains a ball $B(x,r)$ with $r = \frac{\delta}{4} d(P)$ and thus, since $P \supset P'$, $P$ is $\frac{\delta}{4}$-rotund.  

\vspace{-0.2cm}
 (a2)  Otherwise, we choose a concave vertex $v \in \mathcal{V}_{P'}'$ which minimizes the distance to $[q_1; q_2]$. This implies that the triangle  with vertices $v,q_1,q_2$ is contained in $P'$ (see Figure \ref{rot1}). Understanding the vertices as complex numbers we define the phases $\varphi_1 = {\rm arg}(q_1-v) \in [0,2\pi)$ and $\varphi_2  = {\rm arg}(q_2-v) \in [0,2\pi)$, where possible after reflection of $P'$ along $\R \times \lbrace  0\rbrace$ and a rotation we can suppose that $0 \le \varphi_2 < \varphi_1 < 2\pi$ with $\varphi_1 - \varphi_2 < \pi$ (cf. Figure \ref{rot1}). 

We define the function $f:  [\varphi_2,\varphi_1] \to P$ so that $f(\varphi)$ denotes the closest point to $v$ on $(v + \R_+ e^{i\varphi}) \cap \partial P$. Observe that for $\varphi \in [\varphi_2,\varphi_1]$  each $[v;f(\varphi)]$  induces a partition $P = Q^\varphi_1 \cup Q^\varphi_2$ according to Definition \ref{def: induce} with $p_k,q_k \in Q_k^\varphi$ for $k=1,2$. Consequently,  by definition of $\gamma$ and $q_k$, $k=1,2$, we get $\min_{k=1,2} d(Q^\varphi_k) \ge \delta d(P)$ for all $\varphi \in [\varphi_1,\varphi_2]$ and then we obtain
\begin{align}\label{eq:n3}
|[v;f(\varphi)] | \ge \vartheta \min_{k=1,2} d(Q_k^\varphi) \ge \vartheta\delta d(P)
\end{align}
 since $P$  is $\vartheta$-semiconvex. Consequently, we derive that the circular sector
 $$T := \lbrace x \in \R^2: {\rm arg}(x-v) \in [\varphi_2,\varphi_1], \ |[x;v]| \le \vartheta\delta d(P) \rbrace $$
 is contained in $P$.  Clearly, we have $\dist(v,[q_1;q_2]) \le d(P)$, which in view of $|[q_1;q_2]| \ge \frac{1}{2} d(P)$ (see \eqref{eq:n1}) implies $\varphi_1 - \varphi_2 \ge \arctan(1/2)$ by elementary trigonometry. Then it is not hard to see that there is a ball  $B(x,r) \subset T \subset P$ with $r \ge c \vartheta \delta d(P)$ for a universal $c>0$ small enough. This yields that $P$ is $\omega$-rotund for some $\omega$ only depending on $\vartheta$, $\lambda$.   \smallskip

  (b) We now suppose that $ \gamma([\delta,1-\delta])$ is not a segment or $q_k \in \partial P$ for some $k=1,2$, i.e. we find $v \in \mathcal{V}_P'$ with $v \in \gamma([\delta,1 -\delta])$. Choose $v_-,v_+ \in \gamma$ such that $[v_-;v]$ and $[v;v_+]$ are contained in $\gamma$. Let $\varphi_- = {\rm arg}(v_- - v)$, $\varphi_+ = {\rm arg}(v_+-v)$ and without restriction, possibly after a rotation and reflection, we can  assume that $0 \le \varphi_+ < \varphi_- < 2\pi$ with $\varphi_- - \varphi_+ > \pi.$   We now proceed as in (a): 
  
  We see that $[v;f(\varphi)]$ induces a partition $P = Q^\varphi_1 \cup Q^\varphi_2$ for all $\varphi \in (\varphi_+,\varphi_-)$  with $p_k \in Q_k^\varphi$, $k=1,2$, and $|[v;f(\varphi)]| \ge \vartheta\min_{k=1,2} d(Q^\varphi_k) \ge  \vartheta\min_{k=1,2} d_P(v,p_k) \ge \vartheta \delta d(P)$ (cf. \eqref{eq:n3}). Then as before the set $ \lbrace x \in \R^2: {\rm arg}(x-v) \in (\varphi_+,\varphi_-), \ |[x;v]| \le \vartheta\delta d(P) \rbrace$ is contained in $P$. Since $\varphi_- - \varphi_+> \pi$, we conclude that $P$ contains a ball with radius larger than $c\vartheta\delta d(P)$. \eop

The result shows that if $\max_{R \in SO(2)} |P|_{\Pi,R}$ and $\min_{R \in SO(2)} |P|_{\Pi,R}$ are comparable, the polygon $P$ already has the desired properties. Otherwise, we will perform a partition of semiconvex polygons into semiconvex and rotund polygons as described in Section \ref{sec: part-semrot} below. To this end, it is crucial to characterize the position of concave vertices in a semiconvex polygon. The following result shows that for a semiconvex polygon, which is not already rotund, one can  identify (at most) two regions which contain the concave vertices.

\begin{theorem}\label{th: new}
 Let $0 < \vartheta <1$. Then there is a constant $C=C(\vartheta)>0$ such that the following holds for all $\vartheta$-semiconvex polygons $P$: There are two segments $S_1,S_2$ inducing a partition of $P = P_1 \cup P' \cup P_2$ with  $P_i \cap P' = S_i$ for $i=1,2$ such that $P'$ is a convex polygon and the polygons $P_i$ satisfy
\begin{align}\label{eq: bound sum}
\begin{split}
(i)& \ \ \mathcal{H}^1(S_i) \le \vartheta  \mathcal{H}^1(\partial P),\\
(ii)& \ \ \max_{R \in SO(2)}|P_i|_{\Pi,R} \le C\min_{R \in SO(2)}|P_i|_{\Pi,R},\\
(iii)& \ \ \max_{R \in SO(2)}|P_i|_{\Pi,R}  \le C\dist(v,S_i)  \ \text{ for all } \ v \in \mathcal{V}_{P_i}'.
\end{split}
\end{align}
\end{theorem}  

We remark that the choice $P_i = \emptyset$, $i=1,2$,   is admissible. (In this case also the corresponding segment is empty.) Moreover, also the choice $P_1 = P$, $P' = P_2 = \emptyset$ is possible, where Theorem \ref{th: only two} and \eqref{eq: bound sum}(ii) then imply that  $P$  is  rotund. Later, condition \eqref{eq: bound sum}(iii) will be crucial to show that $P_i$ are semiconvex using Lemma \ref{lemma: part-seg}. Theorem \ref{th: only two} together with \eqref{eq: bound sum}(ii) will then yield that the polygons $P_i$ are rotund.

\Proof Possibly after rotation we have $\min_{R \in SO(2)} |P|_{\Pi,R} = |P|_{\Pi,2}$.    Without restriction we can assume that $\vartheta^2|P|_{\Pi,1} > 12|P|_{\Pi,2}$ as otherwise the claim holds for $P_1 = P$, $P' = P_2 = \emptyset$ and $S_1 = S_2 = \emptyset$ with $C= 12/\vartheta^2 +1$, where  \eqref{eq: bound sum}(ii) for $P_1$ follows from $\max_{R \in SO(2)}|P|_{\Pi,R} \le  |P|_{\Pi,1} +  |P|_{\Pi,2}$  and \eqref{eq: bound sum}(i),(iii) are trivial. Moreover, possibly after another infinitesimal rotation we can suppose $\vartheta^2|P|_{\Pi,1} > 12|P|_{\Pi,2}$ and
\begin{align}\label{eq: rotation}
v_1 \e_1 \neq v_2 \e_1 \ \ \  \text{for all} \ \ v_1,v_2 \in \mathcal{V}_P, \ v_1 \neq v_2.
\end{align}
Choose $p_1,p_2 \in \partial P$ with $(p_2-p_1)\e_1 = |P|_{\Pi,1}$ and let $\gamma: [0,l(\gamma)] \to P$ be the piecewise affine curve between $p_1$, $p_2$, parametrized by arc length, with $d_P(p_1,p_2) = l(\gamma)$. Define $\mathcal{U}_1 = \lbrace v\in \mathcal{V}_P': d_P(v,p_1) \le d_P(v,p_2)\rbrace$ and $\mathcal{U}_2 = \mathcal{V}_P' \setminus \mathcal{U}_1$. We first cut off two small pieces near $p_1, p_2$ to obtain an auxiliary convex polygon. Afterwards, we define $P'$ and show \eqref{eq: bound sum}.\\

\noindent  \smallskip
\textit{Step 1: Definition of an auxiliary polygon}

\begin{figure}[H]
\centering

\begin{overpic}[width=0.76\linewidth,clip]{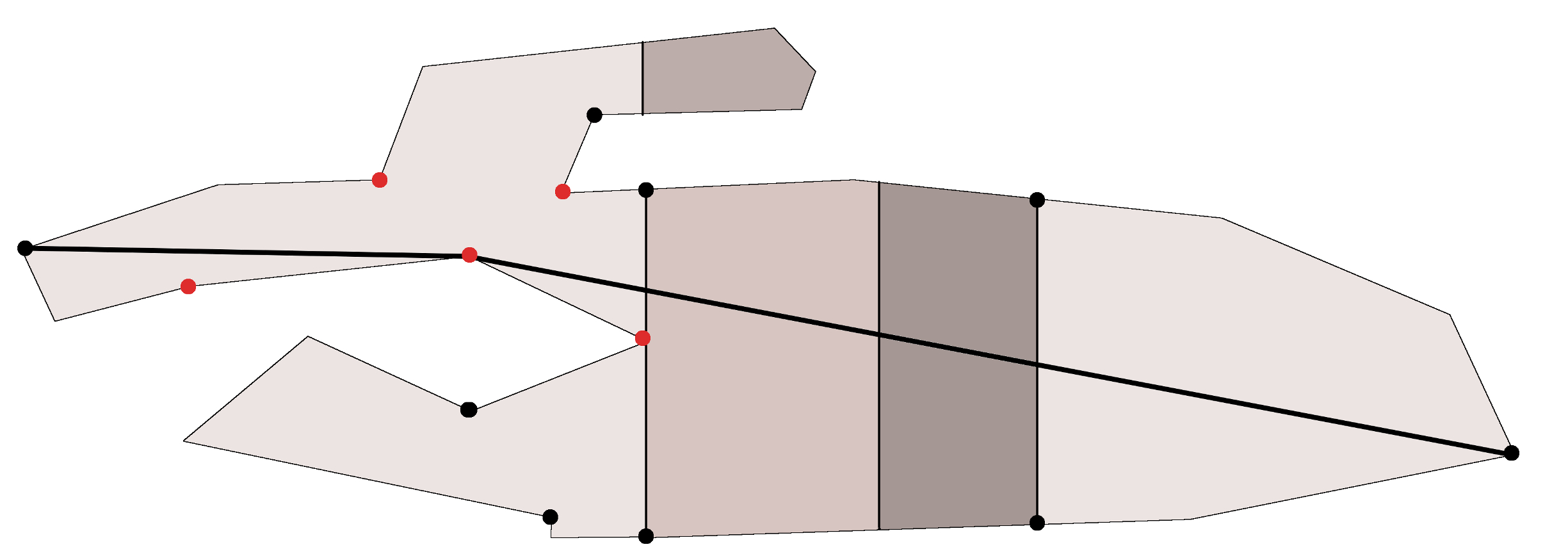}
\put(200,-2){{$q_1$}}
\put(200,74){{$q_2$}}
\put(132,90){{$R$}}
\put(212,47){{$S_1$}}
\put(229,18){{$\gamma$}}
\put(231,25) {\line(0,1){7}}
\put(186,20){{$A$}}
\put(149,20){{$B$}}
\put(54,22){$Q_1^{(1,2)}$}
\put(54,82){$Q_1^{(1,1)}$}
\put(68,82) {\line(1,0){13}}
\put(130,38){{$v_1$}}
\put(130,63){{$w^1_1$}}
\put(130,8){{$w^2_1$}}
\put(249,36){{$P'$}}
\put(203,49) {\line(1,0){7}}
\put(-6,65){{$p_1$}}
\put(288,11){{$p_2 = v_2$}}
\put(95,23){{$v'$}}
\put(107,11){{$v''$}}
\put(104,85){{$v'''$}}
\end{overpic}
\caption{We illustrate a case with  $I = \lbrace 1 \rbrace$.  In red the vertices $\mathcal{V}_*$ are depicted, where $v',v'', v''' \notin \mathcal{V}_*$. Moreover, we have $P_2^* = P' \cup A = P \cap T'$ and  $P_1^* = \hat{P}= P_2^* \cup B$. Note that $R \subset T \cap P$, but $R \cap P_1^* = \emptyset$.} \label{rot2}
\end{figure}
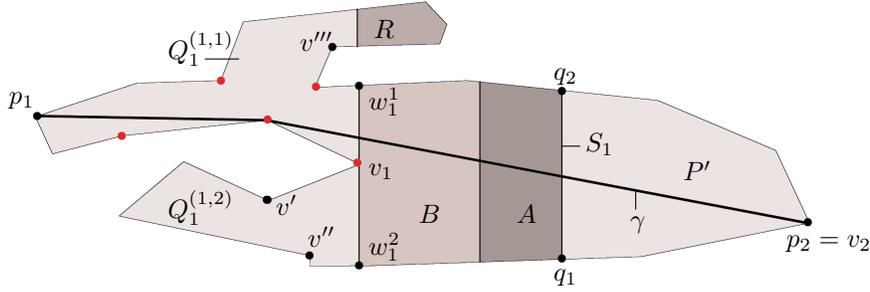

\noindent Let $\mathcal{V}_* \subset \mathcal{V}_P'$ be the vertices $v$ for which there is some $w \in \partial P$ such that $[v;w]$ is parallel to the $\e_2$-axis, $\gamma \cap [v;w] \neq \emptyset$ and  $[v;w]$ induces a partition of $P$ according to Definition \ref{def: induce} (see Figure \ref{rot2}). Note particularly that $v \in \mathcal{V}_*$ for each $v \in \mathcal{V}'_P$ with $v \in \gamma$.  Let $I = \lbrace i=1,2: \mathcal{U}_i \cap \mathcal{V}_* \neq \emptyset \rbrace$. For $i \in I$ we choose $v_i \in \mathcal{U}_i \cap \mathcal{V}_*$ with
 $$|(p_i - v_i)\e_1| = \max_{v \in \mathcal{U}_i \cap \mathcal{V}_*} |(p_i - v)\e_1|.  $$
 For $v_i$ we find a corresponding $w^1_i$ such that $[v_i;w^1_i]$ is parallel to the $\e_2$-axis, intersects $\gamma$ and induces a partition $P =  Q^{(i,1)}_1 \cup Q^{(i,1)}_2$, where the sets are labeled such that $p_k \in Q^{(i,1)}_k$ for $k =1,2$. Note that there may exist a second segment $[v_i;w^2_i]$ parallel to the the $\e_2$-axis inducing a partition  $P =  Q^{(i,2)}_1 \cup Q^{(i,2)}_2$ with $p_1,p_2 \notin Q_i^{(i,2)}$, cf. Figure \ref{rot2}. (If such a segment does not exist, we set $Q_i^{(i,2)} = \emptyset$ and $[v_i;w^2_i] = \emptyset$.) Note that $w^1_i,w^2_i \notin \mathcal{V}_P$ by \eqref{eq: rotation}.  Since $P$ is $\vartheta$-semiconvex, we derive using $\vartheta^2|P|_{\Pi,1} \ge 12|P|_{\Pi,2}$
$$\min_{k=1,2} d(Q^{(i,j)}_k) \le \vartheta^{-1}|[v_i;w^j_i]| \le  \vartheta^{-1}|P|_{\Pi,2} \le  \frac{1}{12}\vartheta|P|_{\Pi,1}.$$ 
As for $l =1,2$, $l \neq i$, we have $\dist(v_i,p_l) \ge \frac{1}{2} d_P(p_1,p_2) \ge \frac{1}{2}|P|_{\Pi,1}$ by definition of $\mathcal{U}_i$, we get $d(Q_l^{(i,j)}) \ge \frac{1}{2}|P|_{\Pi,1}$ and thus obtain 
\begin{align}\label{eq:endl1}
r_{i,j} := d(Q^{(i,j)}_i) \le \vartheta^{-1}|[v_i;w^j_i]| \le   \frac{1}{12}\vartheta|P|_{\Pi,1}.
\end{align}
If $i \notin I$, we set $r_{i,j}  = 0$ for $j=1,2$, $v_i = p_i$ and introduce the trivial partitions $P = Q^{(i,j)}_1 \cup Q^{(i,j)}_2$ with $Q_i^{(i,j)} = \lbrace p_i \rbrace$ and $Q_k^{(i,j)}  = P$ for $k \neq i$.  For shorthand we define $\bar{r}_i = \max\nolimits_{j=1,2}r_{i,j}$ for $i=1,2$. 

By the fact that $d_P(p_i,v_i) \le d(Q_i^{(i,1)}) \le \bar{r}_i$ and \eqref{eq:endl1} we have that the sets 
\begin{align}\label{eq:endl7}
T  := [v_1 \e_1 , v_2 \e_1] \times \R, \ \  \ \  T'  := [p_1 \e_1 + 2\bar{r}_1, p_2 \e_1 - 2\bar{r}_2] \times \R
\end{align}
satisfy $\emptyset \subsetneq   T' \subset T$. Consider the polygon $\hat{P} := Q^{(1,1)}_2 \cap Q^{(1,2)}_2\cap Q_1^{(2,1)} \cap Q_1^{(2,2)}$, which is confined, if existent, by the segments $[v^j_i;w^j_i]$, $i,j=1,2$. Moreover,  let $P^*_1$  be the connected component of  $P \cap T$ contained in $\hat{P}$. (See Figure \ref{rot2}. Below we will see that $P_1^* = \hat{P}$.)

As $P^*_1$ is connected, it is a polygon. We now show that $P^*_1$ is convex. Note that $v  \in {\rm int}(T)$ for all $v \in \mathcal{V}_{P^*_1}'$ since $P^*_1 \subset T$. Moreover, $v \notin {\rm int}(T)$ for all $\mathcal{V}_{P}' \cap \mathcal{V}_*$ by definition of $v_1,v_2$.  Consequently, $\gamma \cap P^*_1$ does not contain a concave vertex of $P$ and is thus a segment. Assume $\mathcal{V}_{P^*_1}' \neq \emptyset$ and choose a vertex $v \in \mathcal{V}_{P^*_1}' \subset \mathcal{V}_P'$ minimizing $\dist(v, \gamma \cap P^*_1)$.  Then there is $p \in \gamma$ such that $[v;p] \subset P$ and $[v;p]$ parallel to the $\e_2$-axis. This then implies  $v \in \mathcal{V}_*$, which gives a contradiction and shows that $P^*_1$ is convex. 

The convexity of $P^*_1$ together with the fact that $w^j_i \notin \mathcal{V}_P$ (see \eqref{eq: rotation}) also implies   $P^*_1 \cap \partial T = \bigcup_{i,j=1,2} [v_i^j;w_i^j]$ and this  yields $\hat{P} = P_1^*$. Moreover, we derive  
$$
P^*_2 := P \cap T' = P^*_1 \cap T'.
$$
Indeed, as  $P^*_1 = \hat{P}$, we obtain $P \setminus P_1^* \subset \bigcup_{j=1}^2 (Q^{(1,j)}_1 \cup Q_2^{(2,j)})$.  Then the definition of $r_{i,j}$ (see \eqref{eq:endl1}) together with  $r_{i,j} \le \bar{r}_i$,   \eqref{eq:endl7}   and $Q^{(i,j)}_i \cap\partial T \neq \emptyset$, $i,j=1,2$, yields $Q^{(i,j)}_i \cap T' = \emptyset$ for $i,j=1,2$. This implies the claim.  Since $P^*_1$ is convex, also $P^*_2$ is convex. \smallskip

\noindent 
\textit{Step 2: Definition of $P'$}\\
\noindent  We are now in a position to define $P'$. As $3\bar{r}_i \le \frac{1}{4}|P|_{\Pi,1}$ by   \eqref{eq:endl1},  we can choose $t_1<t_2$ with 
\begin{align}\label{eq:eeextra2}
p_1 \e_1 + 3\bar{r}_1 \le t_1 \le  p_1 \e_1 + \frac{1}{4}|P|_{\Pi,1}, \ \ \ p_2 \e_1 - \frac{1}{4}|P|_{\Pi,1} \le t_2 \le  p_2 \e_1 - 3\bar{r}_2 
\end{align} 
such that $S_i := P \cap (\lbrace t_i \rbrace \times \R)$ satisfy
\begin{align}\label{eq:eeextra}
\mathcal{H}^1(S_i) \le |t_i - p_i\e_1| \le \max\lbrace \mathcal{H}^1(S_i), 3\bar{r}_i \rbrace.
\end{align}
This follows from a continuity argument taking $|P|_{\Pi,1} \ge 12|P|_{\Pi,2}$ into account. Clearly, $P' := P \cap ([t_1,t_2] \times \R)$ is a convex polygon (cf. again Figure \ref{rot2}). Denote the closures of the at most two connected components of $P \setminus P'$ by $P_1, P_2$, where $P_i =\emptyset$ if and only if $i \notin I$ and note that indeed $S_i = P' \cap P_i$ for $i \in I$. It remains to confirm \eqref{eq: bound sum}. As a preparation, we show that  there is a universal $C>0$   such that for $i \in I$
\begin{align}\label{eq:endl6}
d(P_i) \le C\vartheta^{-1} (|t_i - p_i \e_1| + \mathcal{H}^1(S_i)).
\end{align}
We confirm the claim e.g. for $i=1$. Let $q_1$, $q_2$ be the endpoints of the segment $S_1$. As $P^*_1$ is convex, we have that the closed triangle $\triangle$ with vertices $q_1,q_2$ and $v_1$ is contained in $P^*_1 \subset P$. If $[q_j;v_1]$ is not completely contained in $\partial P^*_1$ for $j=1,2$, it induces a partition $P = R^{(j)}_1 \cup R^{(j)}_2$, where the sets are labeled such that $R^{(j)}_1 \subset P_1$. If $[q_j;v_1] \subset \partial P^*_1$, we set  $R^{(j)}_1 = \emptyset$. We obtain $P_1 = R^{(1)}_1 \cup R^{(2)}_1 \cup \triangle$. Note that $d(\triangle) \le |t_1 - p_1 \e_1| + \mathcal{H}^1(S_1)$ and due to the fact that $P$ is $\vartheta$-semiconvex,  we get 
$$\min_{k=1,2} d(R^{(j)}_k)  \le \vartheta^{-1}|[q_j;v_1]| \le \vartheta^{-1} (|t_1 - p_1 \e_1|  + \mathcal{H}^1(S_1)) \le \frac{1}{3}|P|_{\Pi,1},$$
 where in the last step we used \eqref{eq:endl1}, \eqref{eq:eeextra} and that by assumption $\mathcal{H}^1(S_1)\le|P|_{\Pi,2} \le \frac{1}{12}\vartheta^2 |P|_{\Pi,1}$. Since $d(R^{(j)}_2) \ge \frac{3}{4}|P|_{\Pi,1}$ by \eqref{eq:eeextra2}, we derive $d(R^{(j)}_1)  \le d(R^{(j)}_2)$ and then \eqref{eq:endl6} follows.

\vspace{-0.4cm}

We now show \eqref{eq: bound sum}.  First, (i) follows from $\mathcal{H}^1(S_i) \le |P|_{\Pi,2}$ and $12|P|_{\Pi,2} \le \vartheta^2|P|_{\Pi,1} \le \vartheta^2 \mathcal{H}^1(\partial P)$. If $R \e_1$,  $R \in SO(2)$,  encloses an angle  smaller than $\frac{\pi}{4}$ with the $\e_2$-axis, we find $ |P_i|_{\Pi,R} \ge \frac{1}{\sqrt{2}} \max\lbrace |[v_i;w^1_i]| + |[v_i;w^2_i]|, \mathcal{H}^1(S_i) \rbrace$. Likewise, if $R \e_1$ encloses an angle  smaller than $\frac{\pi}{4}$ with the $\e_1$-axis, we get in view of \eqref{eq:eeextra2} 
   $$ |P_i|_{\Pi,R} \ge |\gamma \cap (P_i \cap T)|_{\Pi,R} \ge c|t_i - v_i \e_1|\ge c(|t_i - p_i \e_1| - \bar{r}_i) \ge c|t_i - p_i \e_1|$$
    for a universal $c$ small enough,    where we used that $\gamma \cap (P_i \cap T)$ is a segment enclosing a small angle with $\e_1$ since $\vartheta^2|P|_{\Pi,1} \ge 12|P|_{\Pi,2}$ (cf. Figure \ref{rot2}). By \eqref{eq:endl1} and \eqref{eq:eeextra} this implies 
$$\min_{R \in SO(2)} |P_i|_{\Pi,R} \ge c\min \lbrace \max\lbrace \vartheta \bar{r}_i, \mathcal{H}^1(S_i) \rbrace, |t_i - p_i \e_1|\rbrace \ge c\vartheta |t_i - p_i \e_1|.$$
Consequently, \eqref{eq:eeextra} and \eqref{eq:endl6}  yield
\begin{align*}
\max_{R \in SO(2)}|P_i|_{\Pi,R} &\le d(P_i) \le C\vartheta^{-1} 2|t_i - p_i \e_1|  \le  C \vartheta^{-2} \min_{R \in SO(2)} |P_i|_{P,R}.
\end{align*}
This gives (ii). Finally, to see (iii), we recall that  $P^*_2 = P \cap T'$ is a convex polygon and thus in view of \eqref{eq:endl7}, \eqref{eq:eeextra2}, we get $\dist(v,S_i) \ge |t_i - p_i \e_1| -2 \bar{r}_i \ge   \frac{1}{3}|t_i - p_i \e_1| =  \frac{1}{3}|P_i|_{\Pi,1}$ for all $v \in \mathcal{V}_{P_i}'$. The claim  now follows from \eqref{eq: bound sum}(ii).  \eop

\subsection{Partitions into semiconvex and rotund polygons}\label{sec: part-semrot}

We now show that semiconvex polygons can be partitioned into semiconvex and rotund polygons. We start with the partition of   convex polygons   into rotund polygons by introducing segments parallel to the direction of shortest extend.

\begin{lemma}\label{lemma: convpart}
Let $\theta>0$. Then there is $\omega = \omega(\theta)>0$ such that for all convex polygons $P$, satisfying $\sphericalangle(v,P) \ge \frac{\pi}{4}$ for all vertices $v \in \mathcal{V}_P$, there is a partition $P = P_1 \cup \ldots \cup P_N$ with
\begin{align}\label{eq: length bound convexi}
\mathcal{H}^1\Big( \bigcup\nolimits^N_{j=1}\partial P_j \setminus \partial P\Big) \le \theta\mathcal{H}^1(\partial 
P)
\end{align}
and the polygons $(P_j)_{j=1}^N$ are $\omega$-rotund.
\end{lemma}

\Proof After rotation we may assume that $\min_{R \in SO(2)} |P|_{\Pi,R} = |P|_{\Pi,2}$.  Clearly, it is not restrictive to suppose that $\theta \le \theta_0$ for some $\theta_0 \le 1$ to be specified below.  If $|P|_{\Pi,1} < 7\theta^{-1} |P|_{\Pi,2}$, we obtain $\max_{R \in SO(2)} |P|_{\Pi,R} \le |P|_{\Pi,1}+|P|_{\Pi,2} \le (1+7\theta^{-1} )|P|_{\Pi,2} $ and $P$ is $\omega$-rotund by Theorem \ref{th: only two} for $\omega$ only depending on $\theta$.

\begin{figure}[H]
\centering
\begin{overpic}[width=0.49\linewidth,clip]{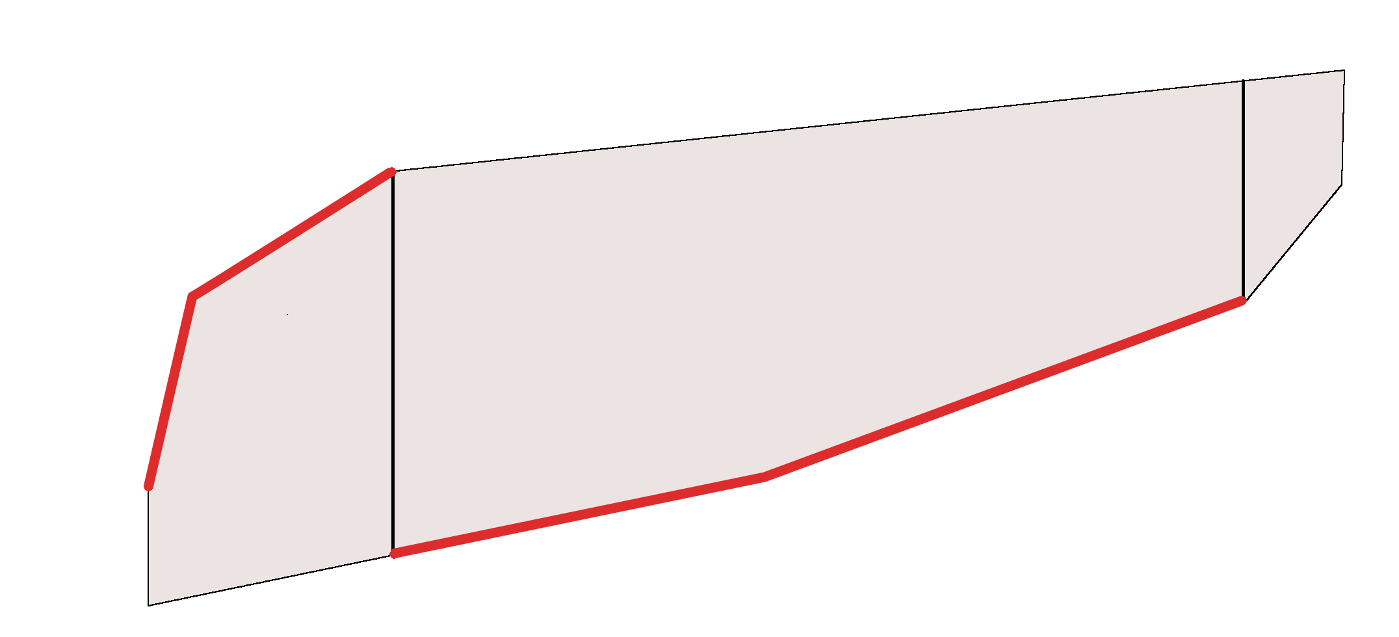}
\put(64.5,37){{$S_{s_1}$}}
\put(58,60){{$u$}}
\put(58,15){{$v$}}
\put(55.5,40) {\line(1,0){5}}
\put(160,56){{$S_{s_2}$}}
\put(172,59) {\line(1,0){5}}
 \put(92,28){{$\gamma$}}
\put(96,20) {\line(0,1){5}}
\put(32,37){{$\gamma_1$}}
\put(35,17){{$Q_1^{s_1}$}}
\put(26,40) {\line(1,0){5}}
\put(197,66){{$Q_2^{s_2}$}}
\put(189,70) {\line(1,0){5}}

\end{overpic}
\caption{We depicted the curves $\gamma_1$ and $\gamma$ considered above in the proof of \eqref{eq: angle1}. Note that similar arguments involving the angle between tangent vectors and $\e_1$ are also used in \eqref{eq: arguments} and \eqref{eq: angle10}.}

  \label{conv1}
\end{figure}
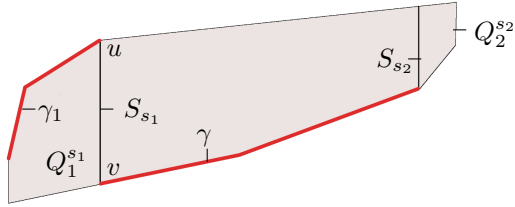

Now assume $|P|_{\Pi,1} \ge 7\theta^{-1} |P|_{\Pi,2}$. For $t \in \R$ we denote by $S_t$ the segments  $S_t = (\lbrace t\rbrace  \times \R ) \cap P$ which induce  partitions $P = Q^t_1 \cup Q^t_2$, where $Q^t_2 \subset \lbrace x_1 \ge t \rbrace$. For shorthand we write $\varphi_\theta = \arctan\theta$.  Choose the smallest $s_1$ and the largest $s_2$ such that the polygon $P' := \overline{P \setminus (Q^{s_1}_1 \cup Q^{s_2}_2)}$ with $\mathcal{V}_{P'} = (u_1,\ldots,u_{n})$ satisfies 
\begin{align}\label{eq: angle3}
 [u_i; u_{i+1}] \neq S_{s_1}, S_{s_2} \   \Rightarrow  \ {\rm arg}(u_{i+1} - u_i) \in    \big(\lbrace 0, \pi \rbrace  + [- \varphi_\theta, \varphi_\theta]\big) {\rm mod}2\pi,
 \end{align}
 i.e. $u_{i+1} - u_i$ and $\e_1$ enclose an angle smaller than $\varphi_\theta$. We  show that for $j=1,2$ one has 
\begin{align}\label{eq: angle1}
\begin{split}
(i)& \ \ |Q^{s_1}_1|_{\Pi,1}  + |Q^{s_2}_2|_{\Pi,1} \le \frac{1}{2}|P|_{\Pi,1}, \\
(ii)& \ \ |Q^{s_j}_j|_{\Pi,1} \le 3\theta^{-1}\mathcal{H}^1(S_{s_j}), \\
(iii)& \ \ 0<|Q^{s_j}_j|_{\Pi,2} \le 4\mathcal{H}^1(S_{s_j}).
\end{split}
\end{align} 
By convexity of $P$ and the choice in \eqref{eq: angle3} we find curves $\gamma_j$ in $\partial Q_j^{s_j}$ with $|\gamma_j|_{\Pi,1} = |Q_j^{s_j}|_{\Pi,1}$  such that the angle enclosed by $\e_1$ and the tangent vector $\gamma'_j$ of $\gamma_j$ is larger than $\varphi_\theta$ (see Figure \ref{conv1}).  By $|P|_{\Pi,1} \ge 7\theta^{-1} |P|_{\Pi,2}$ a short calculation then yields 
$$
\sum\nolimits_j \theta|Q_j^{s_j}|_{\Pi,1}  = \sum\nolimits_j \tan(\varphi_\theta)|Q_j^{s_j} \cap \gamma_j|_{\Pi,1} \le \sum\nolimits_j  |Q_j^{s_j}|_{\Pi,2} \le  2|P|_{\Pi,2} \le \frac{2\theta}{7}  |P|_{\Pi,1},
$$
which gives (i). The first inequality in (iii) follows from the fact that $Q^{s_j}_j$ cannot be degenerated to a single vertex as for $\theta_0$ small in view of \eqref{eq: angle3} this would contradict the lower bound $\frac{\pi}{4}$ on the interior angles of $P$. Note, however, that $|Q^{s_j}_j|_{\Pi,1} =0$ is possible, where in this case we have $Q^{s_j}_j = S_{s_j}$. 

We now show (ii)-(iii), e.g. for $S_{s_1}$. Recalling the previous observation we see that the claim is clear if $|Q^{s_1}_1|_{\Pi,1} =0$ and we therefore assume $|Q^{s_1}_1 |_{\Pi,1} >0$. For convenience define $[u;v] := S_{s_1}$.  We now get 
$$(a) \   \min\lbrace\sphericalangle(u,Q_1^{s_1}),\sphericalangle(v,Q_1^{s_1})\rbrace \le \frac{\pi}{2} - \varphi_\theta, \ \ \ \ \ (b) \  \max\lbrace\sphericalangle(u,Q_1^{s_1}),\sphericalangle(v,Q_1^{s_1})\rbrace \le \frac{\pi}{2} + \frac{\varphi_\theta}{2}.$$
If (b) was false, as above using the convexity of $P$ we would find a curve $\gamma$ in $\partial P'$ with $|\gamma|_{\Pi,1} = |P'|_{\Pi,1}$  such that the angle enclosed by $\e_1$ and $\gamma'$ is larger than  $\frac{\varphi_\theta}{2}$ (see Figure \ref{conv1}). But then  similarly as in the proof of \eqref{eq: angle1}(i) we would find, 
by a Taylor expansion for $\theta_0$ small, and $|P'|_{\Pi,1} \ge  \frac{1}{2}|P|_{\Pi,1}$ (see \eqref{eq: angle1}(i)) 
$$|P|_{\Pi,2} \ge |P'|_{\Pi,2} \ge \tan\Big(\frac{\varphi_\theta}{2}\Big)|P'\cap\gamma|_{\Pi,1} \ge  \frac{2\theta}{5}|P'|_{\Pi,1} \ge  \frac{\theta}{5}|P|_{\Pi,1},$$ which contradicts the assumption. To see (a), we observe that the construction of $P'$ implies that, up to changing the roles of $u$ and $v$, $|\sphericalangle(u,Q_1^{s_1})-\frac{\pi}{2}| \ge \varphi_\theta$. As in the proof of (b), we derive that it is not possible that the angle is obtuse. This gives (a). 

Combining (a) and (b) and recalling the convexity of $P$ we derive
$$\mathcal{H}^1(S_{s_1}) + \tan(\varphi_\theta/2)|Q^{s_1}_1|_{\Pi,1} - \tan(\varphi_\theta)|Q^{s_1}_1|_{\Pi,1}\ge 0$$
and then for $\theta_0$ small  by a Taylor expansion $\mathcal{H}^1(S_{s_1}) \ge \frac{\theta}{3}|Q^{s_1}_1|_{\Pi,1}$, i.e. (ii) holds. 

\vspace{-0.3cm}
 
\vspace{-0.1cm}
To see the second inequality in (iii), we again use  (a),(b) and (ii) to obtain for $\theta$ small 
\begin{align}\label{eq: arguments}
\hspace{-0.2cm}|Q^{s_1}_1|_{\Pi,2} \le \mathcal{H}^1(S_{s_1}) + \tan(\varphi_\theta/2)|Q_1^{s_1}|_{\Pi,1} \le \mathcal{H}^1(S_{s_1}) + \theta|Q_1^{s_1}|_{\Pi,1}\le 4\mathcal{H}^1(S_{s_1}).
\end{align} 
For later purpose note that \eqref{eq: angle1}(ii) and the assumption $|P|_{\Pi,1} \ge 7\theta^{-1} |P|_{\Pi,2}$  also imply 
 \begin{align}\label{eq: P}
 |P'|_{\Pi,1} \ge  |P|_{\Pi,1} -3\theta^{-1}(\mathcal{H}^1(S_{s_1}) + \mathcal{H}^1(S_{s_2}))\ge |P|_{\Pi,1} - 6\theta^{-1}|P|_{\Pi,2} \ge \theta^{-1}|P|_{\Pi,2}.
 \end{align}
We are now in a position to partition $P'$ with vertical segments: we assume segments $S_{t_1},S_{t_2}, \ldots, S_{t_n}$ with $s_1= t_1 < t_2 < \ldots < t_n$ and $P'_j = \overline{Q^{t_{j}}_2 \setminus Q^{t_{j+1}}_2}$ for $j=1,\ldots,n-1$, $n \in \N$, have been constructed with 
\begin{align}\label{eq: second bound2}
\hspace{-0.2cm}0<|P'_j|_{\Pi,1}  = \theta^{-1}\mathcal{H}^1(S_{t_{j+1}}), \, j=1,\ldots,n-1, \ \ \ \ \  |Q^{t_{n}}_2 \cap P'|_{\Pi,1}  \ge \theta^{-1}\mathcal{H}^1(S_{s_2}).
\end{align}
We observe that the latter condition in \eqref{eq: second bound2} holds in the case $n=1$, where no set has been constructed yet. In fact, we have  $Q^{t_1}_2 \cap P' = P'$ and then $|Q^{t_1}_2 \cap P'|_{\Pi,1} \ge \theta^{-1}|P|_{\Pi,2} \ge  \theta^{-1}\mathcal{H}^1(S_{s_2})$ by \eqref{eq: P}. Moreover, recall $\mathcal{H}^1(S_{t_{1}})>0 $ by \eqref{eq: angle1}(iii).

If  $|Q^{t_{n}}_2 \cap P'|_{\Pi,1}  \le 4\theta^{-1}\mathcal{H}^1(S_{s_2})$, we set $P'_{n} = Q^{t_{n}}_2 \cap P'$, $t_{n+1} = s_2$, $S_{t_{n+1}} = S_{s_2}$ and stop. For later reference we note that in this case
\begin{align}\label{eq: last thing}
\theta^{-1}\mathcal{H}^1(S_{s_{2}}) \le |P'_{n}|_{\Pi,1}  \le 4\theta^{-1}\mathcal{H}^1(S_{s_{2}}).
\end{align}
Otherwise, we have $|Q^{t_{n}}_2 \cap P'|_{\Pi,1}  > 4\theta^{-1}\mathcal{H}^1(S_{s_2})$. As $\mathcal{H}^1(S_t)$ is continuous in $t$ and $\mathcal{H}^1(S_{t_n})>0$ by \eqref{eq: second bound2}, we apply the intermediate value theorem and find some $t_{n+1} \in [t_{n}, s_2]$ such that $P'_{n} := \overline{Q^{t_{n}}_2 \setminus Q^{t_{n+1}}_2}$ satisfies $|P'_{n}|_{\Pi,1} =  \theta^{-1}\mathcal{H}^1(S_{t_{n+1}})$. This gives
$$|Q^{t_{n+1}}_2 \cap P'|_{\Pi,1}  = |Q^{t_{n}}_2 \cap P'|_{\Pi,1} - |P'_{n}|_{\Pi,1}\ge \theta^{-1}(4\mathcal{H}^1(S_{s_2}) -  \mathcal{H}^1(S_{t_{n+1}})).$$
Then using $\mathcal{H}^1(S_{t_{n+1}}) \le \mathcal{H}^1(S_{s_2}) + 2\theta|Q^{t_{n+1}}_2 \cap P'|_{\Pi,1}$ by  \eqref{eq: angle3} (see \eqref{eq: arguments} for a similar argument), we get
\begin{align*}
3|Q^{t_{n+1}}_2 \cap P'|_{\Pi,1}  \ge \theta^{-1}(4\mathcal{H}^1(S_{s_2}) -  \mathcal{H}^1(S_{t_{n+1}})) + 2|Q^{t_{n +1}}_2 \cap P'|_{\Pi,1} \ge  3\theta^{-1}\mathcal{H}^1(S_{s_2}),
\end{align*}
which gives the second part of \eqref{eq: second bound2}. We now proceed with the next iteration step and observe that the construction stops after a finite number of steps with a partition $P' = P'_1 \cup \ldots P'_N$ since by convexity of $P$ we have $\mathcal{H}^1(S_{t}) \ge \min_{i=1,2} \mathcal{H}^1(S_{s_i})>0$ for all $t \in [s_1,s_2]$.

Note that by \eqref{eq: second bound2}, \eqref{eq: last thing} each $P'_n$ contains a triangle with a base of length $\mathcal{H}^1(S_{t_{n+1}})$ and a height with length in the intervall $ [{\theta}^{-1}\mathcal{H}^1(S_{t_{n+1}}),4\theta^{-1}\mathcal{H}^1(S_{t_{n+1}})]$. By \eqref{eq: angle3} it is not hard to see that each of these triangles contains a ball with radius larger than $C\mathcal{H}^1(S_{t_{n+1}})$ for a universal $C>0$ small enough. Likewise, again arguing as in \eqref{eq: arguments}, by \eqref{eq: angle3}, \eqref{eq: second bound2}, \eqref{eq: last thing} we also find 
\begin{align}\label{eq: angle10}
|P'_n|_{\Pi,2} \le \mathcal{H}^1(S_{t_{n+1}}) + 2\theta |P'_n|_{\Pi,1} \le 9 \mathcal{H}^1(S_{t_{n+1}}). 
\end{align}
Consequently, again by \eqref{eq: second bound2}, \eqref{eq: last thing} we derive $d(P_n') \le (4\theta^{-1}+9)\mathcal{H}^1(S_{t_{n+1}})$ and thus we conclude that $(P'_n)_{n=1}^N$ are $\omega$-rotund for  some $\omega = \omega(\theta)$ small enough.

Define $P_1 = P'_1 \cup Q^{s_1}_1$, $P_N = P'_N \cup Q^{s_2}_2$ and $P_n = P_n'$ for $n=2,\ldots,N-1$. Clearly, $P_2, \ldots, P_{N-1}$ are $\omega$-rotund. Applying \eqref{eq: angle1}(ii),(iii) and \eqref{eq: second bound2} we get 
$$d(P_1) \le d(P_1') + C\theta^{-1}\mathcal{H}^1(S_{s_1}) \le d(P_1') + C\theta^{-1}(\mathcal{H}^1(S_{t_2}) + 2\theta|P'_1|_{\Pi,1} )\le Cd(P_1'),$$
where in the penultimate step we once again exploited \eqref{eq: angle3}. A similar expression holds for $P_N$. Consequently, possibly passing to a smaller $\omega$ also $P_1,P_N$ are $\omega$-rotund. Finally, to see \eqref{eq: length bound convexi}, we compute by \eqref{eq: second bound2} 
$$\mathcal{H}^1\Big( \bigcup\nolimits^N_{j=1}\partial P_j \setminus \partial P\Big) = \sum^{N}_{j=2} \mathcal{H}^1(S_{t_{j}}) \le \theta \sum^N_{j=1} |P'_{j}|_{\Pi,1} \le \theta |P|_{\Pi,1}  \le \theta \mathcal{H}^1(\partial P).$$
 \eop

We now finally show that semiconvex polygons  can be partitioned into semiconvex and rotund polygons up to an arbitrary small set.

\begin{theorem}\label{th: qussc}
Let $\theta, \vartheta, \epsilon >0$. Then there are $\omega=\omega(\vartheta,\theta)$, $\bar{\vartheta} = \bar{\vartheta}(\vartheta,\theta)$  and a universal constant $C>0$ such that the following holds: For all  $\vartheta$-semiconvex polygons $P$ there is a  partition  $P = P_0 \cup \ldots \cup P_N$ with $\mathcal{H}^1(\partial P_0) \le \epsilon$ and 
\begin{align}\label{eq: addition}
\sum^N_{i=1}\mathcal{H}^1(\partial P_i) \le (1+ C\theta)\mathcal{H}^1(\partial P)
\end{align}
such that the polygons $(P_i)_{i=1}^N$ are $\bar{\vartheta}$-semiconvex and  $\omega$-rotund.
\end{theorem}

\Proof Possibly by passing to a smaller $\vartheta$, we can assume that $\vartheta \le \theta$ in the following since \eqref{eq: semi}  still holds for a smaller value of $\vartheta$. We apply Theorem \ref{th: new} to obtain a partition $P = P_1 \cup P' \cup P_2$ such that $P'$ is a convex polygon and $P_i$ satisfy \eqref{eq: bound sum} with $S_i := P_i \cap P'$ for $i=1,2$. (Recall that some of the polygons may be empty.)

 We now first concern ourselves with $P'$. By $\mathcal{V}_\sphericalangle$ we denote the vertices $v \in \mathcal{V}_{P'}$  with $\sphericalangle(v,P') < \frac{\pi}{4}$. For each $v \in \mathcal{V}_\sphericalangle$ we choose a (closed) isosceles triangle $\triangle_v \subset P'$ with $v \in \triangle_v$ such that $\sphericalangle(v,\triangle_v) = \sphericalangle(v,P')$ is the only angle smaller than $\frac{\pi}{4}$ and we obtain $\mathcal{H}^1(\partial P_0) \le \epsilon$ as well as $\sphericalangle(u,P'') \ge \frac{\pi}{4}$ for all $u \in \mathcal{V}_{P
''}$, where $P_0 = \bigcup_{v \in \mathcal{V}_\sphericalangle} \triangle_v$ and  $P'' = \overline{P' \setminus P_0}$. We notice that by the triangle inequality
\begin{align}\label{eq: addition*}
\mathcal{H}^1(\partial P'') \le \mathcal{H}^1(\partial P') = \mathcal{H}^1(\partial P) + \sum^2_{i=1} (\mathcal{H}^1(S_i) - \mathcal{H}^1(\partial P_i\setminus S_i)) \le  \mathcal{H}^1(\partial P).
\end{align}
We apply Lemma \ref{lemma: convpart} on $P''$  to obtain a partition  $P'' = P_3 \cup \ldots \cup P_N$ with 
$$
\mathcal{H}^1 \Big(\bigcup\nolimits^N_{j=3}\partial P_j \setminus \partial P''  \Big) \le \theta\mathcal{H}^1(\partial P'')\le \theta\mathcal{H}^1(\partial P)
$$
such that the polygons $(P_j)_{j=3}^N$ are convex and $\omega$-rotund for some $\omega$ only depending on $\theta$. Since each $x \in \bigcup\nolimits^N_{j=3}\partial P_j \setminus \partial P''$ is contained in exactly two components, we compute by \eqref{eq: bound sum}(i), \eqref{eq: addition*} and $\vartheta \le \theta$
\begin{align*}
\sum\nolimits^N_{j=1}\mathcal{H}^1(\partial P_j) &  \le  \sum\nolimits^2_{i=1}\mathcal{H}^1(\partial P_i)  + \mathcal{H}^1(\partial P'') + 2\theta\mathcal{H}^1(\partial P)\\&
\le  \mathcal{H}^1(\partial P) + \sum\nolimits^2_{i=1} (\mathcal{H}^1(\partial P_i) + \mathcal{H}^1(S_i) - \mathcal{H}^1(\partial P_i \setminus S_i)) + 2\theta\mathcal{H}^1(\partial P) \\
& = \mathcal{H}^1(\partial P) + 2\theta\mathcal{H}^1(\partial P) + 2\sum\nolimits^2_{i=1} \mathcal{H}^1(S_i) \le  \mathcal{H}^1(\partial P) + 6\theta\mathcal{H}^1(\partial P).
\end{align*}
This gives \eqref{eq: addition}. It remains to show that $P_1$ and $P_2$, if existent, are semiconvex and rotund. We denote the endpoints of the segments $S_i$ by $u^i_1,u^i_2$, $i=1,2$. First, by \eqref{eq: bound sum}(iii) each $v \in \mathcal{V}_{P_i}'$ satisfies $\mathcal{H}^1(S_i) \le C\dist(v,S_i)$ for $C=C(\vartheta)$ and therefore an elementary geometric argument implies that there is an angle $\alpha=\alpha(\vartheta)>0$ such that $\max_{k=1,2}\sphericalangle(\triangle_v,u^i_k ) \ge \alpha$ for all $v \in \mathcal{V}_{P_i}'$, where $\triangle_v$ denotes the triangle formed by $u_1^i$, $u^i_2$ and $v$. Thus, recalling that $P$ is $\vartheta$-semiconvex, we get that $P_i$ are $\bar{\vartheta}$-semiconvex by Lemma \ref{lemma: part-seg} for $\bar{\vartheta}$ only depending on $\vartheta$. 

Finally, the fact that $P_i$ is $\bar{\vartheta}$-semiconvex together with \eqref{eq: bound sum}(ii) yields that $P_i$ is $\omega$-rotund by Theorem \ref{th: only two} for $\omega$ only depending on $\vartheta$. \eop

\begin{rem}\label{rem: vertices2}

{\normalfont

As in Remark \ref{rem: vertices1} we note that by the partition no additional concave vertices are introduced.

 }
\end{rem} 

\section{Equivalence of John domains and semiconvex, rotund polygons}\label{sec: equi}

In this section we study the relation of semiconvex, rotund polygons and John domains. This together with the partitions introduced in the last sections will allow us to give the proof of Theorem \ref{th: main part}. In the following for convenience we will say that a polygon $P$ is a $\varrho$-John domain if ${\rm int}(P)$ is a $\varrho$-John domain. We first observe that polygons, which are $\varrho$-John domains, are semiconvex and rotund.

\begin{lemma}\label{eq: one direc}
Let $0 < \varrho \le 1$. Each polygon $P$ which is a  $\varrho$-John domain is $\vartheta$-semiconvex and $\omega$-rotund  for $\vartheta,\omega$   only depending on $\varrho$.
\end{lemma}

\Proof Since there is $x \in P$ with $P \setminus B(x,d(P)/2) \neq \emptyset$, Lemma \ref{lemma: plump} implies that $P$ is $\frac{1}{4}\varrho$-rotund. If $P$ was not $\vartheta$-semiconvex for $\vartheta = \frac{\varrho}{4}$, there would be $u_1,u_2 \in \partial P$, $u_1 \in \mathcal{V}_P'$, inducing a partition $P = Q_1 \cup Q_2$ such that 
$$|[u_1;u_2]| < \vartheta \min_{k=1,2} d(Q_k) \le \frac{1}{4}\min_{k=1,2} d(Q_k).$$
We can choose  $v_k \in Q_k$ such that 
$d_P(v_k,w) \ge \frac{1}{4} d(Q_k) \ge 
\frac{1}{4}\vartheta^{-1}|[u_1;u_2]| $ for all $w \in [u_1;u_2]$. Let $\gamma$ be a John curve  between $v_1$, $v_2$ (see Remark \ref{rem: john}) and let  $w_*$ be an intersection point of $\gamma$ with $[u_1;u_2]$. As   ${\rm cig}(\gamma, \varrho) \subset P$, we derive  $B(w_*,\frac{\varrho}{4\vartheta}|[u_1;u_2]| ) \subset P$. In view of $\vartheta = \frac{\varrho}{4}$, this gives a  a contradiction.\eop

We now show that semiconvex and rotund polygons are John domains with controllable  John constant. Recall the notation  $x \e_j$ for the $j$-th component of points $x \in \R^2$ and that sometimes points are understood as complex numbers (see Section \ref{sec: prep}).

\begin{theorem}\label{th: John}
Let $\vartheta,\omega >0$. Then there is $\varrho=\varrho(\vartheta,\omega)$ such that each $\vartheta$-semiconvex and $\omega$-rotund polygon $P$ is a $\varrho$-John domain. 
\end{theorem}

\Proof By $0<c<1,C\ge 1$ we denote generic constants which are always independent of $\vartheta, \omega$. Possibly by passing to smaller $\vartheta, \omega$ we can assume that $\vartheta, \omega$ are sufficiently small with respect to $C$ and $\vartheta$ is small with respect to $\omega$  in the following proof since the properties in Definition \ref{def: semi}(i) and Definition \ref{def: qussc} still hold for smaller values of $\vartheta, \omega$.  As $P$ is $\omega$-rotund, we find some $p \in P$ and $r \ge \omega d(P)$ such that $B(p,r) \subset P$. Let $x \in {\rm int}(P)$ arbitrary. The goal is to construct a curve $\gamma$ between $x$ and $p$ such that for $\vartheta$ small enough
\begin{align}\label{eq: carot}
{\rm car}(\gamma, \vartheta^3) \subset P,
\end{align}
where ${\rm car}(\gamma,\vartheta^3)$ as in \eqref{eq: car}. This then shows that ${\rm int}(P)$ is a $\vartheta^3$-John domain. The construction will involve several steps. \smallskip\\

\noindent  \smallskip
\textit{Step 1: Preparations}

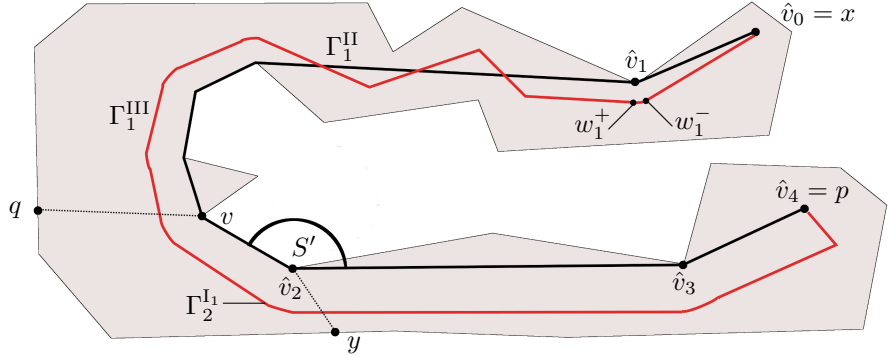
\begin{figure}[H]
\centering
\begin{overpic}[width=0.847\linewidth,clip]{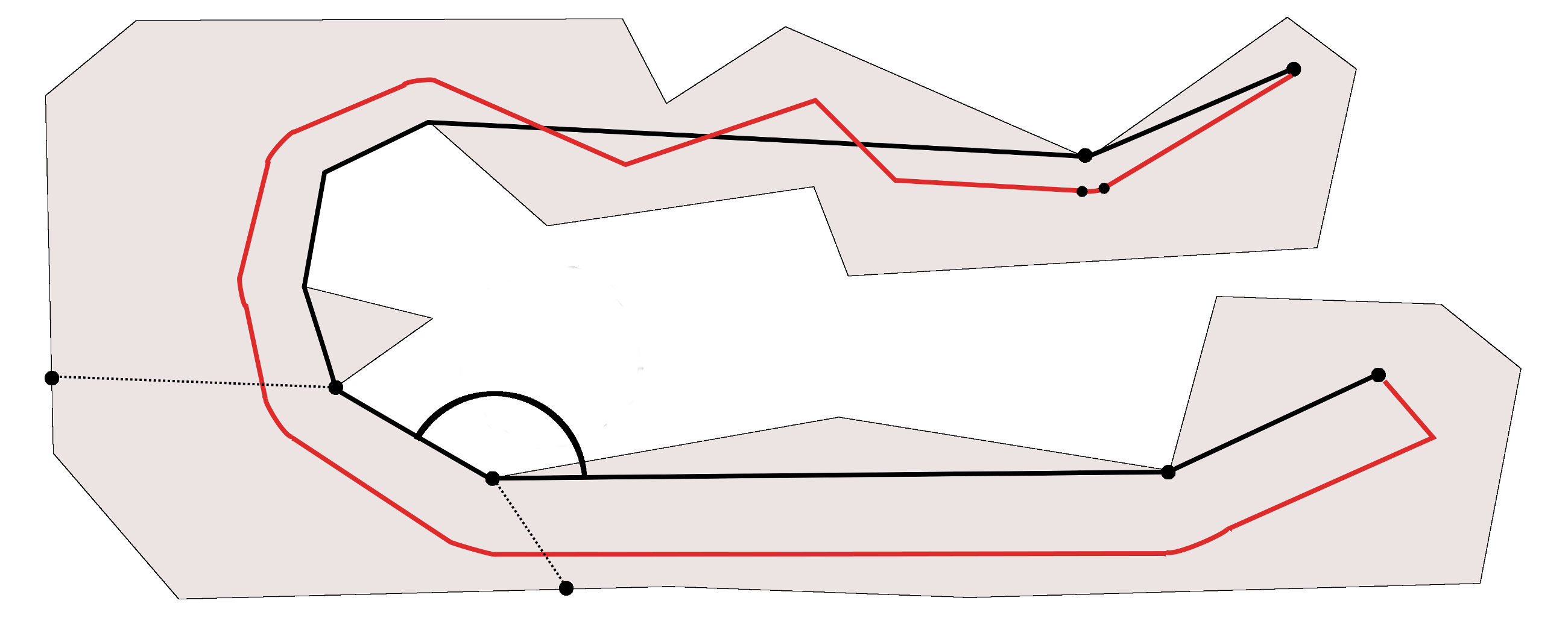}
\put(292,125){{$\hat{v}_0=x$}}
\put(67,14){{$\Gamma_2^{\rm I_1}$}}
\put(81,18) {\line(1,0){16}}
\put(38,86){{$\Gamma_1^{\rm III}$}}
\put(120,112){{$\Gamma_1^{\rm II}$}}
\put(233,107){{$\hat{v}_1$}}
\put(252,84){{$w^-_1$}}
\put(214,83){{$w^+_1$}}
\put(241,95) {\line(1,-1){10}}
\put(237,94) {\line(-1,-1){10}}
\put(102,21){{$\hat{v}_2$}}
\put(128,1){{$y$}}
\put(80,48){{$v$}}
\put(0,52){{$q$}}
\put(107,36){{$S'$}}
\put(251,23){{$\hat{v}_3$}}
\put(290,58){{$\hat{v}_4=p$}}
\end{overpic}
\caption{ The path $\gamma_0$ is depicted in black and the curve $\gamma$ in red, which `changes the side of the boundary' on $\Gamma^{\rm II}_1$ and has a signed curvature on $\Gamma^{\rm III}_1$. The segments $[v;q]$, $[\hat{v}_2;y]$ separate $x$ and $p$. The segment $[\hat{v}_2; \hat{v}_3]$ does not separate $x$ and $p$, but satisfies \eqref{eq: seg len}. We have also illustrated a circular sector $S'$ contained in the cone $S$ and the part of the circle $\Gamma_2^{\rm I_1}$ defined in Step 2. }

  \label{john1}
\end{figure}

\noindent 
Choose the (unique) curve $\gamma_0:[0,l({\gamma_0})] \to P$ with $\gamma_0(0) = x$ and $\gamma_0(l(\gamma_0))= p$ such that $d_P(x,p) = l({\gamma_0})$ (see Figure \ref{john1}). As observed in Section \ref{sec: prep} there are $0 = t_0 < t_1 < \ldots < t_n = {l(\gamma_0)}$ such that $\gamma_0$ is piecewise affine on $[t_{i}, t_{i+1}]$ and $v_i := \gamma_0(t_i) \in \mathcal{V}_P'$ are concave vertices for $i=1, \ldots, n-1$. Moreover, define $v_0 = x$ and $v_n = p$.   We consider a concave vertex $v \in \mathcal{V}_P'$ and  $q \in \partial P$ such that $[v;q]$ induces a partition $P = Q^{(v,q)}_1 \cup Q^{(v,q)}_2$ according to Definition \ref{def: induce} with $x \in Q^{(v,q)}_1$ and $p \in Q^{(v,q)}_2$. For convenience we will call such a segment $[v;q]$  in the following a segment  which \emph{separates $x$ and $p$} (cf. Figure \ref{john1}). Since $P$ is semiconvex, we have
\begin{align}\label{eq: semico-pro}
|[v; q]| \ge \vartheta \min_{k=1,2} d(Q^{(v,q)}_k)\ge \vartheta \min\lbrace  \max\lbrace d_P(v,x), d_P(q,x)\rbrace, \  \omega d(P) \rbrace,
\end{align}
where we used that $d(Q^{(v,q)}_2) \ge r \ge \omega d(P)$.   In particular, if $v=v_i$ we note that $d_P(v,x) = t_i$ and thus for $\vartheta$ small with respect to $\omega$
\begin{align}\label{eq: semico-pro2}
|[v_i; q]| \ge  \omega \vartheta  t_i \ge 4\vartheta^2 t_i.
\end{align}
Likewise, if $v=v_i$, $q = v_{i+1}$ and $[v_i; v_{i+1}]$ separates $x$ and $p$, we find 
\begin{align}\label{eq: semico-pro2-2}
|[v_i; v_{i+1}]|  \ge 4\vartheta^2 t_{i+1}.
\end{align}
Consider the subset  
$$\lbrace 0 = i_0, 1 = i_1, i_2, \ldots, i_{m} = n-1 \rbrace \subset \lbrace 0,\ldots,n-1 \rbrace$$
with corresponding vertices $\hat{v}_j := \gamma_0(t_{i_j})$, $\hat{v}'_j := \gamma_0(t_{i_j+1})$, $j=0,\ldots,m$  such that for $j=1,\ldots, m-1$ the segments $[\hat{v}_j; \hat{v}'_j]$     separate $x$ and $p$   or satisfy
\begin{align}\label{eq: seg len}
|[\hat{v}_j;\hat{v}_j']| = t_{i_j+1} - t_{i_j} \ge 4\vartheta^2 t_{i_j+1}.
\end{align} 
(Observe that the first and the last segment $[\hat{v}_0; \hat{v}'_0]$ and $[\hat{v}_m; \hat{v}'_m]$ do  not induce a partition.) Note that $i_j + 1= i_{j+1}$ is possible, namely if a pair of directly consecutive segments separate $x$ and $p$ or satisfy \eqref{eq: seg len}. We then obtain
\begin{align}\label{eq: seg len2}
|[\hat{v}_j;\hat{v}_j']| = t_{i_j+1} - t_{i_j} \ge 4\vartheta^2 t_{i_j+1}
\end{align}
for all  $j=0,\ldots,m$. If \eqref{eq: seg len} holds, this follows directly. For $j \in \lbrace 0,m\rbrace$ we observe $t_{0} = 0$ and $t_{n} - t_{n-1} \ge r \ge \omega d(P) \ge 4\vartheta^2 d(P)$ for $\vartheta$ small with respect to $\omega$. Otherwise, $[\hat{v}_j; \hat{v}'_j]$ separates $x$ and $p$ and the assertion follows from \eqref{eq: semico-pro2-2} with $i=i_j$, i.e. $v_i = \hat{v}_j$ and $v_{i+1}=\hat{v}_j'$. This property will essentially be important to estimate the length of the curve $\gamma$ defined in Step 6 (cf. \eqref{eq: carC} below).

For each $1 \le i \le n-1$ choose the unique $\nu^-_i,\nu^+_i \in S^1 = \lbrace x \in \R^2: |x|=1 \rbrace$ such that $\nu^-_i \bot v_{i} -v_{i-1}$, $\nu^+_i \bot v_{i+1} -v_{i}$ and $v_i + \eps \nu^\pm_i \in P$ for $\eps>0$ small. Define
\begin{align}\label{eq: point def}
& w_i^- =  v_i + 2\vartheta^2 t_{i} \nu^-_{i}, \  \ \  \ \ \ \ \ w_i^+ =  v_i + 2\vartheta^2 t_{i} \nu^+_{i}.
\end{align}
Moreover, we set $w_0^- = w_0^+ = x$ and $w_n^- = p + 2\vartheta^2 t_n \nu^+_{n-1}$. The goal is to construct a curve $\gamma: [0,l(\gamma)] \to P$ with $\gamma(0) = x$, $\gamma(l(\gamma)) = p$ with ${\rm car}(\gamma,\vartheta^3) \subset P$, where we essentially connect the points $w^\pm_i$ defined above.  We have to construct curves between 
\begin{align}\label{eq:4,5}
\hspace{-.2cm}
{\rm (I)} \
 w^-_{i_j} \text{ and } w^+_{i_j}, \,  w^-_{i_j+1} \text{ and } w^+_{i_j+1}, \ \,    {\rm (II)} \  w^+_{i_j} \text{ and } w^-_{i_j+1}, \  \, {\rm (III)} \ w^+_{i_j+1} \text{ and } w^-_{i_{j+1}}
\end{align} 
for $j=0,\ldots,m-1$ and at the  end of the curve a path between  $w^-_{n-1}$ and $p$.

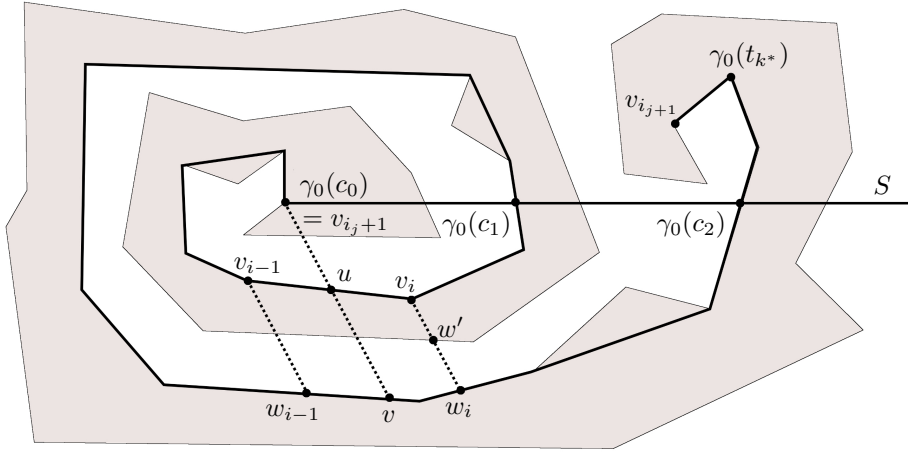
\begin{figure}[H]
\centering
\begin{overpic}[width=0.854\linewidth,clip]{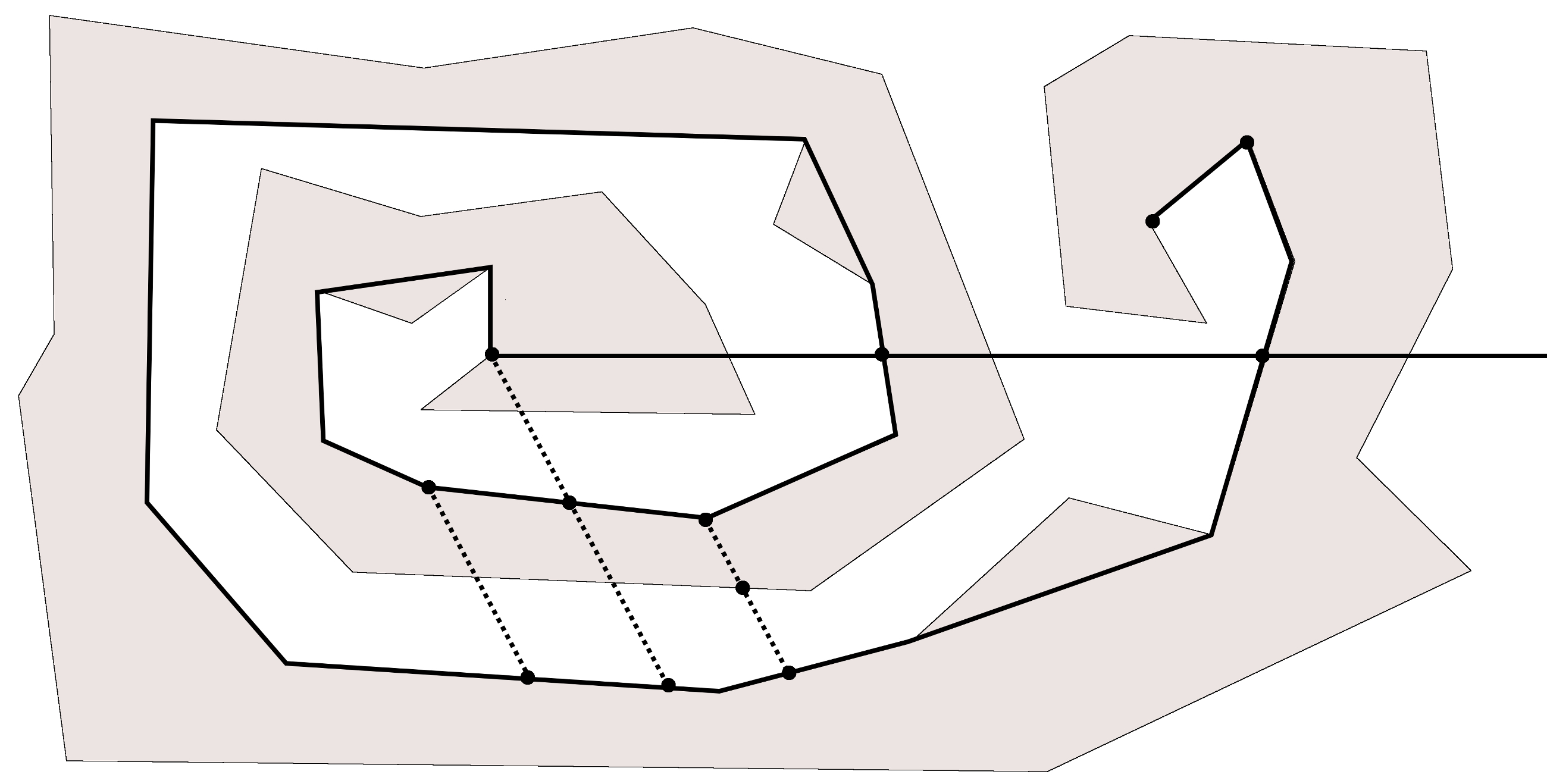}
\put(115,100){{$\gamma_0(c_0)$}}
\put(116,87){{$ = v_{i_j+1}$}}

\put(168.8,85){{$\gamma_0(c_1)$}}
\put(250,85){{$\gamma_0(c_2)$}}
\put(332,99){{$S$}}
\put(270,149){{$\gamma_0(t_{k^*})$}}
\put(151,63.6){{$v_i$}}
\put(167,45){{$w'$}}
\put(90,70.6){{$v_{i-1}$}}

\put(129,66){{$u$}}

\put(170,16){{$w_i$}}
\put(102,15){{$w_{i-1}$}}

\put(146,12.8){{$v$}}
\put(238,132){{$v_{i_{j+1}}$}}

\end{overpic}
\caption{The part of $\gamma_0$ between $t_{i_j+1}$ and $t_{i_{j+1}}$ has been depicted in black (for illustration purposes with $m=2$) and we sketched a segment $[u;v]$ as considered below \eqref{eq:: for k}. Note that the outward helix ends in $\gamma_0(t_{k^*})$. }

  \label{john3}
\end{figure}

The most delicate cases are (II) and (III), where in (II) $\gamma_0$ typically `changes the side of the boundary' and in (III) the part of $\gamma_0$ has signed curvature, possibly  the form of a `helix' (cf. Figure \ref{john1}, Figure \ref{john3}). In Step 2--Step 5 we construct the various parts of the curve, where one has to ensure that (1) the length of $\gamma$ is comparable to the length of $\gamma_0$ (see  \eqref{eq: (I)}(i), \eqref{eq: length gamma}, \eqref{eq: length gamma2}, \eqref{eq: length4}(i)) and (2) the distance of $\gamma$ from the boundary is sufficiently large (see \eqref{eq: (I)}(ii),  \eqref{eq: dist2}, \eqref{eq: dist3}, \eqref{eq: length4}(ii)). In Step 6 we finally show that the constructed curve satisfies the property stated in Definition  \ref{def: John}.  \smallskip\\

\vspace{-0.5cm}

\vspace{-0.5cm}

\noindent  \smallskip \\
\textit{Step 2: Construction of curves {\rm (I)}}

\noindent 
Let $j \in \lbrace 0,\ldots,m-1 \rbrace$ and recall \eqref{eq: point def}-\eqref{eq:4,5}. Let $\Gamma^{\rm I_1}_{j}$ and ${\Gamma}^{\rm I_2}_{j}$ be the parts of the two circles with midpoints $v_{i_j} = \hat{v}_j$, $v_{i_j+1} = \hat{v}_j'$ and radii $2\vartheta^2 t_{i_j}$, $2\vartheta^2 t_{i_j+1}$ connecting $w^-_{i_j}, w^+_{i_j}$ and $w^-_{i_j+1}, w^+_{i_j+1}$  respectively. (Note that $\Gamma^{\rm I_1}_0 = \emptyset$.) We have  
\begin{align}\label{eq: (I)}
\begin{split}
(i)& \ \  l({\Gamma^{\rm I_1}_{j}}) \le 4\pi\vartheta^2 t_{i_j} \le \pi(t_{i_j+1}-t_{i_j}), \ \ \ l({{\Gamma}^{\rm I_2}_{j}}) \le 4\pi\vartheta^2 t_{i_j+1}\le \pi(t_{i_j+1}-t_{i_j}),\\
(ii) & \  \ 
\dist(\partial P,\Gamma^{\rm I_1}_{j}) \ge  \vartheta^2 t_{i_j}, \ \ \ \ \ \dist(\partial P,{\Gamma}^{\rm I_2}_{j}) \ge  \vartheta^2 t_{i_j+1}.
\end{split}
\end{align}
Indeed, the first inequality in (i) is clear and the second follows from \eqref{eq: seg len2}. We show (ii) for $i=i_j$, the proof for $i_j+1$ is similar. Let $\varphi_- = {\rm arg}(v_{i-1}-v_{i})$, $\varphi_+ = {\rm arg}(v_{i+1}-v_{i})$ and suppose that possibly after rotation and reflection we have $0 \le \varphi_+ <\varphi_- < 2\pi$ and $\varphi_- - \varphi_+ < \pi$. We define the (infinite) cone $S=\lbrace x \in \R^2: {\rm arg}(x-v) \in [\varphi_+,\varphi_-]\rbrace$ and note that $\dist({\Gamma^{\rm I_1}_{j}}, S) \ge 2\vartheta^2 t_{i_j}$ by \eqref{eq: point def}. If (ii) was wrong, we would find some $y \in \partial P \setminus S$ such that $[y; \hat{v}_j] \subset P$ and $|[y;\hat{v}_j]| < 3\vartheta^2 t_{i_j}$. As $y \notin S$ and  $\gamma_0$ is the shortest path between $x$ and $p$,  we get that $[y; \hat{v}_j]$  separates $x$ and $p$  (cf. Figure \ref{john1}). This contradicts \eqref{eq: semico-pro2}.\smallskip\\

\noindent  \smallskip
\textit{Step 3: Construction of curves {\rm (II)}}

\noindent 
Let $j \in \lbrace 0,\ldots,m \rbrace$ and recall \eqref{eq: point def}-\eqref{eq:4,5}. (Because of Step 5 below we also consider $j=m$.) To simplify the notation we assume $v_{i_j} = \hat{v}_j= 0 $, $v_{i_j+1} = \hat{v}'_j = (0,d)$, where $d := t_{i_j+1} - t_{i_j}$, and define the rectangle  $R = [-2\vartheta^2 t_{i_j+1}, 2\vartheta^2 t_{i_j+1}] \times [0,d]$. Observe that  $w^+_{i_j},w^-_{i_j+1} \in  \partial R$. For notational convenience we will write $w = w^+_{i_j}$ and $w' = w^-_{i_j+1}$ in the following. Recall that $\hat{v}_j, \hat{v}'_{j} \in \mathcal{V}_P'$ for $j \in \lbrace 1,\ldots,m-1\rbrace$. In the other cases we have (recall $t_0 = 0$, $v_n = p$) 
\begin{align}\label{eq: extra case}
\begin{split}
&B(v_0, t_0) =  B(\hat{v}_0, t_0) = \emptyset \subset P, \  \ \  v_1= \hat{v}_0' \in \mathcal{V}_P', \\ 
&v_{n-1} = \hat{v}_{m} \in \mathcal{V}_P', \ \ \  B(v_n, \omega d(P)) = B(\hat{v}'_{m}, \omega d(P)) \subset P. 
\end{split}
\end{align}
Define the set of vertices $\mathcal{U}_R := \lbrace v \in \mathcal{V}_P': v \in R\rbrace \cup \lbrace \hat{v}_j,\hat{v}_j' \rbrace$. For convenience we now first treat the case $j \in \lbrace 1,\ldots,m-1 \rbrace$ and indicate the minor adaptions for $j \in \lbrace 0,m\rbrace$,  necessary due to $\hat{v}_0,\hat{v}_m' \notin \mathcal{V}_P'$, at the end of Step 3.

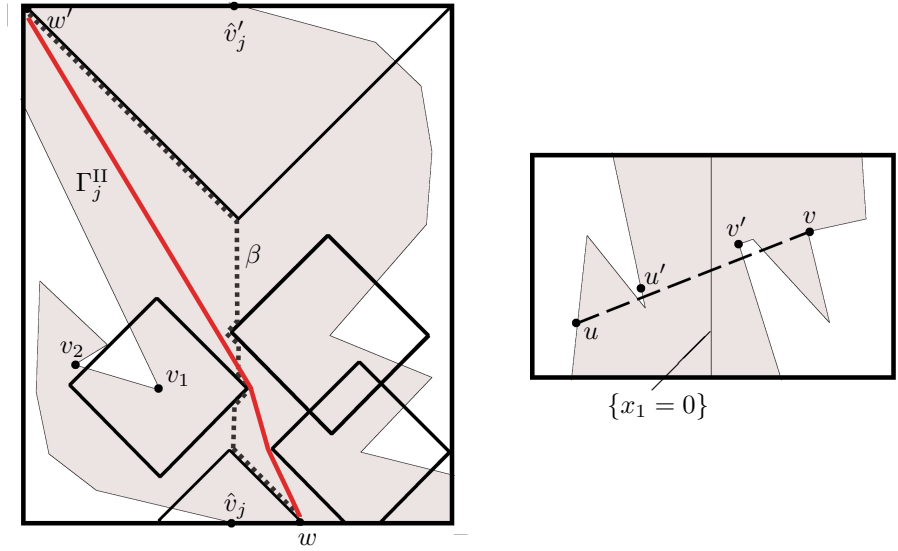
\begin{figure}[H]
\centering
\begin{overpic}[width=0.853\linewidth,clip]{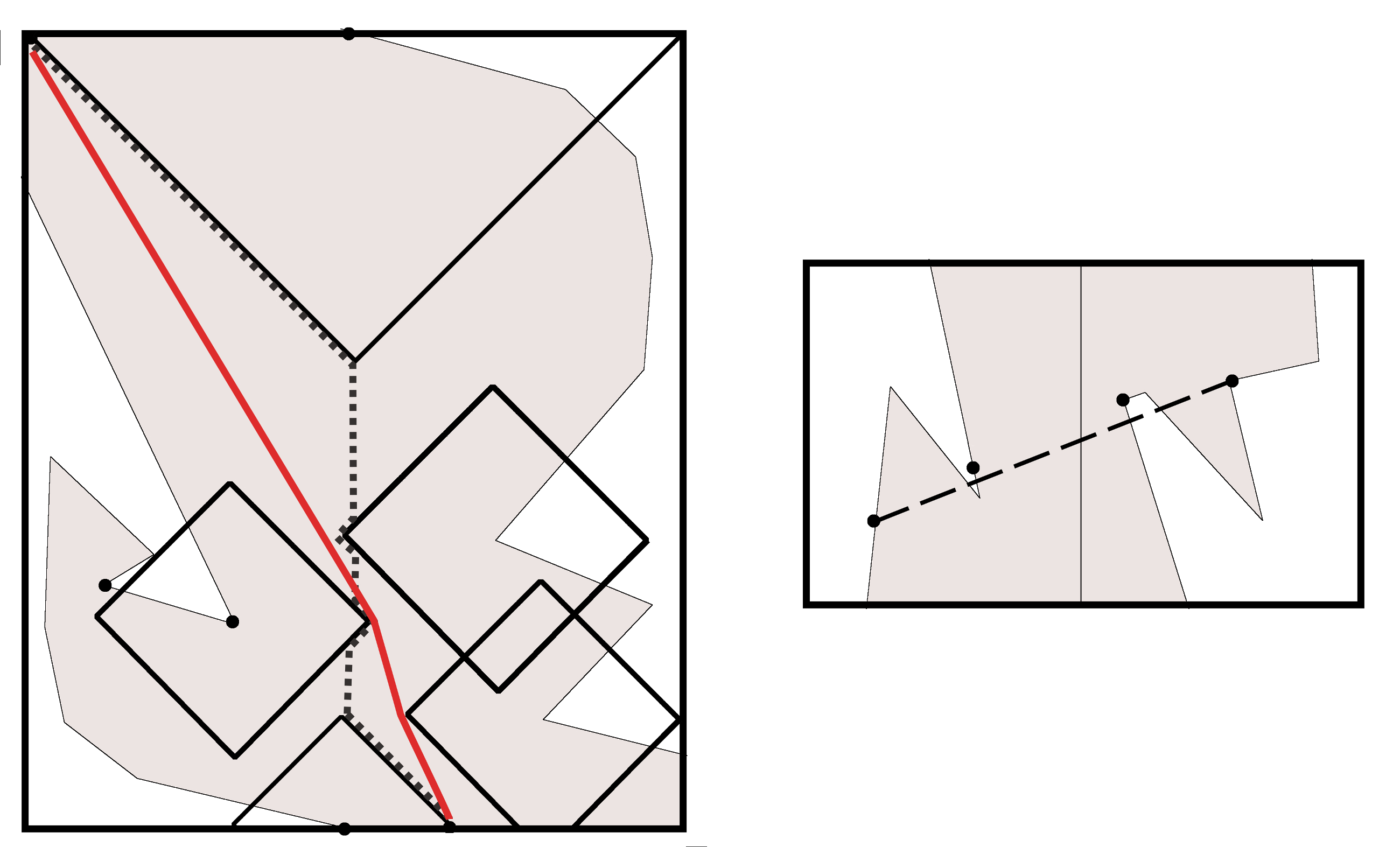}
\put(82,9){{$\hat{v}_j$}}
\put(110,-4){{$w$}}
\put(82,187){{$\hat{v}'_j$}}
\put(15,192){{$w'$}}
 
\put(26,130){{$\Gamma_j^{\rm II}$}}
\put(90,103){{$\beta$}}

\put(60,58){{$v_1$}}

\put(20,69){{$v_2$}}

\put(301,118){{$v$}}

\put(272,113){{$v'$}}

\put(218,74){{$u$}}

\put(241,95){{$u'$}}

\put(244,55) {\line(1,1){22}}
\put(227,46){{$\lbrace x_1 = 0\rbrace$}}

\end{overpic}
\caption{ On the left we have depicted $\beta$ in dotted lines and ${\Gamma^{\rm II}_j}$ in red. Note ${\rm sgn}(\hat{v}_j \e_1) = -1$, ${\rm sgn}(\hat{v}'_j \e_1) = 1$, $C_- = C(\hat{v}_j) \cup C(v_1)$ and $v_2 \notin \mathcal{V}_-$. The main idea in the construction is that ${\Gamma^{\rm II}_j}$ is `not too close to concave vertices'. On the right the situation of \eqref{eq: 3cond} is illustrated.}

  \label{john2}
\end{figure}

Note that $x \e_1 \neq 0$ for all $x \in (\partial P \cap R) \setminus \lbrace \hat{v}_j,\hat{v}_j' \rbrace$ as $[\hat{v}_j;\hat{v}_j']$ induces a partition of $P$.   Let ${\rm sgn}(y) = 1$ for $y>0$  and ${\rm sgn}(y) = -1$ for $y <0$. By convention we set ${\rm sgn}(\hat{v}_j \e_1) = - {\rm sgn}(w\e_1)$ and  ${\rm sgn}(\hat{v}'_j\e_1) = - {\rm sgn}(w'\e_1)$.  We let 
\begin{align}\label{eq: Vpm}
\mathcal{V}_\pm = \lbrace v \in \mathcal{U}_R: \pm{\rm sgn}(v \e_1) >0, \ [v; (0,v\e_2)] \subset P \rbrace 
\end{align}
 and show that 
\begin{align}\label{eq:collecte}
|[v;u]| \ge 8\vartheta^2 (t_{i_j}+v\e_2) \ \  \text{ for all } v\in \mathcal{V}_- \cup \mathcal{V}_+, u \in \partial P \text{ with } {\rm sgn}(v\e_1) \neq {\rm sgn}(u\e_1).
\end{align}
To see this, assume e.g. that $v \in \mathcal{V}_+$ and   suppose first $u\e_2 < v\e_2$. Clearly, $[v; u]$ does not necessarily induce a partition of $P$ as possibly $[v; u] \not\subset P$. However, due to the fact that $\lbrace 0 \rbrace \times [0,d], [v;(0,v\e_2)] \subset P$, we see that there have to exist  $v'\in \mathcal{V}_+$ and $u' \in \partial P$ with
\begin{align}\label{eq: 3cond}
0 \le v'\e_1  \le v\e_1, \ \  \,   \frac{v' \e_2 - u\e_2}{v\e_2-u\e_2} \ge \frac{|[v';u]|}{|[v;u]|}, \ \ \,  v'\e_2  \le v\e_2,  \ \ \, u'\e_1 < 0, \ \ \,  |[v';u']| \le |[v';u]|
\end{align} 
 such that $[v';u']$ induces a partition of $P$ (see Figure \ref{john2}). In fact,  choose $v'$ as a concave vertex in $[0,v\e_1] \times [u\e_2,v\e_2]$ lying on or above the segment $[v;u]$ with minimal distance to $\lbrace 0 \rbrace \times [0,d]$  (note that possibly $v'=v$) and let $u'$ be the point on $\partial P \cap [v';u]$ closest to $v'$. The second property in \eqref{eq: 3cond} follows from the fact that $v'$ lies on or above the segment $[v;u]$.  Since ${\rm sgn}(v'\e_1) \neq {\rm sgn}(u'\e_1)$,  $[v';u']$  separates $x$ and $p$. Now suppose the statement was wrong. We then obtain using    $v\e_2 > u \e_2 \ge 0$, $v' \e_2 \le v \e_2$  as well as   \eqref{eq: 3cond}
$$|[v'; u']| \le  |[v'; u]| \le |[v;u]|\frac{v' \e_2 - u\e_2}{v\e_2-u\e_2} < 8\vartheta^2 (t_{i_j} +  v\e_2) \frac{v' \e_2}{v\e_2} \le  8\vartheta^2 (t_{i_j} +  v'\e_2).
$$
Consequently, we have $|[v';u']| < \vartheta \omega  d(P)$ for $\vartheta$ small with respect to $\omega$. Moreover,  for $\vartheta$ small  we get $2\vartheta^{-1} \le (8\vartheta^2)^{-1}$ and thus by $|v'\e_1| \le |[v'; u']|$
\begin{align*}
\begin{split}
\max\lbrace d_P(v',x), d_P(u',x)\rbrace  &\ge  t_{i_j} +  v'\e_2-   |v'\e_1| > \big(2\vartheta^{-1}-1\big)|[v';u']|   \ge \vartheta^{-1}|[v';u']|. 
\end{split}
\end{align*}
The last two estimates contradict  \eqref{eq: semico-pro}. This shows \eqref{eq:collecte} in the case $u\e_2 < v\e_1$. For $u\e_2 \ge v\e_1$ we proceed similarly, where the second and third property in \eqref{eq: 3cond} are replaced by $v' \e_2 \ge v\e_2$ and $|[v';u]| \le |[v;u]|$.

 Recalling \eqref{eq: Vpm} we let $C(v)$ be the closed square with midpoint $v \in \mathcal{V}_- \cup \mathcal{V}_+$ and diagonal $4\vartheta^2 l(v)$ with faces parallel to $\e_1 + \e_2$ and $\e_1 - \e_2$, where 
 \begin{align}\label{eq: point def2}
l(v) = t_{i_j}+ v \e_2 \in [t_{i_j}, t_{i_j+1}].
\end{align} 
 Moreover, define $C_\pm = \bigcup_{v \in \mathcal{V}_\pm } C(v)$ and let $H_+$, $H_-$ be the closed half space right and left of $\lbrace 0 \rbrace \times \R$, respectively. We  show
\begin{align}\label{eq: particular}
\begin{split}
&(i) \ \ C_{+} \cap C_{-} = \emptyset, \ \ \ \ \ (ii) \ \ (C_{+} \cap H_-) \cup (C_{-} \cap H_{+}) \subset P,\\
& (iii) \ \  w,w' \in (\partial C_{+} \cap H_-) \cup (\partial C_{-} \cap H_{+}).
\end{split}
\end{align}
To see  (i), note that \eqref{eq:collecte} implies $|[v_1; v_2]| \ge 8\vartheta^2 \max\lbrace l(v_1),l(v_2) \rbrace$ for  $v_1\in \mathcal{V}_+, v_2 \in \mathcal{V}_-$  and thus $C(v_1) \cap C(v_2) = \emptyset$. Likewise,  if  (ii)  was wrong, there would be, e.g., $v \in \mathcal{V}_+$ and $u \in \partial P$ with ${\rm sgn}(u \e_1) < 0$ such that $|[v; u]| \le 2\vartheta^2 l(v)$, which contradicts \eqref{eq:collecte}. Finally, we always have  $C(\hat{v}_j) \cap C(\hat{v}_j') = \emptyset$ by \eqref{eq: seg len2}. Consequently, \eqref{eq: point def}, \eqref{eq: point def2}, and the convention ${\rm sgn}(\hat{v}_j \e_1) = - {\rm sgn}(w\e_1)$, ${\rm sgn}(\hat{v}'_j\e_1) = - {\rm sgn}(w'\e_1)$  show that  each of the points  $w,w'$  is contained in  $\partial C_+ \cap H_-$ or $\partial C_- \cap H_+$.  \smallskip

\vspace{-0.40cm}

We define 
$P_R= \overline{(P\cap R) \setminus (C_{+} \cup C_{-})}$. Then $w,w' \in P_R$ by  \eqref{eq: particular}(iii). Moreover, by \eqref{eq: particular} we find a   continuous, piecewise affine  path $\beta$ between $w$,  $w'$ in the set 
$$\big(\lbrace 0 \rbrace \times [0,d] \cup  (\partial C_+ \cap H_-) \cup (\partial C_- \cap H_+) \big) \cap P_R$$
such that the tangent vector of $\beta$ is a.e. contained in $\lbrace \frac{1}{\sqrt{2}} (-1,1), (0,1), \frac{1}{\sqrt{2}} (1,1) \rbrace$ (see Figure \ref{john2}). In particular, $w$, $w'$ are in the same connected component of $P_R$, which will be denoted by $P_R^{\rm con}$ in the following. Let $\Gamma^{\rm II}_j:[0,l({\Gamma^{\rm II}_j})] \to P_R^{\rm con}$ be the shortest curve between $w$ and $w'$ parametrized by arc length. The goal will be to establish \eqref{eq: length gamma} and \eqref{eq: dist2} below. Since $P_R^{\rm con}$ is a polygon, we find that  $\Gamma^{\rm II}_j$ is piecewise affine and changes its direction only in concave vertices of $P_R^{\rm con}$. We    show that for all $0\le s < s' \le l({\Gamma^{\rm II}_j})$  one has 
 \begin{align}\label{eq: length gamma***}
 (i) \ \ [\Gamma^{\rm II}_j(s);  (0,\Gamma^{\rm II}_j(s) \e_2)] \subset P, \ \ \ \ \ \ (ii) \ \ {\rm arg}(\Gamma^{\rm II}_j(s') - \Gamma^{\rm II}_j(s)) \in \big[\frac{\pi}{4},\frac{3\pi}{4}\big].
 \end{align}
 First, \eqref{eq: length gamma***}(i) holds for $\beta$ in place of $\Gamma^{\rm II}_j$ by \eqref{eq: particular}. Since $P_R^{\rm con}$ is simply connected, the bounded connected components of $\R^2\setminus(\beta \cup \Gamma^{\rm II}_j)$ are contained in $P_R^{\rm con}$. Then the fact that $\Gamma^{\rm II}_j$ is the shortest path between $w$, $w'$ in $P_R^{\rm con}$ together with \eqref{eq: length gamma***}(i) for $\beta$ implies \eqref{eq: length gamma***}(i) for $\Gamma^{\rm II}_j$.

From \eqref{eq: length gamma***}(i) we deduce that ${\rm arg}(\Gamma^{\rm II}_j(s') - \Gamma^{\rm II}_j(s)) \in [0,\pi]$  for all $s < s'$ since $\Gamma^{\rm II}_j$ is the shortest path between $w$, $w'$. Select $s_1 < s_2$ and $t_1 < t_2$ such that $\beta(t_k) = \Gamma^{\rm II}_j(s_k)$ for $k=1,2$ and $\Gamma^{\rm II}_j((s_1,s_2)) \cap \beta((t_1,t_2)) = \emptyset$. Let $P_*$ be the polygon with boundary $\Gamma^{\rm II}_j([s_1,s_2]) \cup \beta([t_1,t_2])$. Since $P_* \subset P_R^{\rm con}$ and  $\Gamma^{\rm II}_j$ is the shortest path between $w,w'$ in $P_R^{\rm con}$, $P_*$ only has concave vertices on $\Gamma^{\rm II}_j((s_1,s_2))$. Recalling that the tangent vector of $\beta$ is a.e. contained in  $\lbrace \frac{1}{\sqrt{2}} (-1,1), (0,1), \frac{1}{\sqrt{2}} (1,1) \rbrace$ and that the paths $\Gamma^{\rm II}_j([s_1,s_2])$, $\beta([t_1,t_2])$ have a common start and endpoint,  we derive \eqref{eq: length gamma***}(ii) for $s_1 \le s < s'  \le s_2$. Herefrom we also deduce  
 \begin{align}\label{eq: length gamma}
l({\Gamma^{\rm II}_j}) \le \sqrt{2}(t_{i_j+1} - t_{i_j}).
\end{align}
Additionally, we obtain that if $ \Gamma^{\rm II}_j$ changes its direction in $s$, then
\begin{align}\label{eq:changi}
 \Gamma^{\rm II}_j(s)  = v + 2\vartheta^2 l(v)\e_1  \ \text{ for } v \in \mathcal{V}_- \ \ \ \  \text{or} \ \  \ \  \Gamma^{\rm II}_j(s)  = v - 2\vartheta^2 l(v)\e_1  \ \text{ for } v \in \mathcal{V}_+.
\end{align}
In fact, $\Gamma^{\rm II}_j$ changes its direction only in concave vertices of $P_R^{\rm con}$. First,  if $ \Gamma^{\rm II}_j(s) \in \mathcal{U}_R$ (recall definition below \eqref{eq: extra case}), then $\Gamma^{\rm II}_j(s) \in \mathcal{V}_- \cup \mathcal{V}_+$ by \eqref{eq: length gamma***}(i) and thus $\Gamma^{\rm II}_j(s) \notin P_R$ since $C(\Gamma^{\rm II}_j(s)) \cap P_R =  \emptyset$. This gives a contradiction. Consequently, $\Gamma^{\rm II}_j(s)$ is a corner of $C(v)$ for some $v \in \mathcal{V}_- \cup \mathcal{V}_+$. Then \eqref{eq: length gamma***} together with the geometry of $C(v)$ and the fact that $\Gamma^{\rm II}_j$ is a shortest path implies that $\Gamma^{\rm II}_j(s)$ has to be the left or right corner of $C(v)$, respectively. This yields \eqref{eq:changi}. We now finally show
\begin{align}\label{eq: dist2}
\dist(\partial P, \Gamma^{\rm II}_j(s)) \ge c\vartheta^2(t_{i_j} + s) \ \ \ \text{for} \ s\in [0,l(\Gamma^{\rm II}_j)]
\end{align} 
for some universal $c >0$  small. First, in view of \eqref{eq: (I)} we observe that \eqref{eq: dist2} holds for $s=0, l({\Gamma^{\rm II}_j})$ since $\hat{v}_j = \Gamma^{\rm II}_j(0) \in \Gamma^{{\rm I}_1}_{j}$ and $\hat{v}_j' = \Gamma^{\rm II}_j(l({\Gamma^{\rm II}_j})) \in {\Gamma}^{\rm I_2}_{j}$. For each $s \in [0,l({\Gamma^{\rm II}_j})]$ we denote by $q^\pm(s)$ the nearest point to $\Gamma^{\rm II}_j(s)$ on $\partial P \cap \big(\Gamma^{\rm II}_j(s) \pm \R_+ \e_1\big)$. Moreover, we set $f^\pm(s) = |q^\pm(s)- \Gamma^{\rm II}_j(s)|$. For later note that $f^\pm$ is a lower semicontinuous function and it is possibly discontinuous at $s$  only if $q^\pm(s)$ is a concave vertex. The fact that \eqref{eq: dist2} holds for $s=0, l({\Gamma^{\rm II}_j})$ and \eqref{eq: length gamma***}(ii) show that it suffices to prove  
\begin{align}\label{eq: dist2***}
f^\pm(s) > \vartheta^2(t_{i_j} + s)
\end{align} 
for  $s\in [0,l(\Gamma^{\rm II}_j)]$ as herefrom \eqref{eq: dist2} follows for $c>0$ sufficiently small. Consider, e.g., $f^+$. First, we show that \eqref{eq: dist2***} holds for $s$ with 
$$(a) \ \ q^+(s) \in \mathcal{U}_R \ \ \ \ \text{ or } \ \ \ \ (b) \ \ \Gamma^{\rm II}_j \text{ changes its direction in } s.$$
In fact, in case (a), \eqref{eq: Vpm}, \eqref{eq: length gamma***}(i) imply $q^+(s) \in \mathcal{V}_+$ and thus by \eqref{eq: point def2}  and $\Gamma^{\rm II}_j \subset P_R$ we  have $f^+(s) \ge 2 \vartheta^2(t_{i_j} + \Gamma^{\rm II}_j(s)\e_2)$.  If (a) does not hold, we consider case (b) and  recall \eqref{eq:changi}. If we had $\Gamma^{\rm II}_j(s) = v - 2\vartheta^2 l(v) \e_1$ for some $v \in \mathcal{V}_+$, (a) would be satisfied since then $q_+(s) = v \in \mathcal{U}_R$.   Consequently, we have $\Gamma^{\rm II}_j(s) = v + 2\vartheta^2 l(v) \e_1$ for some $v \in \mathcal{V}_-$   and then $f^+(s) \ge 6 \vartheta^2(t_{i_j} + \Gamma^{\rm II}_j(s)\e_2)$ by \eqref{eq:collecte}  as ${\rm sgn}(v \e_1) \neq {\rm sgn}(q^+(s)\e_1)$. In all cases \eqref{eq: dist2***} follows from the fact that $\Gamma^{\rm II}_j(s)\e_2 \ge \frac{s}{\sqrt{2}}$ by \eqref{eq: length gamma***}(ii).

We now show \eqref{eq: dist2***} by contradiction. Choose the largest value $0<s < l({\Gamma^{\rm II}_j})$ such that \eqref{eq: dist2***} is violated. Then 
neither (a) nor (b) hold. Since (a) does not hold, $f^+$ is continuous in a neighborhood of $s$ and thus
\begin{align}\label{eq: dist2***X}
f^+(s) = \vartheta^2(t_{i_j} + s)
\end{align} 
by the choice of $s$. Choose the largest value $s' < s$ such that for $s'$ one of the conditions (a), (b) holds. (If this is not possible, set $s' = 0$.) We now show that \eqref{eq: dist2***X}  implies $f^+(s') \le  \vartheta^2(t_{i_j} + s')$ which  contradicts the fact that \eqref{eq: dist2***} holds for $s'$. This will conclude the proof of \eqref{eq: dist2***} and then \eqref{eq: dist2} is proved. \smallskip

Let $t'= \Gamma^{\rm II}_j(s')\e_2$, $t= \Gamma^{\rm II}_j(s)\e_2$ and $T:[0,d] \to [0,l({\Gamma^{\rm II}_j})]$ be the (increasing) function with $\tau = \Gamma^{\rm II}_j(T(\tau))\e_2$ for $\tau \in [0,d]$. Due to the fact that $\Gamma^{\rm II}_j$ does not change its direction on  $(s',s]$ we observe that $\tau \mapsto T(\tau)$ and $\tau \mapsto \Gamma^{\rm II}_j(T(\tau))\e_1$ are affine on $(t',t+\eps)$ for $\eps$ small enough. Moreover, as (a) does not hold, we get that 
$\tau \mapsto q^+(T(\tau))\e_1$ is concave in $(t',t+\eps)$ (cf. upper part in Figure  \ref{john2}).  Define $g:[0,d] \to [0,\infty)$ 
by 
$$g(\tau) = q^+(T(\tau))\e_1 - \Gamma^{\rm II}_j(T(\tau))\e_1$$
and observe that $g$ is concave in $(t',t+\eps)$.  More precisely, $g$ is   differentiable up to a finite number of points. To avoid further notation involving the superdifferential of concave functions, we will for simplicity assume that $g$ is smooth. In fact, this can be always obtained by a slight modification of $g$ on $(t',t)$  without affecting the following arguments.

Since $g$ is concave and $T$ is affine on $(t',t+\eps)$, we get $\bar{T} >0$ such that 
\begin{align}\label{eq: convex}
g(\tau) \le g(t) +  g'(t)(\tau - t), \ \ \ \ \ T(\tau) = T(t) +  \bar{T}(\tau - t)
\end{align}
 for $\tau \in (t',t]$. The function $h:[0,d] \to [0,\infty)$ defined by  $h(\tau) = g(\tau)(t_{i_j} + T(\tau))^{-1}$ satisfies $h(t) = \vartheta^2$ by \eqref{eq: dist2***X} and $T(t) =s$. Note also that $h'(t) \ge 0$ due to the maximal choice of $s$. Consequently, $(t_{i_j} + T(t)) g'(t) - \bar{T}g(t) \ge 0$ and  this together with \eqref{eq: convex} and $\bar{T}>0$ yields for  $\tau \in (t',t]$
$$g(\tau)\le g(t) + g'(t)(\tau - t) \le   g'(t)\bar{T}^{-1}(t_{i_j} +T(t))  + g'(t)(\tau - t) = g'(t)\bar{T}^{-1}(t_{i_j} +T(\tau)).$$
Using that $g'$ is  non-increasing and $\bar{T}>0$ we then find for $\tau \in (t',t]$
\begin{align*}
h'(\tau) &= (t_{i_j} + T(\tau))^{-2} \big((t_{i_j} + T(\tau))g'(\tau) - \bar{T}g(\tau) \big) \\&
\ge (t_{i_j} + T(\tau))^{-1} (g'(\tau) - g'(t))  \ge 0
\end{align*}
and thus $h(\tau) \le  \vartheta^2$ on $(t',t)$. This yields $f^+(\sigma) \le \vartheta^2(t_{i_j}+\sigma)$ for all $\sigma \in (s',s]$. As $f^+$ is lower semicontinuous, we get the desired contradiction $f^+(s') \le \vartheta^2(t_{i_j}+s')$.

 To conclude Step 3, it remains to treat the cases announced in \eqref{eq: extra case}. First, for $j=0$, \eqref{eq:collecte} trivially holds for $v=v_0 = 0$ and $t_0 =0$ and \eqref{eq: dist2} is true for $s=0$ since $t_0=0$. For $j=m$, \eqref{eq:collecte} follows from  $v=v_n=p$ and \eqref{eq: extra case} with $\vartheta$ small with respect to $\omega$. Finally,  \eqref{eq: dist2} is satisfied for $s= l({\Gamma^{\rm II}_m})$ again by $B(p,\omega d(P)) \subset P$.   The rest remains unchanged.\smallskip\\

\noindent  \smallskip
\textit{Step 4: Construction of curves {\rm (III)}}

\noindent 
Let $j \in \lbrace 0,\ldots,m-1 \rbrace$ and recall \eqref{eq: point def}-\eqref{eq:4,5}. Let us first observe that if $i_{j} + 1= i_{j+1}$, then $w^-_{i_j+1}$ and $w^-_{i_{j + 1}}$  coincide. Therefore, in this particular case we set ${\Gamma}^{\rm I_2}_{j} = \Gamma^{\rm III}_j = \emptyset$. Now suppose $i_{j} + 1< i_{j+1}$. First, as $[v_i;v_{i+1}]$ do not  separate $x$ and $p$ for $i_j +1\le i \le i_{j+1}-1$ (recall definition before \eqref{eq: seg len}),  we see that $\gamma_0([t_{i_j+1}, t_{i_{j+1}}])$ has the form of a helix, i.e. $\gamma_0$ has in $[t_{i_j+1}, t_{i_{j+1}}]$ a signed curvature. (Clearly, a `degenerated helix' with less than a full winding is possible.) More precisely, $\gamma_0$ may consist of an \emph{outward helix} and an \emph{inward helix} in the following sense: define 
\begin{align*}
\varphi_k = {\rm arg}(v_{k+1}-v_{k}) \ \ \ \text{ for } \ i_j+1 \le k \le i_{j+1}-1
\end{align*}
and let $S_k = v_k + \R_+ e^{i\varphi_k}$ with $\R_+ = (0,\infty)$. Let $k^*$ be the smallest index such that
$$S_{k^*} \cap \gamma_0([t_{i_{j}+1},t_{k^*}]) \neq \emptyset$$
and let $\gamma_0([t_{i_{j}+1},t_{k^*}])$, $\gamma_0([t_{k^*},t_{i_{j+1}}))$ be the outward and inward part of the helix, respectively. Indeed, beyond $v_{k^*}$ the helix can not further growth outwardly as this would inavoidably imply self-intersection of the polygon (see Figure \ref{john3}).

Recalling \eqref{eq: point def} we let $\Gamma^{\rm III}_j: [0,l({\Gamma^{\rm III}_j})] \to  \R^2$ be the arc length parametrized curve with $\Gamma^{\rm III}_j(s^\pm_i) = w^\pm_{i_j + i}$ for suitable $0=  s_1^+ < s^-_2 < s^+_2 < \ldots < s^+_{N-1} < s^-_N = l({\Gamma^{\rm III}_j})$, $N=i_{j+1} - i_j$, which is   affine on $[s^+_{i-1}, s^-_i]$ and on $[s_i^-,s_i^+]$ a part of a circle with midpoint $v_{i_j + i}$ and radius $2\vartheta^2 t_{i_j + i}$ (see also Step 2 and Figure \ref{john1}). The crucial point is to show that the length of $\Gamma^{\rm III}_j$ is comparable to $t_{i_{j+1}}- t_{i_j}$ (cf. \eqref{eq: length gamma2} below). To this end, we have to ensure that up to a finite number of `windings' of the helix, the `radius of a winding' can be suitably bounded from below.

We first concentrate on the part $\gamma_0([t_{i_{j}+1},t_{k^*}])$. Possibly after a translation, rotation and reflection we can assume   $S_{i_j+1} = \lbrace 0 \rbrace \times (0,\infty)$ (i.e. $v_{i_j+1} = 0$ and $v_{i_j+2} \in \lbrace 0 \rbrace \times (0,\infty)$) and ${\rm arg}(v_{i_j+3} - v_{i_j+2}) \in (\frac{\pi}{2},\frac{3\pi}{2})$. Let  $S = [0,\infty) \times \lbrace 0 \rbrace$. Let $c_k \in [t_{i_j+1},t_{k^*}]$, $t_{i_j+1} = c_0 < c_1 < c_2 < \ldots < c_m$ be the points for which $\gamma_0(c_k) \in S$. Note that the number of points $\# (c_k)_k$ can be interpreted as the winding number of the outward helix.  We now show for $ 3 \le k \le m-1$
\begin{align}\label{eq:: for k}
\gamma_0([c_{k-1},c_k]) \cap  B(0,r_k) = \emptyset, 
\end{align}
where $B(0,r_k) $ denotes the open ball with radius $r_k = \vartheta^2 c_{k-3}$. If $m \le 3$, there is nothing to show. Therefore, we suppose $m \ge 4$.  Choose an arbitrary $v \in \gamma_0([c_{k-1},c_{k}])$ for $3 \le k \le m-1$. Let  $u \in \gamma_0([c_{k-2},c_{k-1}])$ be the (unique) point on $[0;v]$. In particular, we have $|[v;u]| \le |v|$. Select the index $i$ such that $u \in [v_{i-1}; v_i] \subset \gamma_0([c_{k-3},c_{k}])$. 

Denote the intersection points  on  $(v_{i-1} + \R_+(v-u)) \cap \gamma_0$  and $(v_{i} + \R_+(v-u)) \cap \gamma_0$ nearest to $v_{i-1}$ and $v_i$, respectively, by  $w_{i-1}$ and $w_i$. Due to the geometry of   $\gamma_0([t_{i_j+1},t_{k^*}])$  we have $w_{i-1},w_i \in \gamma_0([c_{k-3},c_{k+1}])$ and $\min_{l=i-1,i} |[v_l;w_l]| \le |[u;v]|$. Suppose, e.g., $|[w_{i};  v_{i}]| \le |[u;  v]| \le |v|$. Then we find some $w' \in \partial P$ with $w' \in [v_{i};w_{i}]$ such that $[v_{i};w'] \subset P$ and $[v_{i};w']$  separates $x$ and $p$.  Now in view of  $|[w';  v_{i}]|  \le |v|$ and $t_{i} \ge c_{k-3}$ (recall $v_{i} \in \gamma_0([c_{k-3},c_k])$),  \eqref{eq: semico-pro2} implies
\begin{align}\label{eq: semico-pro3}
|v| \ge|[v_{i};w']| \ge 4 \vartheta^2 t_{i} \ge \vartheta^2 c_{k-3} = r_k,
\end{align}
which gives \eqref{eq:: for k}. For the part $\gamma_0([t_{k^*}, t_{i_{j+1}}])$  we proceed analogously. Let $\hat{S} = v_{i_{j+1}} + [0, \infty) e^{i\varphi}$ for $\varphi \in [0,2\pi)$ such that $\gamma_0(c_m) \in \hat{S}$. Let $\hat{c}_k \in [c_m,t_{i_{j+1}}]$, $c_m = \hat{c}_0 < \hat{c}_1 < \hat{c}_2 < \ldots < \hat{c}_{\hat{m}} = t_{i_{j+1}}$ be the points for which $\gamma_0(\hat{c}_k) \in \hat{S}$. Similarly as before we can show that for $1\le k \le \hat{m}-2$ one has
\begin{align}\label{eq:: for k2}
\gamma_0([\hat{c}_{k-1},\hat{c}_k]) \cap  B(v_{i_{j+1}}, \hat{r}_{k}) = \emptyset, 
\end{align}
where $\hat{r}_k = \vartheta^2 \hat{c}_{k-1}$. In fact, select some $v \in \gamma_0([\hat{c}_{k-1},\hat{c}_k])$ with  $\hat{m} - k \ge 2$. Let  $u \in \gamma_0([\hat{c}_{k},\hat{c}_{k+1}])$ be the (unique) point on $[v_{i_{j+1}};v]$. In particular, we have $|[v; u]| \le |[v; v_{i_{j+1}}]|$. Assume $u \in [v_{i-1}; v_i]$, where $v_{i-1},v_i \in \gamma_0([\hat{c}_{k-1}, \hat{c}_{k+2}])$. Arguing exactly as before we find some $w' \in \partial P$ such that $[v_{i};w'] \subset P$, $|[w'; v_{i}]|  \le |[v ; v_{i_{j+1}}]|$ and $[v_{i};w']$  separates $x$ and $p$. (Note that as before we possibly have to replace  $v_{i}$ by $v_{i-1}$.) Then as in \eqref{eq: semico-pro3} using $v_{i} \in \gamma_0([\hat{c}_{k-1}, \hat{c}_{k+2}])$ we have $|[v_{i};w']| \ge 4 \vartheta^2 t_{i} \ge \vartheta^2 \hat{c}_{k-1}.$ Consequently, we obtain $|[v;v_{i_{j+1}}]| \ge  \hat{r}_{k}$. 

By \eqref{eq:: for k}, \eqref{eq:: for k2} and the fact that $\gamma_0$ is parametrized by arc length we deduce 
\begin{align}\label{eq: pi}
\begin{split}
&c_k - c_{k-1} \ge 2\pi \vartheta^2 c_{k-3}, \ \ \ \ 3 \le k \le m-1, \\   &\hat{c}_l - \hat{c}_{l-1} \ge 2\pi \vartheta^2 \hat{c}_{l-1},  \ \ \  \ 1 \le l \le \hat{m}-2. 
\end{split}
\end{align}
For $1 \le k \le m$ let $\mathcal{N}_k = \lbrace n \in \N: v_{i_j + n} \in \gamma_0([c_{k-1},c_k]) \rbrace$. By construction of $\Gamma^{\rm III}_j$ and \eqref{eq: point def} we obtain for $\vartheta$ small
\begin{align}\label{eq: pi2}
s_n^- - s_{n-1}^+ &= \sqrt{(t_{i_j +n} - t_{i_j  + n-1})^2 + (2\vartheta^2 (t_{i_j +n} - t_{i_j  + n-1}))^2} \le 2(t_{i_j +n } - t_{i_j  + n-1}),\notag \\
s_n^+ - s_{n}^- &= 2\vartheta^2  t_{i_j +n} \tilde{\varphi}_n,  
\end{align}
where $\tilde{\varphi}_n$ denotes the  angle enclosed by $\nu_{i_j +n}^-$,  $\nu_{i_j +n}^+$ smaller than $\pi$ (recall \eqref{eq: point def}). Let $1 \le k \le m$ and consider $n \in \mathcal{N}_k$. If $5 \le k \le m$, let $n_k$ be the largest index in $\mathcal{N}_{k-4}$. Otherwise, set $n_k = 0$. Using \eqref{eq: pi} we first observe  
\begin{align*}
\sum^{n_k}_{l=1} t_{i_j +l } \, \tilde{\varphi}_{i_j +l } & = \sum^{k-4}_{t=1} \sum_{l \in \mathcal{N}_t} t_{i_j +l } \, \tilde{\varphi}_{i_j +l } \le \sum^{k-4}_{t=1}\sum_{l \in \mathcal{N}_t} c_t \, \tilde{\varphi}_{i_j +l } \le 3\pi\sum^{k-4}_{t=1} c_t \\
&\le \frac{3}{2\vartheta^{2}}\sum^{k-4}_{t=1} (c_{t+3} - c_{t+2}) \le \frac{3}{2\vartheta^2}(c_{k-1} - c_3) \le \frac{3}{2\vartheta^2}(c_{k-1} - c_0),
\end{align*}
where in the third step a calculation yields $\sum_{l \in \mathcal{N}_t}  \, \tilde{\varphi}_{i_j +l } \le 2\pi + 2\frac{\pi}{2} = 3\pi$ for all $t$. Similarly, we one can show $\sum^n_{l=n_k+1} \tilde{\varphi}_{i_j +l } \le 8\pi + 2\frac{\pi}{2}$ and thus  by \eqref{eq: pi2} we derive 
\begin{align*}
s_n^+ &= s_n^+ - s_1^+ = \sum^n_{l=2} s_l^+ - s^+_{l-1} \le 2(t_{i_j +n} - t_{i_j  + 1}) + 2\vartheta^2  \sum^n_{l=2} t_{i_j +l} \, \tilde{\varphi}_{i_j +l}\\
& \le 2 (t_{i_j +n } - t_{i_j  + 1}) + 3(c_{k-1} - c_0)+  18 \pi \vartheta^2 t_{i_{j} + n}.
\end{align*}
Recalling $c_0 = t_{i_j  + 1}$ and $c_{k-1} \le t_{i_j +n }$  since $n \in \mathcal{N}_k$, we then find  
\begin{align}\label{eq: long term}
s_n^+ \le 5 (t_{i_j +n } - t_{i_j  + 1}) + 18 \pi \vartheta^2 t_{i_{j} + n} \le C (t_{i_j +n } - t_{i_j  + 1}) + C \vartheta^2 t_{i_{j} + 1}.
\end{align}
Now letting $\hat{\mathcal{N}}_k = \lbrace n \in \N: v_{i_j+n} \in \gamma_0([\hat{c}_{k-m-1},\hat{c}_{k-m}])\rbrace$ for $m+1 \le k \le \hat{m} + m$ and repeating the above arguments we find for $n \in \hat{\mathcal{N}}_k$ 
\begin{align}\label{eq: long term*}
s_n^+ \le C (t_{i_j +n} - t_{i_j  + 1})  + C \vartheta^2 t_{i_j  + 1},
\end{align}
where we set $s^+_N := s^-_N = l({\Gamma^{\rm III}_j})$.  Thus, in particular for $n = N = i_{j+1} - i_j$ we have by \eqref{eq: seg len2} 
\begin{align}\label{eq: length gamma2}
\begin{split}
l({\Gamma^{\rm III}_j}) &\le C(t_{i_{j+1}} - t_{i_j  + 1}) +  C \vartheta^2 t_{i_j  + 1} \le C(t_{i_{j+1}} - t_{i_j  + 1}) + C(t_{i_{j}+1} - t_{i_j}).
\end{split}
\end{align}
We finally show that for $s \in [0, l(\Gamma^{\rm III}_j)]$ one has
\begin{align}\label{eq: dist3}
\dist(\Gamma^{\rm III}_j(s), \partial P) \ge c\vartheta^2(  t_{i_j+1} + s)
\end{align}
for some $c>0$ sufficiently small. Let $s \in [s^+_{n-1}, s^+_{n}]$ and $P = Q_1 \cup Q_2$  be the partition induced by $[v_{i_j + n-1};v_{i_j + n}]$ with $x,p \in Q_1$ since the segment does not separate $x$ and $p$. (Observe that $Q_2 = \emptyset$ is possible.) As \eqref{eq: seg len} does not hold, we have  $t_{i_j + n} - t_{i_j + n-1} < 4\vartheta^2 t_{i_j + n} = 4\vartheta^2 t_{i_j + n - 1} + 4\vartheta^2 (t_{i_j + n}-t_{i_j + n - 1})$ and then with $\vartheta$ small we get  
$$t_{i_j + n} - t_{i_j + n-1} \le 5\vartheta^2 t_{i_j + n-1} = 5\vartheta^2 (t_{i_j + n-1} - t_{i_j +1}) + 5\vartheta^2  t_{i_j +1}.
$$ 
Consequently, we find by \eqref{eq: long term}, \eqref{eq: long term*} for $C \ge 2$ and $\vartheta$ small (such that $C\vartheta^2\le 1$)
\begin{align}\label{eq: ecco}
\begin{split}
s + t_{i_j+1} &\le  s_{n}^+ + t_{i_j+1} \le C (t_{i_j +n } - t_{i_j +1}) + C\vartheta^2 t_{i_j +1} + t_{i_j+1}  \\& \le C (t_{i_j +n-1 } - t_{i_j +1}) + C\vartheta^2 t_{i_j +1} + t_{i_j+1} \le C (t_{i_j +n-1 } - t_{i_j +1}) + 2t_{i_j+1} \\&
\le C t_{i_j +n-1 }.
\end{split}
\end{align}
Fix $u \in \partial P$. If $u \in Q_2$, we get by construction of $\Gamma_j^{\rm III}$ that
$\dist(u,\Gamma^{\rm III}_j) \ge 2\vartheta^2 t_{i_j + n-1}$
(cf. Step 2 for a similar argument) and therefore by \eqref{eq: ecco}
$$\dist(u,\Gamma^{\rm III}_j) \ge c\vartheta^2 (s + t_{i_j+1})$$
for  $c>0$ small enough. In this case  \eqref{eq: dist3} holds.  On the other hand, if $u \in Q_1$,  we find $u' \in |[u;v_{i_j +n}]| \cap \partial P$ such that  $|[u';v_{i_j +n}]|$ separates $x$ and $p$. If we had $\dist(\Gamma^{\rm III}_j(s), u) <\vartheta^2(  t_{i_j+1} + s)$, we would get  by \eqref{eq: pi2}
\begin{align*}
|[u';v_{i_j +n}]| & \le \dist(u,\Gamma^{\rm III}_j(s)) + s_n^+ - s^+_{n-1} < \vartheta^2 (  t_{i_j+1}   +  s) + s_n^+ - s^-_{n} + s_n^- - s^+_{n-1}\\ & \le \vartheta^2 (  t_{i_j+1}   +  s_n^+) + 2\pi\vartheta^2 t_{i_j + n} + 2(t_{i_j +n } - t_{i_j  + n-1}). 
\end{align*}
Then the fact  that \eqref{eq: seg len} does not hold  and \eqref{eq: ecco} yield for $\vartheta$ small   with respect to $\omega$ and $C$
\begin{align*}
|[u;v_{i_j +n}]|  \le \vartheta^2 (  t_{i_j+1}   +  s_n^+) + (2\pi + 8)\vartheta^2 t_{i_j + n} \le  C\vartheta^2 t_{i_j + n}  <  \min\lbrace \vartheta t_{i_j + n}, \vartheta \omega d(P)\rbrace.  
\end{align*}
This contradicts  \eqref{eq: semico-pro}   and concludes the proof of \eqref{eq: dist3}.\smallskip\\

\noindent  \smallskip
\textit{Step 5: A curve between  $w^-_{n-1}$ and $p$}

\noindent  It remains to define a path between $w^-_{n-1}$ and $p$ (cf. below \eqref{eq:4,5}). Define a path $\Gamma^{\rm I_1}_m$ between $w^-_{n-1}$, $w^+_{n-1}$ as in Step 2 satisfying \eqref{eq: (I)}. Moreover, take a  path $\Gamma^{\rm II}_m$ between $w^+_{n-1}$, $w^-_{n}$ as in Step 3 such that \eqref{eq: length gamma} and \eqref{eq: dist2} hold. Let $\Gamma^{\rm I_2}_m = \Gamma^{\rm III}_m = \emptyset$ and let $\Gamma^{\rm IV}$ be the segment between $w_n^- = p + 2\vartheta^2 t_n \nu^+_{n-1}$ and $p$. Clearly, since  $B(p,\omega d(P)) \subset P$, we have for $\vartheta$ small with respect to $\omega$ 
\begin{align}\label{eq: length4}
\begin{split}
(i) \ \ l(\Gamma^{\rm IV}) = 2\vartheta^2 t_n , \ \ \  (ii)  \ \ \dist(\Gamma^{\rm IV}(s),\partial P) \ge \vartheta^2 d(P) \ge \vartheta^2 t_n = \vartheta^2 l(\gamma_0).
\end{split}
\end{align}

\noindent  \smallskip
\textit{Step 6: The curve $\gamma$ and the carrot condition}

\noindent  Now define $\gamma:[0,l(\gamma)] \to P$ such that $\gamma$ is parametrized by arc length and
$$\gamma([0,l(\gamma)]) = \bigcup^{m}_{j=0} \Big(\Gamma^{\rm I_1}_j \cup \Gamma^{\rm II}_j \cup \Gamma^{\rm I_2}_j \cup \Gamma^{\rm III}_j \Big) \cup \Gamma^{\rm IV}$$
with $\gamma(0) = x$ and $\gamma(l(\gamma)) = p$. We now show that \eqref{eq: carot} holds for $\vartheta$ sufficiently small. We have to derive that  $B(\gamma(\tau),\vartheta^3 \tau) \subset P$ for all $\tau \in [0,l(\gamma)]$. Let $\gamma(\tau) \in \hat{\Gamma}$, where $\hat{\Gamma} \in \lbrace \Gamma^{\rm IV}\rbrace \cup \lbrace \Gamma^{\rm X}_j, {\rm X}= {\rm I}_1,{\rm I}_2, {\rm II}, {\rm III}, j=0,\ldots,m\rbrace$. Choose $\tau_0 \le \tau$ such that $\gamma(\tau_0) = \hat{\Gamma}(0)$ and $i \in \lbrace 0,\ldots,n \rbrace$ such that $\gamma(\tau_0) \in \lbrace w^-_{i}, w^+_i\rbrace$. (Note that $i=i_j$ or $i={i_j}+1$ for some $j=0,\ldots,m$, cf. \eqref{eq:4,5}). Combining \eqref{eq: (I)}(i), \eqref{eq: length gamma}, \eqref{eq: length gamma2} and  \eqref{eq: length4}(i)    we derive by a telescope sum argument
\begin{align}\label{eq: carC}
\tau_0 \le \hat{C}t_{i} 
\end{align}
for some universal $\hat{C}\ge 1$ large enough, i.e.  $\gamma$ is at most $\hat{C}$ times longer than the original curve $\gamma_0$.  Letting $s = \tau - \tau_0$ and using
$$ \dist(\hat{\Gamma}(s),\partial P) \ge c\vartheta^2(t_i + s)$$
by \eqref{eq: (I)}(ii),  \eqref{eq: dist2}, \eqref{eq: dist3}, \eqref{eq: length4}(ii), respectively, we conclude by \eqref{eq: carC} for $\vartheta$ small
\begin{align*}
\dist(\gamma(\tau),\partial P) &= \dist(\hat{\Gamma}(s),\partial P) \ge c\vartheta^2(t_i + s) \ge c\vartheta^2(t_i + \hat{C}^{-1} s)  \\ 
&= c\hat{C}^{-1}\vartheta^2 \tau + c\vartheta^2 \big(t_i - \hat{C}^{-1} \tau_0\big)  \ge c\hat{C}^{-1}\vartheta^2 \tau\ge \vartheta^3 \tau. \quad\quad \quad \quad \quad \quad \quad \quad   \quad \quad  \text{\eop}
\end{align*}

\section{Proof of the main result and application}\label{se: mainproof}

This section is devoted to the proof of   Theorem \ref{th: main result smooth} and an application.  First, we prove the main partition result for polygons, which with the preparations in the last sections is now straightforward. 

\smallskip
\noindent\emph{Proof of Theorem \ref{th: main part}.} Let $\theta,\eps >0$. By Theorem \ref{theorem: partpol} and Theorem \ref{lemma: partpol*}  we first partition $P = P'_1 \cup \ldots P'_{m}$ into $\bar{\vartheta}$-semiconvex polygons with $\bar{\vartheta} = \bar{\vartheta}(\theta)$ such that $\sum^{m}_{j=1} \mathcal{H}^1(\partial P'_j) \le (1+C\theta)\mathcal{H}^1(\partial P)$  for $C>0$ universal.  Applying Theorem \ref{th: qussc} on each $P'_j$ with $\epsilon = \frac{1}{m}\eps$ and $\vartheta = \bar{\vartheta}$  we find $\tilde{\vartheta}=\tilde{\vartheta}(\theta)$, $\omega=\omega(\theta)$ and for each $P'_j$ a partition $P'_j = P^j_0  \cup P_{i_j +1} \cup \ldots \cup P_{i_{j+1}}$  with $i_1 = 0$, $i_{m+1} = N$ such that $\sum^{m}_{j=1} \mathcal{H}^1(\partial P^j_0) \le \eps$ and the polygons $(P_i)_{i=1}^N$ are $\tilde{\vartheta}$-semiconvex and $\omega$-rotund  with
$$\sum^{N}_{i=1} \mathcal{H}^1(\partial P_i) = \sum^{m}_{j=1} \sum^{i_{j+1}}_{k=i_j+1} \mathcal{H}^1(\partial P_k)\le \sum^{m}_{j=1}(1+C\theta)\mathcal{H}^1(\partial P'_j) \le(1+C\theta)\mathcal{H}^1(\partial P).$$
Define $P_0 = \bigcup_{j=1}^{m} P_0^j$. Starting the proof with $\theta C^{-1}$ instead of $\theta$, we obtain \eqref{eq: addition-main}. The fact that the polygons $(P_i)_{i=1}^N$ are $\varrho$-John domains for $\varrho=\varrho(\theta)$  follows from Theorem \ref{th: John}. \eop

The reader more interested in applications of our main result may now skip   Section \ref{sec: proofi} and continue with Section \ref{sec: gen}.

\subsection{Proof of the main result}\label{sec: proofi}

 To derive the result for sets with $C^1$-boundary we will have to combine different John domains. We start with an adaption of Lemma \ref{lemma: Johna}. In the following ${\rm diam}(D)$  denotes the diameter of a set $D \subset \R^2$.

\begin{lemma}\label{lemma: johnni} 
Let $0 < \varrho,c'  < 1$. Then for some $\varrho'=(c',\varrho)>0$ the following holds:

(i) Let $D_1,D_2 \subset \R^2$ be simply connected $\varrho$-John domains with Lipschitz boundary and $D_1 \cap D_2 = \emptyset$ such that $\partial D_1 \cap \partial D_2$ contains a segment $S$ with 
$$\mathcal{H}^1(S) \ge c' \min\lbrace \diam(D_1),\diam(D_2) \rbrace.$$
Then $D={\rm int}(\overline{D_1} \cup \overline{D_2})$ is a   $\varrho'$-John domain.

(ii)  Let $P$ be a polygon and $\triangle$ a closed triangle with ${\rm int}(P) \cap {\rm int}(\triangle) = \emptyset$ such that ${\rm int}(P)$ is a $\varrho$-John domain and $\partial P \cap \partial \triangle$ contains the longest edge of $\triangle$.  Then $D={\rm int}(\triangle \cup P)$ is a $\varrho'$-John domain with $\varrho'|D|\le |P|$. 
 \end{lemma}

\Proof (i) After rotation and translation we suppose $S =[(-2d,0);(2d,0)]$ with $d^2 \ge c\min\lbrace |D_1|,|D_2|\rbrace$ for $c=c(c')$, where the inequality follows from the assumption and  the isodiametric inequality. For $\eta>0$ let $Q_\eta = (-\eta d,\eta d)^2$, $Q^1_\eta = (-\eta d,\eta d) \times (0,\eta d)$ and $Q^2_\eta = (-\eta d,\eta d) \times (-\eta d,0)$. We now show that there is $\eta>0$ only depending on $\varrho$  such that possibly after changing the roles of $D_1$ and $D_2$ we have 
$$Q^i_\eta \subset D_i  \ \ \ \text{for $i=1,2$}.$$
We show the claim for $i=1$. As $\partial D_1$ is Lipschitz, we get that $[-d,d]\times (0,\eps] \subset D_1$ for $\eps$ small enough. Let $\gamma$ be a John curve in $D_1$ connecting $(-d,\eps)$ and $(d,\eps)$ (cf. Remark \ref{rem: john}).  Since $D_1$ is a $\varrho$-John domain, we find $0 <\eta<1$ only depending on $\varrho$  such that $\gamma \cap Q^1_\eta = \emptyset$. As $[-d,d] \times (0,\eps] \subset D_1$ and $D_1$ is simply connected, we then derive $Q^1_\eta \subset D_1$ as desired.

Define $D_i' = D_i \cup Q_\eta$ and $D' = D_1' \cup D_2'$. Each $D_i'$ is a $\bar{\varrho}$-John domain for $\bar{\varrho}= \bar{\varrho}(\varrho)$   by Lemma \ref{lemma: Johna}(i) since $|Q_\eta| = 2|Q^i_\eta|  = 2|Q^i_\eta \cap D_i|$. Moreover, we find $\varrho' = \varrho'(c',\varrho)$ such that  $D'$ is a $\varrho'$-John domain by Lemma \ref{lemma: Johna}(i) as 
$$\min \lbrace |D_1'|, |D_2'| \rbrace \le 2 \min \lbrace |D_1|, |D_2| \rbrace \le 2c^{-1}d^2 \le C|Q_\eta|$$
for $C$ only  depending on $c,\eta$ and thus only depending on  $\varrho$, $c'$.  Finally, possibly passing to a smaller $\varrho'$ also $D$ is a $\varrho'$-John domain since $D \setminus D' \subset \partial D_1 \cap \partial D_2$.

(ii) Note that (i) is not directly applicable as ${\rm int}(\triangle)$ is possibly not a $\varrho$-John domain. Suppose $S= [(-d,0);(d,0)]$ is the longest edge of $\triangle$ and $\triangle \subset [-d,d] \times [0,\infty)$. Arguing as in (i), we find $\eta$ only depending on $\varrho$ such that the closed triangle $\triangle'$ with vertices $(-d,0)$, $(d,0)$ and $(0,-d\eta)$ is completely contained in $P$. As $\mathcal{H}^1(S) ={\rm diam(\triangle)}$, it is not hard to see that  
$B := {\rm int} (\triangle \cup \triangle')$ is a $\varrho'$-John domain with $\varrho'$ only depending on $\varrho$. Moreover, we find
\begin{align}\label{eq: mori}
|\triangle| \le C|\triangle'|
\end{align} 
for $C=C(\varrho)$.  Then by Lemma \ref{lemma: Johna}(i) and \eqref{eq: mori} also $D = {\rm int}(P \cup B) = {\rm int}(P \cup \triangle)$ is a John domain for a possibly smaller $\varrho'$ only depending on $\varrho$. Finally, $\varrho'|D|\le |P|$ follows from \eqref{eq: mori} for $\varrho'$ small enough.   \eop

Before we concern ourselves with sets with $C^1$-boundary we state the following corollary of Theorem \ref{th: main part}. 

\begin{corollary}\label{cor: main}
Let be given the situation of Theorem \ref{th: main part}. If $\sphericalangle(v,P) \ge \frac{\pi}{4}$ for all $v \in \mathcal{V}_P$, one can set $P_0 = \emptyset$.
\end{corollary}

\Proof First, we apply Theorem \ref{th: main part} to get a partition of $P= \bigcup_{j=0}^N P_j$. Recall that by the construction in the proof of Theorem \ref{th: qussc} and Theorem \ref{th: main part} the component $P_0$ is the finite union of closed,  isosceles triangles with exactly one interior angle smaller than $\frac{\pi}{4}$ (see before \eqref{eq: addition*}). We first see that each two triangles $\triangle_1$, $\triangle_2$ do not share a segment. Indeed, otherwise the  corresponding convex polygons, denoted by $P'_1$, $P'_2$, from which the triangles are cut out, share a segment and contain $v \in \mathcal{V}_{P_1'} \cap \mathcal{V}_{P_2'}$ with $v \in \triangle_1 \cap \triangle_2$ and $\sphericalangle(v,P'_i) \le \frac{\pi}{4}$ for $i=1,2$. As the partition can be constructed such that endpoints of introduced segments never coincide unless there are concave vertices of $P$ (see Remark \ref{rem: vertices1}(ii)), we derive $v \in \mathcal{V}_P'$. Then, however, Remark \ref{rem:  vertices1}(iii)  implies $\mathcal{H}^1(\partial P_1' \cap \partial P_2')=0$, which gives a contradiction.

Moreover,  it is not restrictive to assume that each edge of a triangle $\triangle$ is completely contained in  $\partial P$ or some $\partial P_j$ since otherwise we choose an isosceles $\triangle' \subset \triangle$ with the desired property. We then note that the $\triangle \setminus \triangle'$ is a convex polygon with interior angles larger than $\frac{\pi}{4}$ and thus  we can apply Lemma \ref{lemma: convpart} to obtain a refined partition of $\triangle \setminus \triangle'$ consisting of $\varrho$-John domains such that \eqref{eq: addition-main} still holds.

The assumption $\sphericalangle(v,P) \ge \frac{\pi}{4}$ for all $v \in \mathcal{V}_P$ implies that for each $\triangle$ at least one of the two longer edges is contained in some $\partial P_i$. Then ${\rm int}(P_i \cup \triangle)$ is a John domain by Lemma \ref{lemma: johnni}(ii)  for a John constant only depending on $\theta$ and $|{\rm int}(P_i \cup \triangle)| \le C|P_i|$ for $C=C(\theta)$.  Hereby we can define a partition $(P'_i)_i$ of $\Omega$ satisfying \eqref{eq: addition-main} such that each component $P_i'$ is the union of $P_i$ with some triangles adjacent to $P_i$. Now Lemma \ref{lemma: Johna}(ii) (with $D_0 = P_i$,  $D_j = {\rm int}(P_i \cup \triangle)$ for $j \ge 1$) yields that all $P_j'$ are John domains for a John constant only depending on $\theta$.   \eop

We now extend the result to sets with smooth boundary, where we first derive a version without the sharp estimate \eqref{eq: addition-mainXXXXX}.

\begin{theorem}\label{th: main result smoothX}
Theorem \ref{th: main result smooth} holds with $\sum^N_{j=1}\mathcal{H}^1(\partial \Omega_j) \le C\mathcal{H}^1(\partial \Omega)$ in place of \eqref{eq: addition-mainXXXXX} for a universal $C>1$.
\end{theorem}

\Proof As $\Omega$ has $C^1$-boundary and $\partial \Omega$ is connected  due to the fact that $\Omega$ is simply connected, we can find $p_0,\ldots,p_{n-1} \in \partial \Omega$ such that the closed squares $Q_i$ with diagonal $[p_i;p_{i+1}]$ for $i=0,1,\ldots,n-1$ (set $p_{n} = p_0$ and $Q_n = Q_0$) satisfy 
\begin{align}\label{eq: last cond}
(i) & \ \ d:= \min_{i=0,\ldots,n-1}|[p_i;p_{i+1}]| \ge \frac{1}{2}\max_{i=0,\ldots,n-1}|[p_i;p_{i+1}]|,\\
(ii) & \ \ Q_i \cap Q_{i+1} = \lbrace p_{i+1}\rbrace, \,  \ \dist(Q_i, Q_{(i+k)\,{\rm mod}n}) \ge \frac{d}{2} \ \,  \ \text{ for }  i=0,\ldots,n-1,  \,  |k|\ge 2 \notag
\end{align}
and $\partial \Omega \cap Q_i$ is the graph of a $C^1$ function, where the angle enclosed by $p_{i+1}-p_i$ and the tangent vector of $\partial \Omega$ in $\partial \Omega \cap Q_i$ is smaller than $\frac{\pi}{8}$. Moreover, this can be done in the way that all  interior angles of the interior polygon $P_{\rm int} := \overline{\Omega \setminus \bigcup_{i=0}^{n-1}Q_i}$ are larger than $\frac{\pi}{4}$. Define also the sets 
\begin{align}\label{eq: Pi*} 
 P^{\rm out}_i = \Omega \cap {\rm int}(Q_i)
 \end{align}
 for $i=0,\ldots,n-1$.  The geometry of  $P^{\rm out}_i$ implies that $P^{\rm out}_i$ has Lipschitz boundary and is a $c$-John domain for a universal constant $c>0$. Moreover, we observe that $\mathcal{H}^1(\partial P_{\rm int}) + \sum_{i=0}^{n-1} \mathcal{H}^1(\partial P^{\rm out}_i)\le C\mathcal{H}^1(\partial \Omega)$. The claim   follows from Corollary \ref{cor: main} applied on $P_{\rm int}$.  \eop

This together with Lemma \ref{lemma: johnni}   allows  to give the proof of Theorem \ref{th: main result smooth}.

\smallskip
\noindent\emph{Proof of Theorem \ref{th: main result smooth}.} By Corollary \ref{cor: main} we find a partition  $P_{\rm int}=  P_1 \cup \ldots \cup P_N$ of the polygon $P_{\rm int}$ constructed in the proof of Theorem \ref{th: main result smoothX}, where  by \eqref{eq: addition-main}  for $C>0$ universal
\begin{align}\label{eq: other bound}
\sum_{j=1}^N\mathcal{H}^1(\partial P_j\setminus \partial P_{\rm int}) \le \theta\mathcal{H}^1(\partial P_{\rm int}) \le C\theta \mathcal{H}^1(\partial \Omega).
\end{align}
The goal is now to combine each $P_j$ with certain $(P^{\rm out}_i)_{i=0}^{n-1}$ defined in \eqref{eq: Pi*}  such that the resulting sets are still John domains and \eqref{eq: addition-mainXXXXX} holds. Let $J$ be the set of indices such that $j \in J$ if and only if $\diam(P_j) < \frac{d}{4}$ with $d$ as in \eqref{eq: last cond}.  By \eqref{eq: last cond} we see that each $P_j$, $j \in J$, intersects at most two sets $\overline{P_i^{\rm out}}$, $i=0,\ldots,n-1$. Recalling the geometry of  $(P^{\rm out}_i)_i$ and the fact that the interior angles of the polygon $P_{\rm int}$ are larger than $\frac{\pi}{4}$, we find $\mathcal{H}^1(\partial P_j) \le C\mathcal{H}^1(\partial P_j \setminus \partial P_{\rm int})$ for $j \in J$ for a universal constant $C>0$ and thus by \eqref{eq: other bound}
\begin{align}\label{eq: other bound*}
\sum_{j \in J} \mathcal{H}^1(\partial P_j) \le C\sum_{j \in J} \mathcal{H}^1(\partial P_j \setminus \partial P_{\rm int}) \le \theta\mathcal{H}^1(\partial P_{\rm int}) \le C\theta\mathcal{H}^1(\partial \Omega).
\end{align}
Recall the definition of $Q_i$ in \eqref{eq: last cond} and denote by $Q_i'$ the enlarged square with the same center and orientation, but with diagonal length $\frac{5}{4}|[p_i;p_{i+1}]|$. Note that all sets $P_{\rm int} \cap Q_i'$ are Lipschitz and are all related to a square of sidelength $d$ through Lipschitz homeomorphism with Lipschitz constants of  both the homeomorhism itself and its inverse uniformly bounded independently of $i$. Let $\bar{c}>0$ to be specified below in \eqref{eq: porg}-\eqref{eq:impliesX}. We observe that there is $\bar{C}=\bar{C}(\bar{c})>0$ such that
\begin{align}\label{eq: other bound-new}
\# I \le \bar{C}\theta d^{-1}\mathcal{H}^1(\partial \Omega),   \text{ where }  I := \Big\{i: \mathcal{H}^1\Big(  {\rm int}(P_{\rm int} \cap Q_i') \cap \bigcup\nolimits_{j=1}^N \partial P_j \Big) \ge \bar{c}d \Big\}.
\end{align}
Indeed, this follows from \eqref{eq: other bound} and \eqref{eq: last cond}. 

Consider $i \notin I$. For $j=1,\ldots,N$ define the components $A_{j,i} := P_j \cap Q'_i$ of $P_{\rm int} \cap Q_i'$ and denote by $(A_{j,i}^k)_k$ the connected components of $A_{j,i}$. Then the result in \cite[Lemma 4.6]{FriedrichSolombrino}, which essentially relies on the relative isoperimetric inequality, shows that for $\bar{c}$ sufficiently small there is exactly one component $B_i := P_{j_i} \cap Q_i' \subset (A_{j,i})_{j=1}^N$ with $|B_i| >\frac{1}{2}|P_{\rm int} \cap Q_i'|$ and the other components $A_{j,i} \neq B_i$  satisfy
\begin{align}\label{eq-rev}
{\rm diam}(A_{j,i}^k) \le C\mathcal{H}^1(\partial A^k_{j,i} \cap {\rm int}(P_{\rm int} \cap Q_i')) \le   C\bar{c} d \  \ \ \text{for all} \ \ \ A_{j,i}^k
\end{align}
for a universal $C>0$, particularly independent of $i$ and $A_{j,i}^k$. Then using the fact that $\mathcal{H}^1( \partial (P_{\rm int} \cap Q_i') \cap \partial A_{j,i}^k) \le C{\rm diam}(A_{j,i}^k)$ by the geometry of $P_{\rm int} \cap Q_i'$,  we get  by \eqref{eq: other bound-new}-\eqref{eq-rev}
 \begin{align}\label{eq: porg}
\sum_{A_{j,i} \neq B_i}\mathcal{H}^1(\partial A_{j,i}) \le C\sum_{A_{j,i} \neq B_i}\mathcal{H}^1(\partial A_{j,i} \cap {\rm int}(P_{\rm int} \cap Q_i')) + C \sum_{A_{j,i} \neq B_i}\sum_k {\rm diam}(A_{j,i}^k) \le   C\bar{c} d.
\end{align}
Moreover, as $\dist(\partial Q_i', \partial Q_i) \ge \frac{1}{8\sqrt{2}}d$, for $\bar{c}$ small enough we derive  by \eqref{eq-rev}
\begin{align}\label{eq:implies}
A_{j,i} = (P_j \cap Q'_i) \neq B_i \text{ and } \partial A_{j,i} \cap \partial P^{\rm out}_i \neq \emptyset \ \ \  \Rightarrow \ \ \ P_j \subset Q_i' \text{ and } j \in J.
\end{align}
Moreover, we find  
\begin{align}\label{eq:impliesX}
(i) \ \ \mathcal{H}^1(\partial B_i \cap \partial P_i^{\rm out}) \ge \frac{d}{2}, \ \ \ \ \ (ii) \ \  \partial B_i \cap \partial P_i^{\rm out} \ \text{ connected}.
\end{align}
First, (i) follows for $\bar{c}$ small from \eqref{eq: porg}  and the fact that $\mathcal{H}^1(\partial P_i^{\rm out} \cap \partial P_{\rm int}) \ge d$. If (ii) was false, we would  find that the polygon $B_i = P_{i_j} \cap Q_i'$  has at least one concave vertex not lying on $\partial P_{\rm int}$. This, however, contradicts the construction of the partition, cf. Remark \ref{rem: vertices1}(i) and Remark \ref{rem: vertices2}.  Note that \eqref{eq:impliesX}(i) implies $j_i \notin J$. By Lemma \ref{lemma: johnni}(i),  \eqref{eq: last cond} and \eqref{eq:impliesX} we find 
\begin{align}\label{eq: johnni1}
D_i := {\rm int} (P_{j_i} \cup \overline{P_i^{\rm out}}) 
\end{align}
is a $\varrho'$-John domain with Lipschitz boundary for $\varrho' = \varrho'(\theta)$. 

We are now in the position to define the partition of $\Omega$. For all $j \notin J$, let $I_j \subset \lbrace 0,\ldots,n-1\rbrace \setminus I$ be the index set such that $P_j \cap Q_i' = B_i$ if and only if $i \in I_j$, where $B_i = P_{j_i} \cap Q_i'$ as above. Note that the above arguments in \eqref{eq:implies} show that the sets $(I_j)_j$ are pairwise disjoint and also observe that $I_j$ may be empty. Define $P'_j = \bigcup_{i \in I_j} D_i $ for $j \notin J$ with $D_i$ as in \eqref{eq: johnni1} and consider the partition $(\Omega_j)_j$ consisting of the sets
$$ (P'_j)_{j \notin J} \ \cup \ ({\rm int}(P_j))_{j \in J} \ \cup \ (P_i^{\rm out})_{i \in I}.$$ 
Note that the sets cover $\Omega$ up to a set of negligible measure since each $P_i^{\rm out}$, $i \notin I$, is contained in some $P'_j, j \notin J$. Note that by   \eqref{eq:implies} we derive
$$\bigcup\nolimits_{j} (\partial \Omega_j \setminus \partial \Omega) \subset \bigcup\nolimits_{i \in I} (\partial P_i^{\rm out} \cap \partial P_{\rm int}) \cup \bigcup\nolimits_{j=1}^N (\partial P_j \setminus \partial P_{\rm int}) \cup \bigcup\nolimits_{j \in J} (\partial P_j \cap \partial P_{\rm int}).  $$
This together with \eqref{eq: other bound},  \eqref{eq: other bound*}  and $\sum_{i \in I} \mathcal{H}^1(\partial P_i^{\rm out} \cap \partial P_{\rm int}) \le Cd \ \#I \le C\theta \mathcal{H}^1(\partial \Omega)$ (see \eqref{eq: other bound-new}) yields   $\sum_{j}\mathcal{H}^1(\partial \Omega_j \setminus \partial \Omega) \le C \theta\mathcal{H}^1(\partial \Omega)$ and herefrom we indeed derive \eqref{eq: addition-mainXXXXX} since we can replace $\theta$ by $C^{-1}\theta$ in the above proof. Finally, observe that all components are John domains with Lipschitz boundary for a John constant only depending on $\theta$, where for the sets $(P'_j)_{j \notin J}$ we use  \eqref{eq: johnni1} and Lemma \ref{lemma: Johna}(ii). \eop

\subsection{A generalization and application}\label{sec: gen}
We now present a generalized version of Theorem \ref{th: main result smooth} for Lipschitz sets which are not necessarily simply connected. This version will be one of the main ingredients of \cite{Friedrich:15-4}.  For a bounded set $D\subset \R^2$ we introduce the \emph{saturation} of $D$ defined by $\sat(D) = {\rm int}(\R^2\setminus E^0)$, where $E^0$ denotes the unique unbounded connected component of $\R^2 \setminus D$.

 \begin{theorem}\label{th: main part2}
Let $\eps > 0$ and $M \in \N$. Then there is a universal constant $\varrho>0$ and $C=C(M)>0$ such that for all bounded domains  $\Omega \subset \R^2$ with Lipschitz boundary and the property that $\sat(\Omega) \setminus \Omega$ consists of at most $M$ components the following holds: There is a  partition  $\Omega = \Omega_0 \cup \ldots \cup \Omega_N$ such that $|\Omega_0|   \le \eps$ and the sets $\Omega_1,\ldots,\Omega_{N}$ are $\varrho$-John domains with Lipschitz boundary with
$$\sum\nolimits^N_{j=0}\mathcal{H}^1(\partial \Omega_j) \le C\mathcal{H}^1(\partial \Omega).$$
\end{theorem}

\Proof Let $\eps>0$ be given and let $U_1,\ldots,U_m$ be the connected components of $\sat(\Omega) \setminus \Omega$ with $m \le M$. For each $U_j$  we can choose a segment $S_j$ with $\mathcal{H}^1(S_j) \le \diam(\Omega) \le \mathcal{H}^1(\partial \Omega)$ such that $\Theta_j:= \partial U_j \cup S_j \cup \partial (\sat(\Omega))$ is connected. Consequently, we get that each connected component  of $\Omega \setminus \bigcup_{j=1}^m \Theta_j$ is simply connected. For $s >0$ we cover $\R^2$    with squares of the form $Q(p) = p + [-s,s]^2$, $p \in 2s\Z^2$. Let 
$$\mathcal{Q}_s : = \big\{Q(p): Q(p) \cap \Omega \neq \emptyset, \  Q(p) \cap \big(\partial \Omega \cup \bigcup\nolimits_{j=1}^m S_j\big) \neq \emptyset \big\}.$$
Since  $\Omega$ has Lipschitz boundary, we find that for $s$ sufficiently small 
\begin{align}\label{eq:num}
s\# \mathcal{Q}_s \le C  \mathcal{H}^1\big(\partial \Omega \cup \bigcup\nolimits_{j=1}^m S_j\big)\le C\mathcal{H}^1(\partial \Omega) + CM{\rm diam}(\Omega) \le C\mathcal{H}^1(\partial \Omega)
\end{align} 
with $C=C(M)$. 
By $(P_i)_i$ we denote the connected components of $\R^2 \setminus \bigcup_{Q(p) \in \mathcal{Q}_s} {Q(p)}$  having nonempty intersection with $\Omega$. Since each $P_i$ is the union of squares and the connected components  of $\Omega \setminus \bigcup_j \Theta_j$ are simply connected, also $P_i$ is simply connected and thus $\overline{P_i}$ is a polygon with interior angles not smaller than $\frac{\pi}{2}$. Moreover, we find by \eqref{eq:num}
$$\sum\nolimits_i \mathcal{H}^1(\partial P_i) \le 8s \# \mathcal{Q}_s  \le  C\mathcal{H}^1(\partial \Omega).$$
Likewise, if we choose $s$ small enough, we get that $\Omega_0 : = \Omega \setminus \bigcup_i P_i$ satisfies  
$$|\Omega_0| \le 4s^2 \# \mathcal{Q}_s \le Cs^2\mathcal{H}^1(\partial \Omega) \le \eps, \ \ \ \ \  \mathcal{H}^1(\partial \Omega_0) \le C\mathcal{H}^1(\partial \Omega).$$
The result now follows from Corollary \ref{cor: main} applied on each $P_i$ for $\theta=1$. (Note that alternatively one may also apply Theorem \ref{th: main part} on each $P_i$ choosing the occurring exceptional sets $P^i_0$ small enough in terms of $\eps$.)
 \eop

Finally, we derive a piecewise Korn inequality for a certain subclass of $SBD$ (we refer to \cite{Ambrosio-Coscia-Dal Maso:1997,Bellettini-Coscia-DalMaso:98} for more details on this function space). Although this problem will be thoroughly discussed in \cite{Friedrich:15-4}, we include a simplified analysis in the present exposition to give a first   application of the main results of this article. 

Let $1<p<\infty$ and $M \in \N$. For an open, bounded set $\Omega \subset \R^2$ with Lipschitz boundary we let $\mathcal{W}^p_M(\Omega)$ be the set of functions in $SBD^p(\Omega)$  whose jump set $J_{y} = \bigcup^m_{j=1} \Gamma^y_j$  is the finite union of closed connected pieces of Lipschitz curves with at most $M$ components (i.e. $m \le M$)  and $y|_{\Omega\setminus J_y} \in W^{1,p}(\Omega \setminus J_y)$. Note that similar assumptions have been used, e.g., in \cite{Chambolle:2003, DM-Toa, Lazzaroni, NegriToader:2013}.

\begin{theorem}\label{th: main korn}
Let  $p \in (1,\infty)$ and $M \in \N$.  Then there is   $c=c(p)>0$  and $C=C(M)>0$ such that for all $\Omega \subset \R^2$ open, bounded with Lipschitz boundary with the property that $\sat(\Omega) \setminus \Omega$ consists of at most $M$ components the following holds: For  each $y \in \mathcal{W}_M^p(\Omega)$ there is a partition $(\Omega_j)^N_{j=0}$ of  $\Omega$ with 
 \begin{align}\label{eq: main1}
 \sum\nolimits^N_{j=0} \mathcal{H}^1( \partial \Omega_j) \le C(\mathcal{H}^1(J_y) + \mathcal{H}^1(\partial \Omega))
 \end{align}
and corresponding $A_j \in \R^{2\times 2}_{\rm skew}$, $b_j \in \R^2$ such that 
 $u:= y - \sum\nolimits_{j=0}^N (A_j \, \cdot + b_j)  \chi_{\Omega_j}$ satisfies 
$$
(\diam(\Omega))^{-1}\Vert u \Vert_{L^p(\Omega)} +  \Vert \nabla u \Vert_{L^p(\Omega)} \le c \Vert e(y) \Vert_{L^p(\Omega)} ,
$$
where $e(y) = \frac{1}{2}(\nabla y^T + \nabla y)$.  
\end{theorem}
 
  Note that  one essential point is that the constant $c$ does not depend on $\Omega$ and $C$ depends on $\Omega$ only in terms of the number of components of $\sat(\Omega) \setminus \Omega$. Therefore, the result is also interesting in the case of varying domains $\Omega$ and functions $y \in W^{1,p}(\Omega)$. 
 
 \Proof A classical result states that $y$ is piecewise rigid if $\Vert e(y)\Vert_{L^p(\Omega)} = 0$ (see also \cite{Chambolle-Giacomini-Ponsiglione:2007}), so we can concentrate on the case $\Vert e(y) \Vert_{L^p(\Omega)}>0$.  Applying the following results on each connected component of $\Omega$ separately, it is not restrictive to assume that $\Omega$ is connected. Moreover, we may suppose that  $\Omega$ is simply connected as otherwise we consider $\sat(\Omega)$ and define an extension $\bar{y}$ with  $\bar{y}=0$ on $\sat(\Omega) \setminus \Omega$, where we obtain $\bar{y} \in \mathcal{W}_{2M}^p(\sat(\Omega))$. 
 
 We now repeat the arguments in   the proof of Theorem \ref{th: main part2} on $(\Gamma^y_j)_j$ instead of $(U_j)_j$: we introduce segments to obtain simply connected components of $\sat(\Omega)$ and covering the boundary with squares we obtain an estimate of the form \eqref{eq:num}, where the right hand side now also depends on $\mathcal{H}^1(J_y)$. As before this yields a partition $(\Omega_j)_{j=0}^N$ of  $\Omega$  such that $|\Omega_0| \le \eps$ for an arbitrarily small $\eps>0$ and $\Omega_1,\ldots,\Omega_N$ are $\varrho$-John domains for a universal constant $\varrho$. Then \eqref{eq: main1} follows as in Theorem \ref{th: main part2}.

  As Korn's inequality holds on John domains with a constant only depending on the John constant (see e.g. \cite{Acosta}), we get by an elementary scaling argument
 $$\sum\nolimits^N_{j=1} \big( (\diam(\Omega_j))^{-p}\Vert y - (A_j\,\cdot + b_j)\Vert^p_{L^p(\Omega_j)} + \Vert \nabla y - A_j\Vert^p_{L^p(\Omega_j)} \big)\le   c\Vert e(y)\Vert^p_{L^p(\Omega)}$$
 for suitable $A_j \in \R^{2 \times 2}_{\rm skew}$, $b_j\in \R^2$ and $c=c(p)$. Finally, as $y \in L^{p}(\Omega)$, $\nabla y \in L^{p}(\Omega)$ and $|\Omega_0|\le\eps$, we find $(\diam(\Omega))^{-1}\Vert y\Vert_{L^p(\Omega_0)} + \Vert \nabla y\Vert_{L^p(\Omega_0)} \le \Vert e(y) \Vert_{L^p(\Omega)}$ for $\eps$ small enough so that the assertion holds for $u = y - \sum\nolimits^N_{j=1} (A_j \, \cdot + b_j)  \chi_{\Omega_j}$. \eop

 \smallskip

\noindent \textbf{Acknowledgements} This work has been funded by the Vienna Science and Technology Fund (WWTF)
through Project MA14-009.   The support by the Alexander von Humboldt Stiftung is gratefully acknowledged. I am gratefully indebted to the referee for her/his careful reading of the manuscript and helpful suggestions.


 \typeout{References}

\end{document}